\documentclass[envcountsame,envcountchap]{svmono}

\usepackage{amssymb}
\usepackage{amsmath, amsfonts, amssymb,  amscd}
\usepackage{enumerate}

\usepackage{makeidx}         
\usepackage{graphicx}        
\usepackage{multicol}        
\usepackage[bottom]{footmisc}

\makeindex             

\newtheorem{thm}{\bf Theorem}[chapter]
\newtheorem{prop}[thm]{\bf Proposition}
\newtheorem{defin}[thm]{\bf Definition}
\newtheorem{coro}[thm]{\bf Corollary}
\newtheorem{lem}[thm]{\bf Lemma}

\smartqed

\newcommand{\mmmintone}[1]{{\dis{\int\kern -.36cm-}}_{\kern-.21cm\substack{#1}}\;\;}
\newcommand{\mmmintwo}[2]{{\dis{\int\kern -.43cm-}}_{\kern-.21cm\substack{#1}}^{\substack{#2}}\;\;}
\newcommand{\submint}{{\scriptstyle{\int\kern -.66em -}}}
\newcommand{\submintone}[1]{{\scriptstyle{\int\kern -.66em-}}_{\scriptscriptstyle{\kern-.21em\substack{#1}}}}
\newcommand{\fracmint}{{\textstyle{\int\kern -.88em -}}}
\newcommand{\fracmintone}[1]{{\textstyle{\int\kern -.88em
-}}_{\scriptscriptstyle{\kern-.21em\substack{#1}}}\;}

\newcommand{\nn}{\nonumber}

\newcommand{\eps}{\epsilon}

\newcommand{\ga}{\gamma}
\newcommand{\Ga}{\Gamma}
\newcommand{\Om}{\Omega}

\newcommand{\si}{\sigma}

\newcommand{\und}{\underline}
\newcommand{\N}{\mathbb N}
\newcommand{\la}{\lambda}

\newcommand{\dis}{\displaystyle}

\begin{document}

\author{G. Carinci, A. De Masi, C. Giardin\`a, E. Presutti}
\title{Free boundary problems\\
in PDEs and particle systems}
%
%
\vskip3cm
\noindent
{\Large G. Carinci, A. De Masi, C. Giardin\`a, E. Presutti}
\vskip2cm
\noindent
{\Huge Free boundary problems}\\
{\Huge in PDEs and particle systems}
\thispagestyle{empty}

\tableofcontents

\mainmatter

\renewcommand{\theequation}{\thesection.\arabic{equation}}
\setcounter{equation}{0}

\chapter{Introduction}
\label{ch:00}

%
%
We develop here a theory for free boundary problems which applies to
a large class of systems arising from problems in various, even distant, areas of research and which share a  common mathematical structure.  As we shall see in some detail,
these are models for heat conduction, queuing theory, propagation of fire, interface dynamics, population dynamics and evolution of biological systems with selection mechanisms.  We shall consider models in  continuum and interacting
particle systems. Their common mathematical features are the following:

\vspace{0.3cm}
\framebox{
\begin{minipage}{.85\textwidth}
{
\medskip
(1) Microscopic particle dynamics stem from interactions of topological nature. 

\bigskip

(2) Macroscopic evolution is ruled by a free boundary problem.
\smallskip}
\end{minipage}
}

\vspace{0.3cm}
\noindent
In fact in the models we consider the particles move in $d=1$ dimension so that there is a rightmost and a leftmost particle, called {\em boundary particles}.  The rules of dynamics are the usual ones, particles are either free (independent random walks or Brownian motions) or they have some local interaction (for instance simple exclusion) and on top of that there may be creations of new particles or particles may duplicate via a branching process. 
In addition, in order to keep (approximatively) constant the total number of particles,
boundary particles are subject to a death process.

\medskip
The topological nature of the interaction refers to the fact that the boundary particles are special as they may disappear at some given rate, being then replaced by new boundary particles, the
rightmost and leftmost particles among those which have survived.  Thus the ``inside particles'', i.e.\ those in between the boundary particles, evolve in the ``usual'' way, but the inside particles are not fixed a priori and may eventually become boundary particles depending on the evolution itself.

As a consequence of  particle evolution, the spatial domain occupied by the particles varies in time. In particular the location of  the boundary particles changes in the course of time
due to the death process at the boundary.
Correspondingly, as we shall discuss extensively in this volume,  
the macroscopic version of the models is provided 
by a {\em free boundary problem} for a PDE with Dirichlet condition
supplemented by prescribing the boundary flux.
As often occurs, one can relate a macroscopic evolution
to  microscopic dynamics via a {scaling limit} procedure
(hydrodynamic limit).

\vspace{0.3cm}

The basic example that we will study in detail here is given by the linear heat equation \begin{equation*}
 \frac{\partial \rho}{\partial t} = \frac 12
 \frac{\partial^2 \rho}{\partial r^2}
   \end{equation*}
in the time varying domain $[0,X_t]$ with some
initial condition $\rho(r,0)=\rho_0(r)\ge 0$ and boundary conditions
  \begin{equation*}
-\frac 12 \frac{\partial \rho}{\partial r}(0,t) = j >0,\quad
  \rho(X_t,t)= 0 \;. 
   \end{equation*}
The free boundary $X_t$ (also called the {\em edge} in this book) 
is not given a priori but it
should be determined in such a way that
 \begin{equation*}
-\frac 12 \frac{\partial \rho}{\partial r}(X_t,t) = j.
   \end{equation*}
Interpreting $\rho$ as a mass density, the last condition states that the mass flux leaving the system at $X_t$ must be equal to $j$, and since $j$ is also the mass flux entering at 0 (as fixed by the boundary condition at 0), the total mass in the system is preserved.  From this perspective the free boundary problem becomes a control problem: find an edge evolution $X_t$ in such a way that the total mass is constant in time.

%
%
%

\vspace{0.3cm}

\noindent
Well known theorems on the Stefan problem yield a local existence theorem for our basic example when we have 
``classical initial data''. We will define here a weak version of the problem and prove global existence and uniqueness of a {\em relaxed solution} for  general initial data.
 The other models that we will consider in this work have similar structure and the strategies of proof are very close to that in the basic example.  The key point in all of them is:

\vspace{0.3cm}
\framebox{
\begin{minipage}{.85\textwidth}
{
\smallskip
Construct upper and lower barriers that squeeze the solution we are looking for.
\smallskip}
\end{minipage}
}

\vspace{0.3cm}

\noindent
The correct notion of order for these problems is
defined by mass transport.
Referring to a basic example for the sake of definiteness, the barriers are defined
in terms of a simplified evolution where we introduce a time grid of length $\delta$
and the evolution is ruled by the heat equation in $\mathbb R_+$ in the open intervals
$(n\delta,(n+1)\delta)$ with boundary condition
 \begin{equation*}
-\frac 12 \frac{\partial \rho}{\partial r}(0,t) = j\;.
   \end{equation*}
At the times $n\delta$ we remove an amount of mass equal to $j\delta$ so that at these times the mass conservation is restored.  The key point is that we get an upper barrier if we start by removing mass already at time 0 while we get a lower barrier when we   remove mass from   time $\delta$ on:
the order here is in the sense of moving mass to the right.
A key step is to prove that  the barriers have a unique separating element. Once we have this we conclude by showing that the 
solution we are considering is trapped in between the barriers which then identifies the solution as the element separating the barriers.
As we will see this part of the proof exploits extensively probabilistic ideas and techniques based on the well known relation between heat equation and Brownian motion and between the hitting distribution at the boundaries and the Dirichlet condition in the heat equation.

We think it can be useful for the reader to have one case worked out in all details, so that in Part I we prove the above in the context of our basic example by proving global existence and uniqueness of the relaxed solution of the problem; we also show that this is the limit of the empirical mass density of the associated particle system (in the hydrodynamic limit).  In Part II we discuss, in a very sketchy way, several other models, the conjecture being
that the results proved for the basic model extend to these other cases that have been done, at least partially.

%
%
\part{The basic model}

\renewcommand{\theequation}{\thesection.\arabic{equation}}
\setcounter{equation}{0}

\chapter{{\bf Introduction to Part I}}
\label{ch:1}

In Part I of this work we study a  model for mass transport where
 Fick's law is satisfied.   Fick's law is
the analogue for mass of  Fourier's law for heat conduction.
Fourier's law, see  \cite{fourier},  specifies the amount of heat flux in
a metal bar when we heat it
from one side and cool  it from the other.
Its analogue  for mass fluxes
is   Fick's law,
formally  described
by the same equation.
Since the transversal direction to the flow is not
relevant we model our system as one dimensional.
The ideal  experiment of mass transport that we have in mind
is the following: for $t\ge 0$ we confine the system in a time varying space interval $[0,X_t]$,
where $X_t$ is a given positive, continuous and piecewise $C^1$ function; 
for instance we move the edge
 $X_t$ with constant velocity for some time, then we change velocity and so on.
We act on the system by injecting mass
from its left boundary 0 at rate $j>0$ while we remove mass from
the right boundary $X_t$ in such a way as to keep  the mass
density at $X_t$ equal to 0 for all $t\ge 0$.  The evolution of the mass density
$\rho(r,t)$
in the interior of the spatial domain is ruled by
combining the continuity equation and   Fick's law, so that,
supposing  a constant conductivity (set  equal to 1/2), we have
  \begin{equation}
    \label{intro.1}
 \frac{\partial \rho}{\partial t} = -
 \frac{\partial J}{\partial r}\;,
 \qquad
 J = - \frac12 \frac{\partial \rho}{\partial r}
   \end{equation}
where $J(r,t)$ is the local mass-flux and $\rho(r,t)$ the mass density.  Thus
$\rho(r,t)$  solves the heat equation
  \begin{equation}
    \label{intro.2}
 \frac{\partial \rho}{\partial t} = \frac 12
 \frac{\partial^2 \rho}{\partial r^2}
   \end{equation}
in the time varying domain $[0,X_t]$ with some
initial condition $\rho(r,0)=\rho_0(r)$ and boundary conditions
  \begin{equation}
    \label{intro.3}
J(0,t) = j,\quad
  \rho(X_t,t)= 0 \;. 
   \end{equation}
Physically these boundary conditions mean that the system
is in contact with a {\em current reservoir} which  sends in mass at rate $j$ and thus imposes a current $j$ at the origin; instead at the other endpoint $X_t$ there is a
{\em density reservoir} which removes mass as fast as needed to fix the mass density to be
constantly
equal
to zero.
As a consequence, in this setting, the total mass of the system is not 
a conserved quantity.

\medskip
The main question we want to study here arises when we require 
mass conservation at all times. To achieve this, one needs to regard $X_t$ 
as a {control parameter} and one is lead to study the following 
{\em control problem}:
%
%
%

\vspace{0.3cm}
\framebox{
\begin{minipage}{.85\textwidth}
{
\smallskip
Is it possible to choose $X_t$ in such a way that
the total mass in the system
is constant ?
\smallskip
}
\end{minipage}
}

\vspace{0.3cm}

\noindent
We clearly succeed if
we can solve the free boundary problem (FBP) given by \eqref{intro.2} with
initial datum $\rho(r,0)=\rho_0(r)$, $r\in [0,X_0]$, and
  \begin{equation}
    \label{intro.4}
-\frac 12 \frac{\partial \rho}{\partial r}(0,t) = j,\qquad
-\frac 12 \frac{\partial \rho}{\partial r}(X_t,t) = j,\qquad
  \rho(X_t,t)=0\;.
   \end{equation}
In fact the rate at which mass is taken out of the system from $X_t$ is
  \begin{equation*}
J(X_t,t)=-\frac 12 \frac{\partial \rho}{\partial r}(X_t,t)
   \end{equation*}
which, by \eqref{intro.4}, is exactly equal to the rate at which we inject mass
at $0$ so that the total mass is constant.  

As discussed in the next chapter (see Section \ref{ch:2.2}) we can  find in the existing literature on FBP
an affirmative answer
for special initial data and for finite times.  In fact
one can readily check (see Section \ref{ch:2.2} for details) that the current $J(r,t)$ solves
the  classical Stefan problem for which  the theory
(in particular in one dimension) is
very rich with many detailed results available
 \cite{douglas, fasano, friedman, luckhaus, visintin}. As a consequence local existence and uniqueness
of classical solutions can be proved for the FBP \eqref{intro.2}--\eqref{intro.4}
for smooth initial data which satisfy the boundary conditions.  In some cases
the classical solution is global extending to
all times, but this is not true in general as it is known that singularities may develop.

Thus our control problem when stated for an arbitrarily long time interval $[0,T]$  and for general initial data cannot always be solved via the above FBP.  Take for instance $\rho_0 \in  L^1 (\mathbb R_+)$, bounded, continuous and everywhere strictly positive: in such a case the whole problem has to be redefined.
As usual the idea is to study a {\em relaxed}  version: we thus
introduce an accuracy parameter $\eps>0$  and replace
$\rho_0$ by a nice function $\rho^{(\eps)}_0$, smooth, non-negative and with compact support, requiring however that
$\int |\rho_0(r) - \rho^{(\eps)}_0(r) | \; dr  \le \eps$.  We may also ask that $\rho^{(\eps)}_0$
satisfies \eqref{intro.4} so that, for what said above, we have a classical
solution of FBP for some  time $[0,S]$.  However this could be shorter than the interval $[0,T]$ we have fixed initially, in which case the problem
still remains.  Moreover even if $S\ge T$ we have a poor control of the solution
and it is hard to see how this behaves when we remove the relaxation taking $\eps\to 0$.  The idea then is to further
simplify the problem by
relaxing also the boundary condition at the edge.  We refer to the next chapter for a precise definition. Here we just say that in Part I we will prove that any $\eps$-relaxed solution converges to a unique limit when $\eps\to 0$.
This will allow us to define a notion of relaxed solution of the problem which is global in time
and applies to a large class of initial data.

\medskip
In the last chapter of Part I we study a particles version of the above basic model.
The system has
$N$ particles so that  the mass distribution is no longer
continuous but instead concentrated on points (the positions of the $N$ particles).
To simulate an initial condition $\rho_0(r)$ (we assume $\int \rho_0(r)dr=1$ for simplicity), we
distribute the $N$ particles independently of each other and with law
$\rho_0(r)dr$.  We then define the ``empirical mass density measure''
 \begin{equation}
    \label{intro.6}
\pi_0^{(N)}(dr)= \frac 1N \sum_{i=1}^N \delta_{B_i(0) }(r)dr
   \end{equation}
where $B_i(0)$ are the random positions of the $N$ particles and $\delta_a(r)$ is the Dirac delta at $a$. The value $0$ refers to time, so far we have been describing the situation at
time 0.  Thus $\pi_0^{(N)}(dr)$ is a probability measure on $\mathbb R$ which is random as the terms $B_i(0)$ are the random positions of the particles.  If we denote by $E$ the expectation with respect to the law of the $B_i(0)$ and by $f(r)$ a test function, we have
   \begin{equation}
E\Big[ \int \pi_0^{(N)}(dr) f(r)\Big] = \int \rho_0(r) f(r) dr.
   \end{equation}
By the law of large number if $N$ is large we do not need to take the expectation because, with large probability,
$\int \pi_0^{(N)}(dr) f(r)$ is close to $\int \rho_0(r) f(r)dr$.

Let us now make the particles move.  We first consider the free case where the particles are independent Brownian motions $B_i(t)$ on $\mathbb R_+$ with reflection at 0. Call
 $\pi_t^{(N)}(dr)$ the random mass distribution at time $t$ and denote now by $E$
 the joint law of the initial distribution of the particles and of their Brownian evolution.
 We then have
   \begin{equation}
    \label{intro.7}
E\Big[ \int \pi_t^{(N)}(dr) f(r)\Big] = \int \rho(r,t) f(r)dr
   \end{equation}
where $\rho(r,t)$ is the solution of \eqref{intro.2} on $\mathbb R_+$ with
Neumann boundary condition at 0 given by 
$\frac{\partial  \rho}{\partial r} (0,t)=0$.  All that is the
well known relation between heat equation and Brownian motions.

We next go to the injection-removal of mass mechanism.  This is simply done as follows: at exponential times of intensity $jN$ the rightmost particle moves to the origin (which is the same as saying that
we add a new Brownian particle at 0 and simultaneously we take out the  particle which at that time is the rightmost one). In between such actions the particles move as independent Brownian motions (with reflection at the origin). We denote again by $E$ the
expectation with respect to the law of this process (which includes the initial distribution of the particles, their motion and the injection-removal of particles).  Thus the total mass (i.e.\ the
total number of  particles) is conserved but as in the continuum we are injecting mass at 0 and removing mass on the right.
Such a simple action however creates strong correlations among the particles: the choice of the rightmost particle requires knowledge of the positions of all the others. We thus lose the independency property and the analysis of the left-hand side of \eqref{intro.7} in this case becomes highly non-trivial. Existence of the process is easy but the relation with the continuum version is harder. The question becomes simpler if we study the asymptotic behavior of the system as $N\to \infty$, namely its
``hydrodynamic limit''.  We would like that:
   \begin{equation}
    \label{intro.7.1}
\lim_{N\to \infty}E\Big[ \int \pi_t^{(N)}(dr) f(r)\Big] = \int \rho(r,t) f(r) dr
   \end{equation}
where $ \rho(r,t)$ is the solution of the control problem described previously
and in particular of the FBP when this has a classical solution. In  Chapter \ref{ch:nuovo} we prove \eqref{intro.7.1}.

\renewcommand{\theequation}{\thesection.\arabic{equation}}
\setcounter{equation}{0}

\chapter{The basic model, definitions and results}
\label{ch:2}

In this chapter we expand the analysis presented in
the Introduction by giving a detailed definition of the control problem and its relaxed version.  We then show that for special initial conditions  the control problem
is related to a free boundary problem (FBP) which is solved locally in time using the existing literature on the Stefan problem.  We then  present
the main result of Part I (Theorem \ref{CDGPthm2.3.1}) which states that the relaxed
control problem has a unique global solution.  The proof uses inequalities
based on mass transport.  We introduce lower and upper barriers obtained by a time discretization of \eqref{intro.2}--\eqref{intro.3} and
state the other  main theorem of Part I (Theorem \ref{CDGPthm2.4.1}), which
says that there is a unique element which separates the lower and upper barriers.
The proof of Theorem \ref{CDGPthm2.4.1} starts in Chapter
\ref{ch:2bis} and is completed in Chapter \ref{ch:4}.  The proof of  Theorem \ref{CDGPthm2.3.1} is carried out in the remaining chapters of Part I, the essential point is to show that the elements of an optimal sequence are eventually squeezed between the barriers and therefore their limit points    coincide with the unique element which separates the barriers.

 \vskip2cm

\section{The basic problem}
\label{ch:2.1}

As discussed in the Introduction we  consider the heat equation \eqref{intro.2}
in the time varying domain $[0,X_t]$, $X_t$ a positive, continuous and piecewise
$C^1$ function,
with boundary conditions  \eqref{intro.3} and initial datum $\rho_0$.

\medskip
\begin{defin} [Assumptions on $\rho_0$]
\label{defin2.1.0}
 We suppose throughout the sequel that $\rho_0(r)$
 is a non-negative function
belonging to the set       
\begin{equation}
\label{CDGP3.1.1}
\mathcal U= \left\{u\in L^{\infty}(\mathbb R_+,\mathbb R_+)\cap L^{1}(\mathbb R_+,\mathbb R_+)\;:\; \int u > 0 \right\}.
   \end{equation}
\end{defin}
\medskip

\begin{defin} [The basic  problem]
\label{defin2.1.1}
The function $\rho(r,t)$  is a solution of the basic  problem 
in the time interval $[0,T]$ with initial datum $\rho_0$ if there exists $X_t$ positive so that
$\rho(r,t)$ solves \eqref{intro.2}--\eqref{intro.3} with initial condition $\rho_0$ and
\begin{equation}
\label{ch2.1.1}
 \int_{0}^{X_t}\rho(r,t)\,dr  =   \int_{0}^{X_0}\rho_0(r )\,dr \quad \text{for all $t> 0$}.
\end{equation}
\end{defin}

\medskip

\medskip

\begin{defin} [The $\eps$-relaxed problem]
 \label{defin2.1.2}
 For $\eps>0$, the function $\rho^{(\eps)}(r,t)$
is a $\eps$-relaxed solution of the basic problem with initial datum $\rho_0$ in the time interval $[0,T]$ if

\begin{itemize}

\item
$\int |\rho_0(r)-  \rho^{(\eps)}(r,0)| dr \le \eps$,

\item
there exists
$X^{(\eps)}_t$, $t\in [0,T]$, positive, continuous and piecewise
$C^1$ so that for each $t\in [0,T]$, $\rho^{(\eps)}(r,t)$ has support in $[0,X^{(\eps)}_t]$,

\item
$\rho^{(\eps)}(r,t)$ solves 
\eqref{intro.2}--\eqref{intro.3}  in $[0,T]$ with
 $X_t$ replaced by $X^{(\eps)}_t$ and with initial condition $\rho^{(\eps)}(r,0)$,

\item approximate mass conservation is satisfied, i.e.
\begin{equation}
\label{ch2.1.2}
 \Big|\int_{0}^{X^{(\eps)}_t}\rho^{(\eps)}(r,t)\,dr  -  \int_{0}^{X^{(\eps)}_0}\rho^{(\eps)}_0(r )\,dr\Big | \le \eps \quad \text{for all $t\in [0,T]$}.
\end{equation}
\end{itemize}

\end{defin}

\medskip

\begin{defin} [Optimal sequences]
\label{defin2.1.3}
The sequence $\rho^{(\eps_n)}(r,t)$ is an optimal sequence relative to
$\rho_0$ and $T>0$ if for each $n\in\mathbb{N}$ the function 
$\rho^{(\eps_n)}(r,t)$ is an $\eps_n$-relaxed
solution in $[0,T]$ of the basic problem with initial datum $\rho_0$
and if $\eps_n \to 0$ as $n\to\infty$.

\end{defin}

\medskip

\begin{defin} [Relaxed solution]
\label{defin2.1.3.1}
$\rho(r,t)$ is a relaxed solution in $[0,T]$
of the basic problem  with initial datum $\rho_0$   if
it is a weak limit of the elements $\rho^{(\eps_n)}(r,t)$ of an optimal sequence in $[0,T]$ with initial datum $\rho_0$.

\end{defin}

\vskip2cm

\renewcommand{\theequation}{\thesection.\arabic{equation}}
\setcounter{equation}{0}

\section{Stationary solutions}
\label{ch:2.2a}

The basic problem (see Definition \ref{defin2.1.1})
has special global solutions given by the stationary profiles:
\begin{equation}
\label{ch2.2a.1}
\rho^{(st)}(r |M ) = \Big(a(M) - 2j r\Big) \mathbf 1_{a(M) - 2j r \ge 0},\qquad \int\rho^{(st)}(r |M ) dr = M
\end{equation}
Since mass is conserved we have a one parameter family of stationary solutions indexed by the mass (denoted above by $M$).  We conjecture that these are the only stationary solutions but we do not have a proof.

In many problems stationary profiles are helpful because they can be used to ``trap'' trajectories and thus give a-priori estimates. We will prove that the relaxed solutions
of the basic problem (see Definition \ref{defin2.1.3.1}) preserve order
and this together with the knowledge of the stationary solutions will play an important role in the sequel.

\vskip2cm

\renewcommand{\theequation}{\thesection.\arabic{equation}}
\setcounter{equation}{0}

\section{The FBP for the basic model}
\label{ch:2.2}

$\rho_0$ is a {\em classical} initial datum if it is a smooth,
strictly positive function in $[0,X_0)$, $X_0>0$, and it is  such that
\begin{equation}
\label{ch2.2.1}
 \lim_{r\to X_0}\rho_0(r)=0,\quad \lim_{r\to 0}\frac{d\rho_0(r)}{dr}=- 2j,\; \lim_{r\to X_0}\frac{d\rho_0(r)}{dr}=- 2j.
\end{equation}

\medskip

\begin{theorem}[Local classical solutions]
\label{CDGPthm2.2.1}
If  $\rho_0$ is a classical initial datum then the basic problem (of Definition \ref{defin2.1.1})
has a local solution: namely
there exists $T>0$ and $\{X_t, t\in [0,T]\}$, so that \eqref{intro.2}
with initial datum $\rho_0$
has a solution $\rho(r,t)$ which satisfies \eqref{intro.4} for all $t\in[0,T]$.
If  $v(r,0) :=- \frac 12 \frac{\partial \rho_0}{\partial r}(r)  - j \ge 0$ then
the local solution extends to all times.

\end{theorem}

\medskip
The proof of  Theorem \ref{CDGPthm2.2.1} given below follows from the theory of   the Stefan problem as we are going to see.
The   equations  \eqref{intro.2} and  \eqref{intro.4} complemented by the initial datum $\rho_0$ in the unknowns $X_t$, and $\rho(\cdot,t)$ define a free boundary problem, FBP, where
the datum at the free boundary
involves both the value of $\rho$ and its  space derivative. In the Stefan problem, the prototype of FBP's, instead the datum is the speed of the edge:
 \begin{eqnarray}
\label{ch2.2.2}
&&\frac{\partial v}{\partial t}   =
\frac 12 \frac{\partial^2 v}{\partial r^2} \;\; 
\quad  v(r,t)\Big|_{r=0,X_t} = 0,
\nn
 \\&& \frac{dX_t}{dt} = -(2j)^{-1} \frac{\partial v(r,t)}{\partial r}\Big|_{r=X_t}.
\end{eqnarray}
 Local existence for \eqref{ch2.2.2} is proved in
\cite{fasano}--\cite {fp4}.

\medskip

\noindent
{\bf Proof of Theorem \ref{CDGPthm2.2.1}}.%
Given $X_t$ and  $v(r,t)$ satisfying \eqref{ch2.2.2} we set
 \begin{eqnarray}
\label{ch2.2.3}
&&\rho(r,t) = 2\int_{r}^{X_t} \Big(v(r',t)+j\Big)dr'.
\end{eqnarray}
One can then check that  \eqref{intro.2} and \eqref{intro.4} are all satisfied.
The non-negativity of $\rho(\cdot,t)$ follows from the maximum principle. 
Following Fasano and Primicerio, see e.g. \cite{fasano}, we say that  if $v(r,0) \ge 0$
then \eqref{ch2.2.2} has a ``sign specification''.  With a sign specification
the solution
is global hence the last statement in Theorem \ref{CDGPthm2.2.1}.
  \qed

\vskip.5cm

Uniqueness of the local classical solution for the  Stefan problem \eqref{ch2.2.2} is also known.
As mentioned  if there is  a sign specification
the solution of \eqref{ch2.2.2}
is global while if
there is no sign specification in general we only have local existence with
examples where singularities do appear.  The analysis of their structure  is a very interesting
and much studied problem, see for instance \cite{Crank}, \cite{fp5}, \cite{fp6},
\cite{Ockendon}.

\vskip2cm

\renewcommand{\theequation}{\thesection.\arabic{equation}}
\setcounter{equation}{0}

\section{Main theorem: existence and uniqueness}

By default throughout this chapter the initial datum $\rho_0\in \mathcal U$, see
Definition \ref{defin2.1.0}.

\medskip

\begin{theorem}[Existence and uniqueness]
\label{CDGPthm2.3.1}

Let $\rho_0\in \mathcal U$, then for any $T>0$ there exists a unique
relaxed solution of the basic problem  in $[0,T]$ with initial datum $\rho_0$ (see Definition \ref{defin2.1.3.1}).
Moreover:
%
%
%
%

\begin{itemize}
\item[(a)] As implicit in the above statement there exist optimal sequences in $[0,T]$ with initial datum $\rho_0$.

\item[(b)] The elements $\rho^{(\eps_n)}(r,t)$ of
an optimal sequence relative to
$\rho_0$ and $T>0$,
converge  weakly
to a limit $\rho_T(r,t)$. 

\item[(c)]
The limit $\rho_T(r,t)$ is independent of the optimal sequence  and if $S>T$, $\rho_S(r,t)= \rho_T(r,t)$, $t\in [0,T]$.  We denote
by $\rho(r,t)$ the function which agrees with $\rho_T(r,t)$ for all $T>0$.

\item[(d)] For all $t>0$ $\rho(r,t)$ is in $L^1$ and $\int \rho(r,t)dr
=\int \rho_0(r)dr$.

\item[(e)] If $\rho_0\le \rho^*_0$ then  $\rho(r,t)\le \rho^*(r,t)$.

\item [(f)]  $\rho(r,t)$ converges weakly to $\rho_0$ as $t\to 0$.

\end{itemize}
\noindent
Moreover, if $\rho_0(r)$ is continuous and with   support in $[0,X_0]$, then

\begin{itemize}
 \item[(g)]  
 $\rho(r,t)$   is a continuous function in $(r,t)$ which converges
 pointwise to $\rho_0$ as $t\to 0$.

  \item[(h)] If $\rho_0$ is a classical initial datum $\rho(r,t)$ solves  the FBP of Section \ref{ch:2.2}
  locally in time.

\end{itemize}

\end{theorem}

\vskip.5cm

Since any classical solution $\{(X_t, \rho(\cdot,t)), t\in [0,T]\}$, of the FBP \eqref{intro.2}--\eqref{intro.4}
is also an optimal sequence, (choosing $X^{(\eps)}_t=X_t$ and $\rho^{(\eps_n)}(r,t)=\rho(r,t)$ for any $n$), then
$\rho(\cdot,t)$ coincides with the function defined in Theorem \ref{CDGPthm2.3.1} and item (h) follows.

The weak point in the above theorem is the lack of control of the edge. We have only  what was stated in the following Corollary which is an immediate consequence of item (e) of Theorem \ref{CDGPthm2.3.1}
and of the existence of a stationary solution of the classical FBP as discussed in Section \ref{ch:2.2a}. Recall
\eqref{ch2.2a.1} for notation.

\medskip

\begin{coro}
If $\rho^{(st)}(r |M ')\le \rho_0(r) \le \rho^{(st)}(r |M '')$ then
\begin{equation}
\label{ch2.4.1-coro}
 \rho^{(st)}(r |M ')\le \rho(r,t) \le \rho^{(st)}(r |M '')\quad \text{for all $t>0$}.
\end{equation}
In particular if $\rho_0(r)$ has compact support then there exists $X>0$ so that
$\rho(r,t)=0$ for all $r\ge X$.

\end{coro}

We will prove Theorem \ref{CDGPthm2.3.1} using a variational method which is explained in the next sections.


\vskip2cm

\renewcommand{\theequation}{\thesection.\arabic{equation}}
\setcounter{equation}{0}

\section{The upper and lower barriers}

We do not have enough information on the elements $\rho^{(\eps_n)}(r,t)$
in an optimal sequence to directly prove that they converge as $\eps_n\to 0$.  We will
instead introduce a different relaxation procedure where the removal of mass occurs only at discrete times  $n\delta$, $n\in \mathbb N$, $\delta>0$. The evolution in the time intervals
$(n\delta, (n+1)\delta)$ is free, namely given by \eqref{intro.2} with only the boundary condition at $0$, i.e.\ the first one in  \eqref{intro.3}, the other one at $X_t$ is dropped. Therefore in these time intervals the mass density is strictly positive
on the whole $\mathbb R_+$.  At the times $n\delta$ we remove the right amount of mass, equal to $j\delta$, by cutting the right part of the function which after the cut has compact support.  Such evolutions are much simpler than those  in the optimal sequence but they have also the extra advantage  of monotone properties, this is why we  call them upper and lower barriers. Monotonicity
will allow us to control the limit as $\delta$ goes to 0 of the barriers and then to relate this to the limit of the
 $\rho^{(\eps_n)}(r,t)$.   We start here with the definition of the barriers.

  To this end
we introduce a time mesh $\delta>0$ and will define the barriers
at the times $k\delta$, $k\ge 0$. We use the following notation:
     \begin{equation}
\label{CDGP3.1.1789}
\mathcal U_\delta = \{u\in \mathcal U: \int u > j\delta\}
   \end{equation}
where $\mathcal U$ has been defined in \eqref{CDGP3.1.1}, and we introduce two operators, i.e.
the {\em cut operator} $C_\delta$ on $\mathcal U_\delta$ and
the {\em free evolution  operator} $T_\delta$ on $\mathcal U$.  

\medskip

\begin{defin} [The cut operator]
 \label{CDGPdefin2.4.1}
The cut operator $C_\delta$ maps $u\in \mathcal U_\delta$    into $\mathcal U$  as follows:
\begin{equation}
\label{ch2.4.1-def}
 C_\delta u (r)= \mathbf 1_{r\le R_u} u(r), \quad \text{where} \quad R_u: \;\;\int_{R_u}^\infty u(r) dr= j \delta.
\end{equation}
Observe that $\int  C_\delta u   =  \int u -j\delta$.

\end{defin}

\medskip

To define the free evolution operator we use the Green functions:

\medskip

\begin{defin} [The Green function]
 \label{CDGPdefin2.4.4}
Define for $r$, $r'$ and $t>s\ge 0$   the Green function
   \begin{equation}
\label{ch2.4.6}
G_{s,t}^{\rm neum}(r',r) = G_{t-s}(r',r) +
G_{t-s}(r',-r),\quad   G_t(r',r) =  \frac{e^{- \frac{(r-r')^2}{2t}}}{\sqrt{2\pi t}}
   \end{equation}
and write for any $u \in \mathcal U$:
\begin{equation}
\label{ch2.4.6.1}
G_{s,t}^{\rm neum}*u (r)= \int_{\mathbb R_+} G_{s,t}^{\rm neum}(r',r) u(r')\,dr'.
\end{equation}
To simplify notation we shall sometimes write
   \begin{equation}
\label{ch2.4.6.001}
G_{t-s}^{\rm neum}(r',r)=G_{s,t}^{\rm neum}(r',r).
   \end{equation}

\end{defin}

\medskip

The following proposition explains why $G_{t}^{\rm neum}$ is called the Green function.

\medskip

\begin{prop}
\label{CDGPprop2.4.1}
The function
$(r,t) \mapsto G_{t}^{\rm neum}(r',r)$,  $t>  0$, $r',r>0$, solves the heat equation \eqref{intro.2}  and  for any $t>0$,
  \[
 \lim_{r\to 0} \frac{\partial}{\partial r}G_{t}^{\rm neum}(r',r) = 0.
  \]
Moreover if  $u\in \mathcal U$ is a continuous function
  \begin{equation}
\label{ch2.4.7}
T_{t} u(r):= G_{t}^{\rm neum}*u(r)+ j
\int_0^t  G_{s',t}^{\rm neum}(0,r)\,ds'
   \end{equation}
solves   \eqref{intro.2},  converges to $u(r)$ as $t\to 0$ and for any $t>0$
  \begin{equation}
\label{ch2.4.7.1-prop}
\lim_{r\to 0}\frac{\partial}{\partial r}T_{t} u(r) = -2j.
   \end{equation}

\end{prop}

%
%
%
%
%
%
%
\noindent
{\bf Proof.}  The above statements are direct consequence
of the following properties of the Gaussian kernel. For $t>0$:
   \begin{eqnarray}
   \label{ch2.4.7.1-proof}
&& \Big(\frac{\partial }{\partial t} -\frac 12 \frac{\partial^2 }{\partial r^2}\Big)
\frac {e^{-r^2/2t}}{\sqrt{2\pi t}} =0,\\
  \label{ch2.4.7.2}
&& \lim_{t \to 0^+} \int_{-\eps}^{\eps}\frac {e^{-x^2/2t}}{\sqrt{2\pi t}} dx = 1, \quad \text{for any $\eps>0$},\\
&& \lim_{x \to 0^+} \int_0^t
\frac{x}{t-s}\frac {e^{-x^2/2(t-s)}}{\sqrt{2\pi (t-s)}}ds =1.
  \label{ch2.4.7.3}
  \end{eqnarray}
\qed

\medskip
We are now ready for the definition of the free evolution operator:

\medskip

\begin{defin} [The free evolution operator]
 \label{CDGPdefin2.4.2}
The free evolution operator  $T_\delta$ maps $\mathcal U$ into itself and
 $T_\delta u$
is equal to the expression \eqref{ch2.4.7} with $t=\delta$.

\end{defin}

\medskip

It follows directly from the definition that:
\begin{equation}
\label{ch2.4.2}
 \int T_\delta u =  j \delta + \int  u
\end{equation}
and therefore that
the products
\begin{equation}
\label{ch2.4.3}
\text{ $C_\delta T_\delta $ and  $T_\delta C_\delta $ preserve the mass}
\end{equation}
(the latter defined on $\mathcal U_\delta$ ).

\medskip

\begin{defin} [The barriers]
 \label{CDGPdefin2.4.3}
The upper barrier $S^{\delta,+}_{k\delta}u$, $k\ge 0$,
is defined as
\begin{equation}
\label{ch2.4.4}
 S^{\delta,+}_{k\delta} u = (T_\delta C_\delta)^k u,\quad  u \in \mathcal U_\delta
\end{equation}
while the lower  barrier $S^{\delta,-}_{k\delta} u$, $k\ge 0$, is
\begin{equation}
\label{ch2.4.5-barrier}
 S^{\delta,-}_{k\delta} u =( C_\delta T_\delta)^k  u,\quad u \in \mathcal U.
\end{equation}

\end{defin}

\medskip

\begin{prop}
\label{CDGP2.5.0.11}
At all times the barriers have the same total mass as initially:
\begin{equation}
\label{ch2.4.5-mass}
 F(0;S^{\delta,\pm}_{k\delta} u) =F(0,u)
\end{equation}
where
 \begin{equation}
\label{ch2.4.9}
F(r;u)=   \int_r^\infty u(r')\,dr',\quad r \ge 0.
   \end{equation}

\end{prop}

\medskip
In Chapter \ref{ch:2bis} we will prove that the upper
barriers are equi-bounded and equi-continuous as functions of $(r,t)$ on $[\eps,T]\times \mathbb R_+$, for any $\eps>0$ and $T>0$ which yields convergence by subsequences.
To gain full convergence we will use inequalities based on the order by mass-transport,   defined in the next section.

\vskip2cm

\renewcommand{\theequation}{\thesection.\arabic{equation}}
\setcounter{equation}{0}

\section{Mass transport}

We use a notion of order (in the sense of mass transport) under which we will prove that the upper barriers are larger than the lower barriers
and that the convergence as $\delta\to 0$ is monotone.  Inequalities by the above order will be of paramount importance
in the proof of Theorem \ref{CDGPthm2.3.1} as we will   show that the elements
$\rho^{(\eps_n)}(r,t)$ in an optimal sequence are eventually squeezed (as $\eps_n\to 0$) between the upper and the lower barriers.
Observe that the notion of barriers for the construction of solutions
of partial differential equations
is well known   \cite{L, friedman2} (see also
\cite{DG} in the context of motion by mean curvature).
The notion of order that we   use to define  upper and lower barriers is:

\medskip
   \begin{defin}  [Partial order]
\label{CDGPdefin2.4.5}
For any $u,v \in \mathcal U $
we set
   \begin{equation}
\label{ch2.4.8}
u \preccurlyeq v
 \;\;\;\text{iff}\;\;\;
F(r; u)    \le F(r; v) \;\;\;\text{for all $r\ge 0$}
   \end{equation}
where  $F(r;u)$ is defined in \eqref{ch2.4.9}.
\end{defin}

\medskip

When $u$ and $v$ have the same total mass, then $u \preccurlyeq v$ if and only if $v$ can be obtained from $u$ by moving mass to  the right. This statement will be made precise in Proposition \ref{prop4.1}, hence the above partial order is related to mass transport.

%
%

The next theorem justifies the name of upper and lower barriers.  We first
consider a very special case namely the inequality
\begin{equation}
\label{ch2.4.5.1}
 S^{\delta,-}_{ \delta} u \preccurlyeq S^{\delta,+}_{ \delta} u
\end{equation}
whose proof we hope will give a feeling of what is going on.  
Define $u_1$ by writing 
\begin{equation*}
   u = C_\delta u + u_1,\quad u_1= u-C_\delta u.
\end{equation*}
Recalling the definition of $C_\delta$, $u_1$ has mass $j\delta$ which is
to the right of the mass of $C_\delta u$.  By \eqref{ch2.4.7}
  \begin{equation*}
v:=  T_\delta u = G^{\rm neum}_\delta *\{C_\delta u + u_1\} + j
\int_0^\delta  G_{\delta-s}^{\rm neum}(0,r)\,ds.
\end{equation*}
Then  $S^{\delta,-}_{ \delta} u$ is obtained from $v$ by cutting  a mass
$j\delta$ to the right of $v$, while $S^{\delta,+}_{ \delta} u$ is obtained from
$v$ by erasing  the term $G^{\rm neum}_\delta * u_1$: thus $S^{\delta,+}_{ \delta} u$
is obtained from $S^{\delta,-}_{ \delta} u$ by moving mass to the right, hence
\eqref{ch2.4.5.1}. More details can be found in the proof of Lemma \ref{lemma4.5}.

The following theorem is the key step in the proof of Theorem
\ref{CDGPthm2.3.1}. Its content is divided in three parts: inequalities among barriers, convergence theorems and properties of the limit. It is proved in Chapters
 \ref{ch:2bis}, \ref{ch:4n}, \ref{ch:3} and \ref{ch:4}, a summary is given in Section \ref{ottavabis}.

 \medskip

\section{Barrier theorems}

    \begin{thm}[Barriers and separating elements]
           \label{CDGPthm2.4.1}
Let $u\in \mathcal U$ and $t>0$. 

\noindent
Inequalities among barriers:

\begin{itemize}

\item [(1)]   If $u \in \mathcal U_\delta$ then
\begin{equation}
\label{cambi.1}
S_{t}^{\delta,-} u \preccurlyeq S_{t}^{\delta',+} u, \quad t=k\delta=n\delta',\; k,n \in \mathbb N.
\end{equation}

\item [(2)]  For any $\delta>0$, $u \in \mathcal U_\delta$ and $t=k\delta, k\in \mathbb N$
\begin{equation}
\label{cambi.0}
\int |S_{t}^{\delta,-} u(r) - S_{t}^{\delta,+} u(r)| dr \le 2 j \delta.
\end{equation}

\item [(3)] For $n$ so large that
$u\in \mathcal U_{2^{-n}t}$ and for  all $r\ge 0$, $F(r;S_{t}^{2^{-n}t,-} u)$ is a non-decreasing function of $n$  and $F(r;S_{t}^{2^{-n}t,+} u)$ is a non-increasing function of $n$.  Moreover, as proved in \eqref{ch2.4.5-mass}, $F(0;S_{t}^{2^{-n}t,\pm} u)= F(0;u)$.

\end{itemize}

\noindent
Convergence:

\begin{itemize}

\item [(4)] There exists a bounded function $S_t u(r)$ continuous in $(r,t)$ for $t>0$
such that  $S_{t}^{2^{-n}t,+} u(r)$ converges to $S_t u$ uniformly in the compacts of
$(r,t) \in \mathbb R_+\times (0,\infty)$ and for all $t>0$
 $S_{t}^{2^{-n}t,+} u $ converges to $S_t u$ in $L_1$.

\item [(5)]  The convergence is monotone in the mass transport order of
Definition
 \ref{CDGPdefin2.4.5}:
    \begin{eqnarray}
&&\hskip-1cm
F(r;S_{t} u) = \lim_{n\to \infty} F(r;S_{t}^{2^{-n}t,\pm} u),
\label{4.4a}
\end{eqnarray}
hence, by (3), $F(0;S_{t} u) = F(0;u)$.

\end{itemize}

\noindent
Properties of $S_tu$:

\begin{itemize}

\item [(6)]  $S_t u$ separates the barriers:
  \begin{eqnarray}
  &&\hskip-1cm
  F(r;S_{t} u) =\inf_{\delta: t=k\delta, k\in \mathbb N} F(r;S_{t}^{\delta,+}u)
=\sup_{\delta: t=k\delta, k\in \mathbb N} F(r;S_{t}^{\delta,-}u).
\label{4.4b}
  \end{eqnarray}

\item   [(7)] $S_{t} u \to u$ weakly as $t\to 0$ and if $u$ is continuous with compact support then $S_{t} u \to u$ point-wise as $t\to 0$.

\item [(8)] If $u \preccurlyeq v$
then
$S_t u \preccurlyeq S_t v$.

\item   [(9)] If   $ u \le v$ point-wise then
$S_{t} u \le S_{t} v$ point-wise for all $t> 0$.

\end{itemize}

\end{thm}

\vskip.5cm

The first step in the proof of  Theorem \ref{CDGPthm2.3.1} after
Theorem \ref{CDGPthm2.4.1} is the following identification theorem:

\medskip

    \begin{thm}[Identification theorem]
           \label{CDGPthm2.4.1.00}
   For any $T>0$ and $u \in \mathcal U$ there exist relaxed solutions of
   the basic problem in $[0,T]$ with initial datum $u$ and they are all equal to $S_tu$.

\end{thm}

The proof of Theorem \ref{CDGPthm2.4.1.00} is the most original part of this work. It uses extensively probability ideas and techniques as it
relies on the representation of the solution of the heat equation with Dirichlet boundary conditions in terms of Brownian motion and its hitting distribution at the boundary.  After
 showing in Chapter \ref{ch:4b} the existence of optimal sequence, we prove in
  Chapter \ref{ch:7}
 that  given any $\delta>0$
the elements $\rho^{(\eps_n)}$ of an optimal sequence in the limit
as $n\to \infty$ are squeezed in between $S^{\delta,\pm}_{ t} \rho_0$.
By the arbitrariness of $\delta$
this implies that $\rho^{(\eps_n)}$ converges weakly,   its limit being from one side
 equal to $S_tu$ while, from the other side, is
 by definition a relaxed solution, hence Theorem \ref{CDGPthm2.4.1.00}.  Thus the relaxed solution
inherits all the properties of the separating element stated in Theorem \ref{CDGPthm2.4.1}
which allows us to complete the proof of  Theorem \ref{CDGPthm2.3.1}, see Section \ref{ch:7.6}.

\renewcommand{\theequation}{\thesection.\arabic{equation}}
\setcounter{equation}{0}

\chapter{Regularity properties of the  barriers}
\label{ch:2bis}

In this chapter we will prove some regularity properties of
the  barriers $S_{t}^{\delta,\pm}u$,
$u\in \mathcal U_\delta$.
By the smoothness of $G_{t}^{\rm neum}(r,r')$, $t>0$, it is easy
to prove that for any $n>0$,
$S^{\delta,+}_{n\delta}u\in C^\infty(\mathbb R_+)$ while
$S^{\delta,-}_{n\delta}u$
is
$C^\infty$ in the interior of its support.
Such a smoothness however, being inherited from
$G^{\rm neum}_\delta$, depends on $\delta$, while we want
properties which hold uniformly as $\delta\to 0$.
The main results in this section is  that the family
$S_{t}^{\delta,+}u(r)$ is equi-bounded and equicontinuous in space-time for $t$
away from 0, these statements are proved in the following three sections.

\renewcommand{\theequation}{\thesection.\arabic{equation}}
\setcounter{equation}{0}

\vskip2cm

\section{Equi-boundedness}
\label{ch:2bis.1}

We denote by $\|u\|_{\infty}$ and $\|u\|_{1}=F(0;u)$ the $L^\infty$ and $L^1$ norm of $u \in \mathcal U$.

\medskip 

\begin{thm}
\label{thmee7.1}
There is a constant $c$ so that the following holds.
Let $\delta>0$ and $u\in \mathcal U_\delta$, then
 \begin{equation}
	\label{e6.1}
\|S^{\delta,\pm}_{t}u\|_{\infty} \le
c\left\{
\begin{array}{ll}
j+ \|u\|_{\infty} & \text{for all $t \in \delta \N, \: t \le 1$},\\
j+ \|u\|_1  & \text{for all  $t \in \delta \N, \: t>1$}.
\end{array}
\right.
\end{equation}

\end{thm}

\medskip

\noindent
{\bf Proof.}
Let $t=n\delta $, $n$ a positive integer, then
  \[
  S^{\delta,\pm}_{t}u (r) \le
  \int dr'G_{0,\delta}^{\rm neum}(r',r) S^{\delta,\pm}_{t-\delta}u(r') + j \int_{t-\delta}^t ds
G_{s,t}^{\rm neum}(0,r).
  \]
The inequality is because we are neglecting $C_\delta$.
 Iterating we get for  $0\le m< n$,
	\begin{equation}
	\label{e6.2}
S^{\delta,\pm}_{t} u(r) \le \int dr'G_{m\delta ,t}^{\rm neum}(r',r) S^{\delta,\pm}_{m\delta}u(r') + j \int_{m\delta}^t ds
G_{s,t}^{\rm neum}(0,r).
	\end{equation}
Let $t\le 1$, take $m=0$ in \eqref{e6.2} then
	\begin{equation*}
S^{\delta,\pm}_{t} u(r) \le \int dr'G_{0 ,t}^{\rm neum}(r',r) \|u\|_\infty + j \int_{0}^t ds \frac{2}{\sqrt{2\pi (t-s)}}
	\end{equation*}
which proves \eqref{e6.1} when $t\le 1$.

Let
$n_\delta$ be the smallest integer such that
$\tau:=\delta n_\delta\ge 1$.  Let  $t\in [k\tau, (k+1)\tau]$ and $m\delta=(k-1)\tau$ in \eqref{e6.2}.  Then
	\begin{eqnarray*}
S^{\delta,\pm}_{t} u(r) &\le& \int dr'\frac{2}{\sqrt{2\pi \tau}} S^{\delta,\pm}_{(k-1)\tau} u(r') + j \int_{(k-1)\tau}^t ds \frac{2}{\sqrt{2\pi (t-s)}} \\ &\le& c (j+ \|S^{\delta,\pm}_{(k-1)\tau} u\|_1).
	\end{eqnarray*}
By Proposition \ref{CDGP2.5.0.11} $\|S^{\delta,\pm}_{t} u\|_1= \|u\|_1$  hence  \eqref{e6.1}.  \qed

\vskip1cm

\renewcommand{\theequation}{\thesection.\arabic{equation}}
\setcounter{equation}{0}

\section{Space equi-continuity}

In this section we will prove that the family
$\{S^{\delta,+}_{t}u(r)\}$ is equi-continuous in $r$
for any fixed $t>0$.  We need a preliminary lemma  where
we use the following notation:
 	\begin{equation}
	\label{e6.333333}
w_{s,t}^{\delta,+}(r) := \int dr'G_{t-s}^{\rm neum}(r,r') S^{\delta,+}_{s} u(r'),\;\;\;
v_{s,t}^{\delta,+}:=  S^{\delta,+}_{t}u - w_{s,t}^{\delta,+}.
	\end{equation}

\medskip

\begin{lem}
\label{lemmae6.24}
There is a constant
$c$ so that the following holds. For all  $\delta>0$, $u\in \mathcal U_\delta$,
$0\le s< t$, $s,t \in\delta \mathbb N$, $t-s\le 1$,
  	\begin{equation}
	\label{e6.3}
\|\frac{\partial}{\partial r} w_{s,t}^{\delta,+}(r)\|_\infty \le c
\frac{\|u\|_\infty+\|u\|_1+j}{ \sqrt{ t-s} },
	\end{equation}
 	\begin{equation}
	\label{e6.3.0}
\|v_{s,t}^{\delta,+}\|_1 \le 2j (t-s),\quad  \|v_{s,t}^{\delta,+}\|_\infty \le c  j \sqrt{t-s}.
	\end{equation}

\end{lem}

\medskip

\noindent
{\bf Proof.}
By \eqref{e6.1}
	\begin{equation*}
|\frac{\partial}{\partial r} w_{s,t}^{\delta,+}(r)|\le c\|S^{\delta,+}_{s}u\|_{\infty}\int dr' \frac{|r-r'|}{t-s}
G^{\rm neum}_{s,t}(r',r) \le c' \frac{\|u\|_\infty +\|u\|_1 +j}{ \sqrt{ t-s}  }
	\end{equation*}
which proves \eqref{e6.3}.

By \eqref{e6.2} with $s=m\delta$,
	\begin{equation}
	\label{3.2.3bis}
v_{s,t}^{\delta,+} \le  j \int_{s}^t ds'
G_{s',t}^{\rm neum}(0,r).
	\end{equation}
To get a
lower bound we first define, for any $ \tau \in \delta \mathbb N$,
\begin{equation}
\label{a1}
v^{(\delta)}_\tau (r) := \mathbf 1_{r\ge R}S^{\delta,+}_\tau u (r),
\quad R:\;\;
  \int_R^\infty S^{\delta,+}_\tau u (r)= j\delta.
\end{equation}
By \eqref{e6.1}
\begin{equation}
\label{a2}
\|v^{(\delta)}_\tau \|_\infty \le C, \quad C=c\Big(j+ \|u\|_{\infty}+\|u\|_1\Big).
\end{equation}
By neglecting the contribution of the mass injection we get:
 	\begin{equation*}
S^{\delta,+}_{t}u \ge
  G^{\rm neum}_{\delta}* \Big(S^{\delta,+}_{t-\delta} u - v^{(\delta)}_{t-\delta}\Big)
	\end{equation*}
and, calling $s=m\delta$,
 	\begin{equation*}
 S^{\delta,+}_{t}u  \ge    G^{\rm neum}_{t-s}* S^{\delta,+}_{m\delta} u
 - \sum_{k= m} ^{ n-1}  G^{\rm neum}
 _{ (n- k)\delta} * v^{(\delta)}_{k\delta}.
	\end{equation*}
This together with
\eqref{3.2.3bis} gives
 	\begin{equation}
	\label{e6.4aa.1}
|v_{s,t}^{\delta,+}(r)|  \le \sum_{k= m} ^{ n-1}  G^{\rm neum}
 _{ (n- k)\delta} * v^{(\delta)}_{k\delta}(r)+  j \int_{s}^t ds'
G_{s',t}^{\rm neum}(0,r).
	\end{equation}
Recalling that $\int v^{(\delta)}_{k\delta} =j\delta$,
	\begin{eqnarray*}
\| \sum_{k=m}^{n-1}  G^{\rm neum}_{ (n- k)\delta} * v^{(\delta)}_{k\delta}  \|_\infty \le  \sum_{k=m}^{n-1} \frac{1}{\sqrt{2\pi(n-k)\delta}}j\delta
\le c j\sqrt{t-s}.
 	\end{eqnarray*}
Then we have: 
	\begin{eqnarray*}
\| j \int_{m\delta}^t ds'
G_{s',t}^{\rm neum}(0,r)\|_\infty \le 
c j\sqrt{t-s}
 	\end{eqnarray*}
so that
	\begin{equation}
	\label{CDGP4.5.16}
\|v_{s,t}^{\delta,+}\|_\infty \le c  j \sqrt{t-s}
	\end{equation}
and the second inequality in \eqref{e6.3.0} is proved.  To prove the first one
we use \eqref{e6.4aa.1} and
\eqref{e6.333333}  
to write
\begin{eqnarray*}
\|v^{\delta,+}_{s,t}\|_1\le j (t-s)+ \sum_{k=m}^{n-1} \int dr \int dr'  G^{\rm neum}_{ (n- k)\delta}(r,r') v^{(\delta)}_{k\delta} (r') = 2j (t-s)
 	\end{eqnarray*}
which concludes the proof of \eqref{e6.3.0}.

\qed

\begin{thm}
\label{thmZ3.2.0.15}
Given $u \in \mathcal U $, for any $\delta>0$, $t\in \delta \mathbb N$ and   $\zeta>0$
there is $d>0$ which depends on $\zeta,t $, $\|u\|_\infty$ and $\|u\|_1$ so that
  	\begin{equation}
	\label{e6.1aaa.1}
| S^{\delta,+}_{t}u(r)-S^{\delta,+}_{t}u(r')| <\zeta,\quad
|r-r'|<d.
	\end{equation}

\end{thm}

\medskip

\noindent
{\bf Proof.}
  By \eqref{e6.333333} 
 we can write
   	\begin{equation*}
S^{\delta,+}_{t}u (r)- S^{\delta,+}_{t}u (r')= w_{s,t}^{\delta,+}(r')-
 w_{s,t}^{\delta,+}(r) +v_{s,t}^{\delta,+}(r')-v_{s,t}^{\delta,+}(r).
	\end{equation*}
 We then use \eqref{e6.3} and \eqref{e6.3.0} and get
  for any $s\in \delta \mathbb N$, $0\le s <t$,
 	\begin{equation}
	\label{e6.3.0.1.2}
| S^{\delta,+}_{t}u(r)- S^{\delta,+}_{t}u (r')| \le
c (\|u\|_\infty+\|u\|_1+j) \frac{|r-r'|}{ \sqrt{t-s}} + 2 c j\sqrt{t-s}.
	\end{equation}
We choose $s$ so that $2 cj\sqrt{t-s} < \zeta/2$ and for such a value of $s$
we take $d$ so that
\[
 c (\|u\|_\infty+\|u\|_1+j)  \frac{d}{ \sqrt{t-s}} < \frac \zeta 2.
\]
\qed

\vskip1cm

\renewcommand{\theequation}{\thesection.\arabic{equation}}
\setcounter{equation}{0}

\section{Time equi-continuity}

\begin{thm}
Let $u\in \mathcal U$, then for any $\delta>0$, $t \in \delta \mathbb N$,
$t>0$,
$\zeta>0$ there is $\tau=\tau_{\zeta,t}>0$ so that
 	\begin{equation}
	\label{e6.1aaa.1.00}
\| S^{\delta,+}_{t'}u -S^{\delta,+}_{t}u \|_\infty <\zeta, \quad
t' \in  \delta\mathbb N\cap (t,t+\tau_{\zeta,t}).
	\end{equation}

\end{thm}

\medskip

\noindent
{\bf Proof.}
By \eqref{CDGP4.5.16} 
 	\begin{eqnarray*}
&&| S^{\delta,+}_{t'}u(r)  -  G^{\rm neum}_{t,t'}*S^{\delta,+}_{t}u(r)| \le
    cj\sqrt{t'-t}.
	\end{eqnarray*}
Let $\zeta'<\zeta$ and $d'$ the corresponding constant in
Theorem \ref{thmZ3.2.0.15}, then
 	\begin{eqnarray*}
&&| S^{\delta,+}_{t'}u(r)  - S^{\delta,+}_{t}u(r) | \le \|S^{\delta,+}_tu\|_\infty \int_{r':|r-r'|\ge d'}   G^{\rm neum}_{t,t'}(r',r) \,dr'\\&&\hskip2cm +  \zeta' + cj\sqrt{t'-t}.
	\end{eqnarray*}
If
$\delta>\tau$ there is no $t': t<t'<t+\tau$ and
\eqref{e6.1aaa.1.00}
is automatically satisfied.  Let then $\delta\le\tau$.
There is a constant $c_1$ so that
 	\begin{eqnarray*}
\int_{r':|r-r'|\ge d'}  G^{\rm neum}_{t,t'}(r,r') \,dr' \le c_1 e^{- (d ')^2/(4\tau)}
	\end{eqnarray*}
Thus
 	\begin{eqnarray*}
&&| S^{\delta,+}_{t'}u(r)  - S^{\delta,+}_{t}u(r) | \le c_1\|S^{\delta,+}_t u\|_\infty  e^{- (d')^2/(4\tau)}  +  \zeta' + cj\sqrt{\tau}
	\end{eqnarray*}
which concludes the proof of the theorem.

\qed

\renewcommand{\theequation}{\thesection.\arabic{equation}}
\setcounter{equation}{0}

\chapter{Lipschitz and $L^1$ estimates}
\label{ch:4n}

In this  chapter  we first prove some elementary inequalities and then Lipschitz estimates for the operators involved in the definition of barriers.  We finally prove that upper and lower barriers are $L^1$ close, proportionally to $\delta$.


\renewcommand{\theequation}{\thesection.\arabic{equation}}
\setcounter{equation}{0}
\vskip1cm

\section{Elementary inequalities}
\label{ch3.0}

Recall that    $\mathcal U_\delta$
is defined in \eqref{CDGP3.1.1789}.

\medskip

\begin{prop}
\label{propZ4.1.0.18}
Let $u\in \mathcal U_\delta$, then
   \begin{equation}
\label{CDGP4n.0.1}
 u \ge C_\delta u \quad \text{point-wise}.
   \end{equation}
Moreover if $u,v\in \mathcal U_\delta$ and $u \le v$  point-wise,  then
   \begin{equation}
\label{CDGP4n.0.2}
  \quad C_\delta u \le C_\delta v,\quad  T_\delta u \le
  T_\delta v,\quad S^{\delta,+}_{n\delta}u \le S^{\delta,+}_{n\delta}v,  \; \text{point-wise}.
   \end{equation}

\end{prop}

\medskip

\noindent
{\bf Proof.} \eqref{CDGP4n.0.1} follows immediately from the definition of $C_\delta$, namely $ C_\delta u (r)= u\mathbf 1_{r\le R_u}$, $R_u:$ $\int_{R_u}^\infty
u = j\delta$.  If $u\le v$ then $ R_u \le  R_v$ hence $C_\delta u \le C_\delta v$.  Recalling the definition of $T_\delta$ we get
\[
 T_\delta v(r)- T_\delta u(r)= \int_0^{\infty} dr' [v(r')-u(r')] G^{\rm neum}_\delta(r',r)
\]
which is non-negative. The last inequality in \eqref{CDGP4n.0.2} with $n=1$ follows from the previous ones which we have already proved.  By induction
the inequality is then proved for $n>1$ as well.
\qed

\renewcommand{\theequation}{\thesection.\arabic{equation}}
\setcounter{equation}{0}

\vskip2cm

\section{Lipschitz properties}
\label{ch3.1}

 Recall that $\mathcal U$,  $\mathcal U_\delta$ and $F(r;u)$ are defined in \eqref{CDGP3.1.1} , \eqref{CDGP3.1.1789} and
 \eqref{ch2.4.9}, respectively. We also write  $|f|_1= \|f\|_1$ for the $L^1$ norm of $f$.

\medskip

\begin{prop}
\label{prop3.1}

Let $\delta>0$ and let $u\in \mathcal U_\delta$, then
   \begin{equation}
\label{CDGP3.1.3}
|u-C_\delta u|_1 = j\delta,\quad
F(0; u) = F(0;C_{\delta}u)+j\delta=F(0;T_{\delta} u)-j\delta.
   \end{equation}
Let also $v\in \mathcal U_\delta$, then
      \begin{equation}
\label{CDGP3.1.4}
|C_{\delta}u -C_{\delta} v|_1\le |u-v|_1,\quad
|T_ \delta u-T_\delta v|_1\le |u-v|_1,
   \end{equation}
   \begin{equation}
	\label{CDGP3.1.5}
|S^{\delta,\pm}_{k\delta} u - S^{\delta,\pm}_{k\delta}v|_1 \le |u-v|_1,\qquad \text{for
all $k\in \mathbb N$}.
\end{equation}

\end{prop}

\medskip

\noindent
{\bf Proof.} The second and third equalities in \eqref{CDGP3.1.3} follow directly from the definition of the operators $C_\delta$ and $T_\delta$.  The first one follows from
\eqref{CDGP3.1.3} and the second one.
To prove the first inequality
in \eqref{CDGP3.1.4}  we recall \eqref{ch2.4.1-def} for the definition of $R_u$ and,
assuming that $R_u\le R_v$,
      \begin{equation*}
|C_{\delta}u -C_{\delta} v|_1=  \int_0^{R_u} |u-v|
+ \int_{R_u}^{R_v} v.
   \end{equation*}
We can then add $\dis{\int_{R_v}^\infty v}$ and subtract
   $\dis{\int_{R_ u}^\infty u}$ as they are both equal to $j\delta$:
      \begin{equation*}
|C_{\delta}u -C_{\delta} v|_1=  \int_0^{R_u} |u-v|
+ \int_{R_ u}^{\infty}  v-\int_{R_ u}^{\infty}u  \le |u-v|_1.
   \end{equation*}
To prove the second inequality we use \eqref{ch2.4.7} so that
      \begin{eqnarray*}
|T_ \delta u-T_\delta v|_1&\le& \int dr| \int dr' G^{\rm neum}_\delta (r',r)u(r')-G^{\rm neum}_\delta (r',r)v(r')|\\
& \le& \int dr\int dr' G^{\rm neum}_\delta (r',r)|u(r')-v(r')|
   \end{eqnarray*}
which is equal to $|u-v|_1$ because $\int G^{\rm neum}_\delta (r',r) dr=1$.
  \eqref{CDGP3.1.5} is a direct consequence of \eqref{CDGP3.1.4}.
\qed

\renewcommand{\theequation}{\thesection.\arabic{equation}}
\setcounter{equation}{0}

\vskip2cm

\section{$L^1$ estimates}
\label{ch4n.2}

 In this section we prove that the upper and lower barriers are $L^1$ close
 for small $\delta$.

\medskip

\begin{thm}
\label{CDGPmthm5.3}
Let  $u  \in  \mathcal U_{\delta}$, then for all $k\in \N$,
\begin{eqnarray}
&&
|S_{k\delta}^{\delta,+}u - S_{k\delta}^{\delta,-}u |_1\le 2j\delta.
\label{4.3a.00}
  \end{eqnarray}

  \end{thm}

\medskip
\noindent
{\bf Proof.}
We need to bound  $\int |\phi-\psi|$, where
 \[
\phi:= C_{\delta}T_\delta \cdots  C_{\delta}T_\delta u,\qquad \psi:=T_\delta  C_{\delta}\cdots T_\delta C_{\delta} u \qquad {\text{$k$ times}}.
 \]
Call
 \[
 v = C_{\delta}u,\qquad v_k= T_\delta C_{\delta}\cdots T_\delta v,\qquad u_k=T_\delta C_{\delta}\cdots T_\delta u \qquad {\text{$k$ times}}
 \]
so that $\phi= C_{\delta}u_k$ and $\psi=v_k$.  Hence using   \eqref{CDGP3.1.3},
 \begin{eqnarray*}
 |\psi-\phi|_1 &{=}& |C_{\delta}u_k - v_k|_1 \le |C_{\delta}u_k - u_k|_1 +
|v_k-u_k|_1\\
  &{=}& j\delta + |v_k-u_k|_1 \le j\delta +|u-v|_1 = 2j\delta.
 \end{eqnarray*}
\qed


\renewcommand{\theequation}{\thesection.\arabic{equation}}
\setcounter{equation}{0}

\chapter{Mass transport inequalities}
\label{ch:3}

 We present in this  chapter   some   known
facts about mass transport and use them to prove   properties which will then be
extensively used in the sequel.

\renewcommand{\theequation}{\thesection.\arabic{equation}}
\setcounter{equation}{0}

\vskip2cm

\section{Partial order and mass transport}
\label{ch3.2}

In this section we relate the notion of partial order discussed so far to the notion of mass transport.  To define the latter, consider a non-decreasing
map $f:\mathbb R_+ \to
\mathbb R_+$ and interpret $f(r)$ as the position of   $r$ after the ``displacement''.  Moving mass to the right then means that $f(r)\ge r$ for all $r$.  If there was initially a mass $M$ in an interval $[a,b]$, then after the displacement there will be a mass $M$ in the interval $[f(a),f(b)]$.  Thus if the initial mass density is $u$ then the final mass density $v$ is such that for any $a<b$,
\[
\int_a^b u = \int_{f(a)}^{f(b)} v.
\]
As a consequence
\[
F(a;u) =\int_a^\infty u = \int_{f(a)}^{\infty} v \le F(a;v)
\]
(because $f(a)\ge a$).  Thus if $v$ is obtained from $u$ by moving mass to the right then $u\preccurlyeq v$.  The converse is proved next:

\medskip

{\begin{prop}[The mass displacement lemma]
\label{prop4.1}
Given  $u \preccurlyeq  v$ in $ \mathcal U$  with $F(0;u)=F(0;v)$ we define for $r\in \mathbb R_+$:
  \begin{equation}
\label{4.5a}
f(r) := \sup\;\Big\{r': \int_0^{r'} v(z)\,dz= \int_0^r u(z)\,dz\Big\}.
   \end{equation}
Then
  \begin{equation}
\label{4.6a}
f(r) \ge r
   \end{equation}
and for any function $\phi\in L^\infty(\mathbb R_+,\mathbb R)$,
  \begin{equation}
\label{4.7a}
 \int_0^\infty v(r)\phi(r)\,dr = \int_0^\infty u(r) \phi(f(r))\,dr.
   \end{equation}
\end{prop}}

\medskip

\noindent
{\bf Proof.}  Since $F(0;u)= F(0;v)$,
\[
\int_0^r u(z) \,dz + F(r;u)= \int_0^r v(z) \,dz + F(r;v)
\]
and since $F(r;u)\le F(r;v)$,
\[
\int_0^r u(z) \,dz  \ge \int_0^r v(z) \,dz
\]
which yields \eqref{4.6a}.
By a density argument \eqref{4.7a} follows from \eqref{4.6a}.

\qed

\medskip

\begin{coro}
\label{4.2.0.18}
Let  $u  \preccurlyeq v$ in $ \mathcal U$ and $F(0;u)=F(0;v)$, then for all bounded, non-decreasing
functions $h$ on $\mathbb R_+$:
  \begin{equation}
\label{4.7aaa}
  \int_0^\infty u(r) h(r) \,dr  \le \int_0^\infty v(r) h(r)\,dr.
   \end{equation}

\end{coro}

\medskip

\noindent
{\bf Proof.}  Observe that \eqref{4.7aaa} is verified by definition for all
functions $h$ of the form $\mathbf 1_{[R,\infty)}$, $R\ge 0$.  Its validity for
functions $h$ as in the text follows from \eqref{4.7a} because
  \begin{equation*}
 \int_0^\infty v(r)h(r)\,dr =  \int_0^\infty u(r) h(f(r)) \,dr
   \end{equation*}
and $ h(f(r)) \ge h(r)$ by \eqref{4.6a}.  \qed

\renewcommand{\theequation}{\thesection.\arabic{equation}}
\setcounter{equation}{0}

\vskip2cm

\section{A relaxed notion of partial order}
\label{ch3.3}

\noindent
\begin{defin}[Partial order modulo $m$]
\label{def:2.1}
For any $u$ and $v$ in $\mathcal U$ and  $m>0$, we define  \begin{equation}
\label{4.9a}
u  \preccurlyeq v\;\;\text{ modulo $m$  \; iff }\;\;   F(r;u) \le F(r;v) + m
\;\; \text{for all $r\ge 0$}.
   \end{equation}

\end{defin}

\medskip

\begin{lem}
\label{prop4.2}

Let $u  \preccurlyeq v$ modulo $m$ and $v  \preccurlyeq w$ modulo $m'$,  then
  \begin{equation}
\label{4.10a.0-modulo-m+m'}
u  \preccurlyeq w\quad \text{modulo}\; m+m'.
   \end{equation}

\end{lem}

\medskip

\noindent
{\bf Proof.} $ F(r;u) \le F(r;v) + m \le (F(r;w) + m') +m$.  \qed

%
%
%
%
%
%
%

%
\vskip.5cm

\renewcommand{\theequation}{\thesection.\arabic{equation}}
\setcounter{equation}{0}

\vskip2cm

\section{Inequalities for the cut and the free evolution operators}
\label{?}
In this section by default $\delta>0$, $u$ and $v$ are in $\mathcal U$   and if needed in
 $\mathcal U_\delta$ (as when applying the cut operator $C_\delta$).  We first state and prove the following lemma:
 \vskip.5cm

\begin{lem}
\label{lemma4.2}
Let $u  \preccurlyeq v$ and assume that $m:= F(0;v)-F(0;u)>0$.
Define $\tilde R$  so that $\dis{\int_0^{\tilde R} v= m}$,
then
  \begin{equation}
\label{4.8a}
u  \preccurlyeq v \;{\mathbf 1_{[\tilde R,+\infty)}}=:\tilde v,\quad
F(0;u)=F(0;\tilde v).
   \end{equation}

\end{lem}

\medskip

\noindent
{\bf Proof.}   From  $F(0;u) = F(0;v)-m=F(0,\tilde v)$ we get
	\begin{eqnarray*}
	\int_0^r u(z)dz + F(r;u) = \int_0^r v(z)dz + F(r;v) -m.
	\end{eqnarray*}
Since  $F(r;u)\le F(r;v)$, for all  $r\ge \tilde R$,
	\begin{eqnarray*}
	\int_0^r u(z)dz\ge \int_0^r v(z)dz - m=\int_{\tilde R}^r v(z)dz=\int_0^r \tilde v(z)dz.
	\end{eqnarray*}
Also for $r<\tilde R$ we have
$\dis{\int_0^r u\ge 0= \int_0^r \tilde v}$, so that $\dis{\int_0^r u\ge  \int_0^r \tilde v}$
for all $r$.
Since $F(0;u) = F(0;\tilde v)$ the previous inequality implies that
$F(r;u) \le F(r;\tilde v)$.
\qed

\medskip

 \begin{lem}
 \label{lemma4.5}
Let $u$, $v$ and $\delta$ as above and let $u \preccurlyeq v$. Then 
  \begin{equation}
\label{4.11a}
G_\delta^{\rm neum} * u  \preccurlyeq G_\delta^{\rm neum} *  v, \quad
T_\delta u  \preccurlyeq T_\delta v,
   \end{equation}
  \begin{equation}
\label{4.12a}
C_\delta u  \preccurlyeq u,\quad
C_\delta u  \preccurlyeq C_\delta  v.
   \end{equation}

 \end{lem}

\medskip

\noindent
{\bf Proof.} 
By \eqref{ch2.4.7} the second inequality in \eqref{4.11a} follows from the first one.
 To prove the first one we first observe that since $F(0;u) \le F(0;v)$ there exists $\tilde R_v$ (which may be equal to $0$) so that
 \[
 \int_{\tilde R_v}^\infty  v(r) dr= \int_0^\infty u(r)dr.
 \]
Call  $\tilde v(r)= v(r)\mathbf 1_{r\ge \tilde R_v}$, then by Lemma \ref{lemma4.2}
$  u  \preccurlyeq \tilde v$ and $F(0;u)=F(0;\tilde v)$.  For any $R\ge 0$,
  \begin{equation*}
F(R;G_\delta^{\rm neum}* \tilde v) =\int_0^\infty \tilde v(r')\phi_R(r')dr',\quad
\phi_R(r')= \int_{R}^\infty G_\delta^{\rm neum}(r,r')\,dr
   \end{equation*}
and analogously
     \begin{equation*}
F(R;G_\delta^{\rm neum}* u) =\int_0^\infty u(r')\phi_R(r')dr'.
   \end{equation*}
By an explicit computation: $\dis{ \frac{d}{dr'}\phi_R(r') >0}$, so that
by Corollary \ref{4.2.0.18}
     \begin{equation*}
F(R;G_\delta^{\rm neum}* u)  \le
F(R;G_\delta^{\rm neum}* \tilde v)
   \end{equation*}
   hence the first inequality in \eqref{4.11a} because $\tilde v \le v$.

   The inequality $C_\delta u  \preccurlyeq u$ holds trivially because
$C_\delta u  \le u$.
Furthermore  we have
      \begin{equation*}
C_\delta u-C_\delta  v= (u- v) \;{\mathbf 1_{[0,R_u]}}
- v \, {\mathbf 1_{(R_u, R_v]}}
   \end{equation*}
where $R_u$ is such that $\int_{R_u}^\infty u = j\delta$ and $R_v$ is defined similarly. Hence
      \begin{equation*}
F(r;C_\delta u)-F(r; C_\delta v)\le \Big(F(r;u)-F(r;v)\Big){\mathbf 1_{[0,R_u]}} -{\mathbf 1_{(R_u, R_v]}} \int_r^{R_ v}
 v(r') dr'
   \end{equation*}
which is therefore $\le 0$.    \qed
%
%
%

\medskip

 \begin{lem}
 \label{nuovo.lemma5.7}
Let  $u  \preccurlyeq v$ modulo $m$,  then
  \begin{equation}
\label{4.13a}
T_\delta u  \preccurlyeq T_\delta v \;\; \text{modulo }\;\; m,\quad
G^{\rm neum}_t *u   \preccurlyeq  G^{\rm neum}_t *v  \;\; \text{modulo }\;\; m.
   \end{equation}
 \end{lem}

\medskip

\noindent
{\bf Proof.}   By \eqref{ch2.4.7} we just need to prove the
 second inequality which obviously holds if $F(0;u) \le m$.  We thus suppose
$F(0;u) > m$ and define
 \begin{equation}
\label{4.10a.0-erre-m}
u^* := u \;{\mathbf 1_{[0,R_m]}} \;\;\;\;\text{with}\;\; R_m :
\int_{R_m}^\infty u  = m.
   \end{equation}
We are going to show that
  \begin{equation}
\label{4.10a}
u^*   \preccurlyeq v. 
   \end{equation}
In fact $F(r;u^*)= F(r;u)-m$ when $r\le R_m$ so that $F(r;u^*) \le F(r;v)$.
 Since $F(r;u^*)=0$ for $r\ge R_m$ then $F(r;u^*) \le F(r;v)$
 so that \eqref{4.10a} is proved.
\noindent
By \eqref{4.11a}
\begin{eqnarray*}
F(r;G^{\rm neum}_\delta*u) &=& F(r;G^{\rm neum}_\delta *u^*) + F(r;G^{\rm neum}_\delta *(u-u^*))
\\&\le& F(r;G^{\rm neum}_\delta *v) + F(0;G^{\rm neum}_\delta *(u-u^*)) \\&=& F(r;G^{\rm neum}_\delta *v) +m.
\end{eqnarray*}
\qed

\medskip

 \begin{lem}
 \label{nuovo.lemma5.8}
 Let   $u \preccurlyeq v$ modulo $m$,  then
  \begin{equation}
\label{4.14a}
 u  \preccurlyeq C_\delta  v\;\;\; \text{modulo }\;\; m+j\delta.
   \end{equation}

 \end{lem}
 \medskip

\noindent
{\bf Proof.}
%
We have $F(r;u) \le F(r;v)+m$ and  $F(r;v) \le F(r;C_\delta v) + j\delta$
hence
\[
F(r;u) \le F(r;C_\delta v) + m +j\delta
\]
which proves \eqref{4.14a}.
\qed

 \begin{lem}
  \label{nuovo.lemma5.9}
 Let  $u  \preccurlyeq v$ modulo $m$, $m \ge j\delta$, then
  \begin{equation}
\label{4.15a}
 C_\delta  u  \preccurlyeq  v\;\;\; \text{modulo}\;\;\;  m -j\delta.
   \end{equation}

 \end{lem}

 \medskip

\noindent
{\bf Proof.} By the definition of $C_\delta $ for
$r\le R_u$ we have
\[
F(r;C_\delta u) = F(r;u) - j\delta \le F(r;v)+m -j\delta.
\]
If instead $r\ge R_u$,
\[
F(r;C_\delta u) = 0 \le  m -j\delta \le F(r;v)+m -j\delta
\]
hence \eqref{4.15a}.
\qed

\medskip

 \begin{lem}
 \label{nuovo.lemma5.9.0}
  Let $u,v$ in $\mathcal U_\delta$,    $u  \preccurlyeq v$ modulo $m$,  then
  \begin{equation}
\label{4.14ak}
 C_\delta u  \preccurlyeq C_\delta  v\;\;\; \text{modulo }\;\; m.
   \end{equation}

 \end{lem}
 \medskip

\noindent
{\bf Proof.}  By Lemma \ref{nuovo.lemma5.8}
  \begin{equation*}
 u  \preccurlyeq   w:= C_\delta v\;\;\; \text{modulo }\;\; m + j\delta.
   \end{equation*}
By Lemma  \ref{nuovo.lemma5.9}
 \begin{equation*}
C_\delta  u   \preccurlyeq   w\;\;\; \text{modulo }\;\; m
   \end{equation*}
hence \eqref{4.14ak}.  \qed

%
%
%
%

\renewcommand{\theequation}{\thesection.\arabic{equation}}
\setcounter{equation}{0}

\vskip1cm

\section{Inequalities for the barriers}
\label{IFB}
The following theorems are consequence of the inequalities established in the previous section.

\medskip

 \begin{thm}
 \label{IFBthm.0}

 Let $\delta>0$, $u,v\in \mathcal U_\delta$, $u  \preccurlyeq v$ modulo $m\ge 0$. Let $k\in \mathbb N$,  then
  \begin{equation}
\label{IFB.0}
 S^{\delta,\pm}_{k\delta}  u  \preccurlyeq   S^{\delta,\pm}_{k\delta}  v\quad
 \text{\rm modulo}\; m.
   \end{equation}

 \end{thm}

  \medskip

\noindent
{\bf Proof.}  It follows from Lemma \ref{nuovo.lemma5.7} and Lemma \ref{nuovo.lemma5.9.0}.  \qed

\medskip

 \begin{thm}
 \label{IFBthm.1}

 Let $\delta>0$, $u\in \mathcal U_\delta$ and $k\in \mathbb N$,  then
  \begin{equation}
\label{IFB.1}
 S^{\delta,-}_{k\delta}  u  \preccurlyeq   S^{\delta,+}_{k\delta}  u.
   \end{equation}

 \end{thm}

  \medskip

\noindent
{\bf Proof.}  We proceed as in the proof of Theorem \ref{CDGPmthm5.3} and write
 \[
S^{\delta,-}_{k\delta}  u= C_{\delta}T_\delta \cdots  C_{\delta}T_\delta u,\qquad S^{\delta,+}_{k\delta}  u=T_\delta  C_{\delta}\cdots T_\delta C_{\delta} u \qquad {\text{$k$ times}}.
 \]
Call 
 \[
 v = C_{\delta}u,\qquad v_k= T_\delta C_{\delta}\cdots T_\delta v,\qquad u_k=T_\delta C_{\delta}\cdots T_\delta u \qquad {\text{$k$ times}}
 \]
so that $S^{\delta,+}_{k\delta}  u=v_k$ and $S^{\delta,-}_{k\delta}  u=
C_\delta u_k$.  By \eqref{4.14a} with $u=v$,
  \begin{equation*}
 u  \preccurlyeq C_\delta  u\;\;\; \text{\rm modulo }\;\; j\delta.
   \end{equation*}
By Lemmas \ref{nuovo.lemma5.7}--\ref{nuovo.lemma5.8}
 \begin{equation*}
 u_k  \preccurlyeq v_k\;\;\; \text{\rm modulo }\;\; j\delta
   \end{equation*}
and by Lemma \ref{nuovo.lemma5.9}
 \begin{equation*}
S^{\delta,-}_{k\delta}  u=
C_\delta u_k  \preccurlyeq v_k=S^{\delta,+}_{k\delta}  u.
   \end{equation*}
\qed

\medskip

 \begin{thm}
 \label{IFBthm.2}

Let $\delta'>0$, $\delta= n \delta'$, $n \in \mathbb N$, $u \in \mathcal U_\delta$ and $t=k\delta$, then
  \begin{equation}
\label{IFB.2}
 S^{\delta,-}_{t}  u  \preccurlyeq   S^{\delta',-}_{t}  u.
   \end{equation}

 \end{thm}

  \medskip

\noindent
{\bf Proof.} We postpone the proof that
 \begin{equation}
\label{IFB.3}
S_{\delta}^{\delta,-}u  \preccurlyeq
S_{\delta}^{\delta',-}u'\qquad \text{if
$u  \preccurlyeq u'$}.
   \end{equation}
By \eqref{IFB.3} with $u'=u$ we get $S_{\delta}^{\delta,-}u  \preccurlyeq
S_{\delta}^{\delta',-}u$ so that using again \eqref{IFB.3}
 \begin{equation*}
S_{\delta}^{\delta,-}(S_{\delta}^{\delta,-}u)  \preccurlyeq
S_{\delta}^{\delta',-}(S_{\delta}^{\delta',-}u),\quad
S_{2\delta}^{\delta,-}u  \preccurlyeq S_{2\delta}^{\delta',-}u
   \end{equation*}
which by iteration proves  \eqref{IFB.2}.

Proof of \eqref{IFB.3}.
We have
\begin{eqnarray*}
&&S_{\delta}^{\delta',-}u' =S_{n\delta'}^{\delta',-}u' = C_{\delta'} T_{\delta'}\cdots C_{\delta'} T_{\delta'}u' \qquad \text{$n$ times},
\\
&& S_{\delta}^{\delta,-}u = C_{\delta}T_\delta u
=  C_{\delta'}^nT_{\delta'}^nu.
\end{eqnarray*}
We will prove by induction on $k$ that
\begin{equation}
\label{IFB.4}
 C_{\delta'}^kT_{\delta'}^ku  \preccurlyeq
  C_{\delta'} T_{\delta'}\cdots C_{\delta'} T_{\delta'}u \qquad \text{$k$ times}.
\end{equation}
\eqref{IFB.3} will then follow by setting  $k=n$
and using Lemma \ref{lemma4.5}.  We thus suppose that \eqref{IFB.4} holds with $k$ and want to prove that it holds for $k+1$. We preliminarily show that for any integer $h>0$,
\begin{equation}
\label{IFB.5}
C_{\delta'}^hT_{\delta'}v \preccurlyeq T_{\delta'}C_{\delta'}^h v.
\end{equation}
In fact by \eqref{4.14a}
$
v \preccurlyeq C_{\delta'}^h v$ modulo $jh\delta'$.
Then by \eqref{4.11a},
\[
T_{\delta'}v \preccurlyeq T_{\delta'}C_{\delta'}^h v \;\text{modulo}\; jh\delta'
\]
and  by \eqref{4.15a} $C_{\delta'}^hT_{\delta'}v \preccurlyeq T_{\delta'}C_{\delta'}^h v
$.
\eqref{IFB.5} is proved.

Call $v = T_{\delta'}^ku$ then
using \eqref{4.12a} and \eqref{IFB.5},
\begin{equation*}
 C_{\delta'}^{k+1}T_{\delta'}^{k+1}u =  C_{\delta'} C_{\delta'}^{k}  T_{\delta'}
v  \preccurlyeq C_{\delta'} T_{\delta'}  C_{\delta'}^{k} v.
\end{equation*}
By assumption \eqref{IFB.4} holds with $k$ so that  calling $w$ its right-hand side,
$ C_{\delta'}^{k} v  \preccurlyeq w$.  By Lemma \ref{lemma4.5},
\begin{equation*}
  C_{\delta'} T_{\delta'}  C_{\delta'}^{k} v \preccurlyeq  C_{\delta'} T_{\delta'}  w
\end{equation*}
which is \eqref{IFB.4}   with $k+1$.
   \qed

\medskip

\begin{thm}
\label{thm5.2.0.29}
Let  $u  \in  \mathcal U$,  $\delta = h \delta'$, $h$ a positive integer and $t=k\delta$, then
		\begin{equation}
		\label{nna.2}
S_{t}^{\delta',+}u  \preccurlyeq
S_{t}^{\delta,+}u.
			\end{equation}

\end{thm}
\medskip
\noindent
{\bf Proof.}
We postpone the proof that
 \begin{equation}
\label{4.16a}
S_{\delta}^{\delta',+}u  \preccurlyeq
S_{\delta}^{\delta,+}u'\qquad \text{if
$u  \preccurlyeq u'$}.
   \end{equation}
By \eqref{4.16a} with $u'=u$ we get $S_{\delta}^{\delta',+}u  \preccurlyeq
S_{\delta}^{\delta,+}u$ so that, using again \eqref{4.16a},
 \begin{equation*}
S_{\delta}^{\delta',+}(S_{\delta}^{\delta',+}u)  \preccurlyeq
S_{\delta}^{\delta,+}(S_{\delta}^{\delta,+}u),\quad
S_{2\delta}^{\delta',+}u  \preccurlyeq S_{2\delta}^{\delta,+}u
   \end{equation*}
which by iteration yields \eqref{nna.2}.
\noindent
Proof of \eqref{4.16a}.
We have
\[
S_{\delta}^{\delta',+}u =S_{h\delta'}^{\delta',+}u = T_{\delta'} C_{\delta'} \cdots T_{\delta'} C_{\delta'} u \qquad \text{$h$ times},
\]
\[
S_{\delta}^{\delta,+}u' = T_\delta C_{\delta}u
=  T_{\delta'}^h    C_{\delta}u'. 
\]
By \eqref{4.14a} $u' \preccurlyeq  C_\delta u'$  modulo $j\delta$
so that, by Lemma \ref{prop4.2}, $u \preccurlyeq C_\delta u'$ modulo $j\delta$.
By \eqref{4.15a}
 $C_{\delta'} u  \preccurlyeq
C_\delta u'$ modulo  $j\delta-j\delta'$. By \eqref{4.13a},
\begin{eqnarray*}
T_{\delta'} C_{\delta'} u \preccurlyeq
T_{\delta'}C_\delta u' \quad  \text{modulo } \: j\delta-j\delta'.
\end{eqnarray*}
Call $w:=T_{\delta'} C_{\delta'} u$ and $v':=T_{\delta'}C_\delta u'$, then
$w \preccurlyeq  v'$ modulo  $j\delta-j\delta'$.
By   \eqref{4.15a}  $C_{\delta'} w \preccurlyeq  v'$ modulo $j\delta-2 j\delta'$ and by  \eqref{4.13a}
\begin{eqnarray*}
T_{\delta'} C_{\delta'} w \preccurlyeq
T_{\delta'} v' \quad  \text{modulo } \: j\delta-2 j\delta'.
\end{eqnarray*}
 \eqref{4.16a} then follows by iteration.
   \qed

\medskip
\noindent
For general $\delta$ and $\delta'$ we will use the following bound:
\medskip

   \begin{lem}
  \label{lemmae35}
There is $c$ so that for any  $0<\delta<\delta'$,
$u\in \mathcal U_\delta$  and $n\ge 1$,
  \begin{equation}
    \label{e6.20.0}
| S^{\delta,+}_{n\delta}u -S^{\delta',+}_{n\delta'}u |_1 \le c |u|_1 n\frac{\delta'-\delta}{\delta^{3/2}}.
   \end{equation}

    \end{lem}

\medskip

\noindent
{\bf Proof.}
By
\eqref{CDGP3.1.4} for any $u,v \in \mathcal U$, $|C_\delta u - C_\delta v|_1 \le | u -  v|_1$.
We also have
  \begin{equation}
    \label{e6.20.1}
| C_\delta w -C_{\delta'}w |_1 \le j(\delta'-\delta),\quad
|G^{\rm neum}_\delta*w - G^{\rm neum}_{\delta'}*w|_1 \le \frac{c(\delta'-\delta)}{\delta^{3/2} } |w|_1.
   \end{equation}
Recalling that $S^{\delta,+}_{\delta}u =G^{\rm neum}_{\delta}* C_\delta u$, we get that
$| S^{\delta,+}_{\delta}w -S^{\delta',+}_{\delta'}v |_1$ is bounded by
  \begin{eqnarray*}
&& \le |G_\delta^{\rm neum}*C_{\delta'} w +G_\delta^{\rm neum}*(C_\delta-C_{\delta'})w
-G_{\delta'}^{\rm neum}*C_{\delta'} w -G_{\delta'}^{\rm neum}*C_{\delta'}(v-w)|,
%
%
%
%
   \end{eqnarray*}
hence
 \begin{equation}
     \label{e6.20.1.1}
     | S^{\delta,+}_{\delta}w -S^{\delta',+}_{\delta'}v |_1
\le |w-v|_1
+ c\frac{\delta'-\delta}{\delta^{3/2} } |w|_1  + j (\delta'-\delta).
   \end{equation}
We use \eqref{e6.20.1.1} to prove \eqref{e6.20.0} by induction on $n$.
We prove \eqref{e6.20.0} when $n=1$ by setting $w=v=u$ in \eqref{e6.20.1.1}.
We suppose by induction that \eqref{e6.20.0} hods till $n-1$. We then
set $w=S^{\delta,+}_{(n-1)\delta}u$ and
$v=S^{\delta',+}_{(n-1)\delta'}u$ in  \eqref{e6.20.1.1} getting
\eqref{e6.20.0}.  \qed


\renewcommand{\theequation}{\thesection.\arabic{equation}}
\setcounter{equation}{0}

\chapter{The limit theorems on barriers}
\label{ch:4}

In this chapter we will prove Theorem \ref{CDGPthm2.4.1}. An analogous theorem is proved
in \cite{CDGP}  when $G_t^{\rm neum}$ is
replaced by the Green function with Neumann condition both at 0 and at 1.

\renewcommand{\theequation}{\thesection.\arabic{equation}}
\setcounter{equation}{0}

\vskip2cm

\section{The limit function $\psi$}
\label{ottava}

In this section we define a function  $\psi(r,t)$
which in the next section will be proved to be the function $S_t u(r)$  of Theorem \ref{CDGPthm2.4.1}.
We
fix  $T>0$, $u\in  \mathcal U$, $\tau>0$ and $t_0>0$, call
$\Delta_\tau:= \{2^{-n}\tau, n \in  \mathbb N\}$,
$\mathcal T_{\tau,n}=\{t= k2^{-n}\tau,  k\in \mathbb N\}$ and
$\mathcal T_\tau=\{t= k2^{-n}\tau, n \in  \mathbb N, k\in \mathbb N\}$.

\medskip

\subsection{Convergence of the upper barriers}
\label{zsubsec6.2.2}

In Chapter \ref{ch:2bis} we  proved that the family of upper barriers is equi-bounded
and equi-continuous so that it converges by subsequences.  In this subsection we will prove convergence, see \eqref{corr8.1} below.  More precisely we   restrict to
$\delta\in \Delta_\tau$ and
define a function $\psi^{(n)}(r,t)$ on $\mathbb R_+\times [t_0,T]$
by first setting
\[
\psi^{(n)}(r,t) = S^{2^{-n}\tau,+}_t u(r),\quad r\in \mathbb R_+,\;\; t \in [t_0,T] \cap \mathcal T_{\tau,n}
\]
and then extending $\psi^{(n)}(r,t)$ to $t\in [t_0,T]$  by linear interpolation.
As mentioned above we have proved
in Chapter \ref{ch:2bis} that the family $\{\psi^{(n)}\}$ is equi-bounded
and equi-continuous so that, by the Ascoli-Arzel\`a theorem, it converges by subsequences in
sup norm (on the compacts) to a continuous function $\psi(r,t)$  on $\mathbb R_+\times [t_0,T]$.  To prove full convergence we will show that $\{F(r;\psi^{(n)}), r \ge 0\}$ converges, the proof will follow from the monotonicity properties of the barriers and the following a priori bound:

\medskip

\begin{lem}
There is $c>0$ so that
	\begin{equation}
	\label{e6.2.00}
S^{\delta,+}_{t} u(r) \le \int dr'u(r') G_{0,t}^{\rm neum}(r',r) + c j \sqrt t e^{-r^2/(4t)}.
	\end{equation}
\end{lem}

\medskip
\noindent
{\bf Proof.}  It follows from
 \eqref{e6.2} with $m = 0$ and bounding for $s \le t$,
 \[
 G_{0,s}(0,r) \le e^{-r^2/(4t)} \Big( ({2\pi s})^{-1/2}
 e^{-r^2/(4s)}\Big).
 \]
\qed

\medskip

 \eqref{e6.2.00}   guarantees
convergence of   $F(r;\psi^{(n)}(\cdot,t))$.

 \medskip

\medskip
\begin{lem}
Let $\psi$ be any limit point of $\{\psi^{(n)}\}$, then
for any $r\in \mathbb R$ and $t\in [t_0,T]\cap  \mathcal T_\tau$,
	\begin{equation}
	\label{corr8.1.12}
\lim_{n\to \infty} F(r; S^{2^{-n}\tau,+}_tu) = F(r;\psi(\cdot,t)).
 	\end{equation}
As a consequence there is a unique  limit point $\psi$
of $\{\psi^{(n)}\}$ and for any $n$ and $t \in  \mathcal T_{\tau,n}$,
	\begin{equation}
	\label{corr8.1.11}
F(r;S^{2^{-n}\tau,+}_tu) \ge F(r; \psi(\cdot,t)).
 	\end{equation}
 Moreover
 	\begin{equation}
	\label{corr8.1.11.1}
\psi(\cdot,t) \le  \int dr'u(r') G_{0,t}^{\rm neum}(r',r) + c j \sqrt t e^{-r^2/(4t)}.
 	\end{equation}
\end{lem}

\medskip
\noindent
{\bf Proof.}  By Theorem \ref{thm5.2.0.29}
$F(r; S^{2^{-n}\tau,+}_tu)$, $t\in\mathcal T_\tau$, is a non-increasing function of $n$ hence the existence of the limit $n\to \infty$.  To identify the limit we observe that
the right-hand side of \eqref{e6.2.00}
is for each $t\le T$ an $L^1$ function of $r$.
\eqref{corr8.1.12} then
follows using the Lebesgue dominated convergence theorem.
Thus all limit functions
$\psi(r,t)$ agree on $t\in [t_0,T]\cap  \mathcal T_\tau$
and since they are continuous they agree on the whole $[t_0,T]$,
thus the sequence $\psi^{(n)}(r,t)$ converges
in sup-norm  as $n\to \infty$ to a continuous function
$\psi(r,t)$ (and not only by subsequences).

  \eqref{corr8.1.11} follows from \eqref{corr8.1.12} because $F(r; S^{2^{-n}\tau,+}_tu)$ is a non-increasing and  \eqref{corr8.1.11.1} follows from \eqref{e6.2.00} because we have already proved that $S^{2^{-n}\tau,+}_tu$ converges to $\psi(r,t)$.
\qed

\vskip.5cm

By the arbitrariness of $t_0$ and $T$ the function $\psi(r,t)$ extends to
the whole $\mathbb R_+\times (0,\infty)$.
 Thus,   by \eqref{e6.2.00},
	\begin{equation}
	\label{corr8.1}
\lim_{n\to \infty} \|S^{2^{-n}\tau,+}_tu - \psi(\cdot,t)\|_\infty  = 0,\quad
t>0, t\in   \mathcal T_\tau,
 	\end{equation}
the convergence
being uniform in $t\in  \mathcal T_\tau$ when it varies
on the compacts not containing 0.

The drawback of this result is that the function $\psi$ we have defined actually depends on
$\tau$, to underline this we will write it as $\psi_\tau(r,t)$.  We will prove in the next subsection that all $\psi_\tau(r,t)$ are identical to each other.

\vskip.5cm

\subsection{Independence of $\tau$}
\label{zsubsec6.2.3}

   \begin{thm}
  \label{thme35}
$\psi_ \tau$ is independent of $\tau$.

    \end{thm}

\medskip

\noindent
{\bf Proof.}  It  suffices to prove
that for any $\tau$ and $\tau'$,
\[
F(r;\psi_\tau(\cdot,t)) = F(r;\psi_{\tau'}(\cdot, t)),\quad r \ge 0, \;\;t > 0
\]
as $\psi_{\tau'}$ and $\psi_{\tau}$ are continuous.  We suppose that
$\tau' \notin \{ k\tau 2^{-n}, k, n \in \mathbb N\}$ (the case when they are rationally related is proved using Theorem \ref{thm5.2.0.29}).  We fix $t'= n \delta'$, $\delta'= \tau' 2^{-m}$.
Let $\delta = k \tau 2^{-q}$, $\delta <\delta'$.
By Lemma \ref{lemmae35}, for all $r \ge 0$,
 \begin{equation*}
F(r; S^{\delta',+}_{t'}u )\ge F(r;S^{\delta,+}_{n\delta}u) -
c |u|_1 n\frac{\delta'-\delta}{\delta^{3/2} }.
   \end{equation*}
Write $\delta= k_p \tau 2^{-p}$ so that $k_p = k 2^{p-q}$ is a positive integer for $p$ large enough.
Then by Theorem \ref{thm5.2.0.29},
 \begin{equation*}
 F(r;S^{\delta,+}_{n\delta}u )\ge  F(r;S^{\tau 2^{-p},+}_{n\delta}u).
   \end{equation*}
By taking $p\to \infty$:
 \begin{equation*}
F(r; S^{\delta',+}_{t'}u )\ge F(r;\psi_{\tau}(\cdot,n\delta)) -
c |u|_1 n\frac{\delta'-\delta}{\delta^{3/2} }.
   \end{equation*}
We then let $\delta \to \delta'$ on $\{k\tau 2^{-n}, k, n \in \mathbb N\}$.  In this limit
$n\delta \to t'$ and by the continuity of $\psi_\tau(\cdot, s)$ in $s$ we get
 \begin{equation*}
F(r; S^{\delta',+}_{t'}u )\ge F(r;\psi_\tau(\cdot, t')).
   \end{equation*}
We next take $m\to \infty$, recall  $\delta'= \tau' 2^{-m}$, and get
 \begin{equation*}
F(r; \psi_{\tau'}(\cdot,t'')) \ge F(r;\psi_\tau(\cdot, t'')),\qquad \text{for any $t''\in \{k\tau' 2^{-n}, k, n \in \mathbb N\}$}.
   \end{equation*}
In an analogous fashion we get
 \begin{equation*}
F(r;\psi_\tau(\cdot, t))\ge F(r; \psi_{\tau'}(\cdot,t)) ,\qquad \text{for any $t\in \{k\tau 2^{-n}, k, n \in \mathbb N\}$}.
   \end{equation*}
Then $\psi_\tau(\cdot,t)= \psi_{\tau'}(\cdot,t)$ for all $t$ in a dense set, hence they are equal everywhere
being both continuous.\qed

\medskip

We can thus drop $\tau$ and simply write $\psi(r,t)$.  We can then summarize:

   \begin{coro}
  \label{coronuovo35}
There is a continuous function $\psi(r,t)$ on $[0,\infty)\times \mathbb R_+$
which satisfies the bound \eqref{corr8.1.11.1} and such that for any $\tau>0$
$S^{2^{-n}\tau,+}_t u(r)$ converges to $\psi(r,t)$ on the compacts.

    \end{coro}

\vskip.5cm

\subsection{Continuity at 0}

\begin{prop}
\label{propz6.3}
Let  $u\in \mathcal U$, then $\psi(\cdot,t)$ converges weakly to $u$ as $t\to 0$ and
	\begin{equation}
	\label{corr8.2}
\lim_{t\to 0} F(r;\psi(\cdot,t))  = F(r;u).
 	\end{equation}
Suppose further that $u$ is a continuous function with compact support. Then
	\begin{equation}
	\label{corr8.2.00}
\lim_{t\to 0} \|\psi(\cdot,t)   -u\|_\infty = 0.
 	\end{equation}

\end{prop}

\medskip

\noindent
{\bf Proof.}  Let $t\in \mathcal T_\tau$, then
by \eqref{e6.333333}--\eqref{e6.3.0} with $s=0$ we have
	\begin{equation}
	\label{e6.333333.00.0}
|S^{ 2^{-n}\tau,+}_{t}u (r) - \int dr' u(r') G_{t}^{\rm neum}(r',r)| \le cj \sqrt t.
	\end{equation}
By \eqref{corr8.1}, letting $n\to \infty$,
	\begin{equation}
	\label{e6.333333.00.1}
|\psi(r,t) - \int dr' u(r') G_{t}^{\rm neum}(r',r)| \le cj \sqrt t.
	\end{equation}
Since $\psi(r,t)$ is continuous in $t$, \eqref{e6.333333.00.1} holds for all $t>0$.
$G_{t}^{\rm neum}*u$ converges weakly to $u$ as $t\to 0$ , hence also $\psi(r,t)$ converge weakly to $u$ as $t\to 0$.  Analogously, since
	\begin{equation*}
\lim_{t \to 0}  \int_r^R dr'' \int dr' u(r') G_{t}^{\rm neum}(r',r'') =  \int_r^R dr'' u(r'')
	\end{equation*}
then by \eqref{e6.333333.00.1}
	\begin{equation*}
\lim_{t \to 0}  \int_r^R dr'' \psi(r'',t)  =  \int_r^R dr'' u(r'').
	\end{equation*}
Hence
	\begin{equation*}
\liminf_{t \to 0}  \int_r^\infty dr'' \psi(r'',t)  \ge  \int_r^R dr'' u(r'')
	\end{equation*}
and by the arbitrariness of $R$ it is $\ge \int_r^\infty dr'' u(r'')$.  To prove the upper bound we use \eqref{corr8.1.11.1} to say that
for any $\eps>0$ there is $R_\eps$ so that for all $t\le 1$
$\int_ {R_\eps}^\infty  dr'' \psi(r'',t) \le \eps$ and
$\int_ {R_\eps}^\infty  dr'' u(r'') \le \eps$ as well.  Then
	\begin{equation*}
\limsup_{t \to 0}  \int_r^\infty dr'' \psi(r'',t)  \le \eps+
\limsup_{t \to 0}
 \int_r^{R_\eps} dr'' \psi(r'',t) \le \eps +
 \int_r^{\infty} dr'' u(r'').
	\end{equation*}
Thus \eqref{corr8.2} is proved.
By \eqref{e6.333333}
	\begin{equation}
	\label{e6.333333.00}
S^{\delta,+}_{t}u (r) -u(r)= \int dr' [u(r') -u(r)]G_{t}^{\rm neum}(r',r) +
v_{0,t}^{(\delta,+)} (r).
	\end{equation}
Hence by \eqref{e6.3.0} there is a function $\eps(t)$ which vanishes as $t\to 0$ such that
	\begin{equation}
	\label{e6.333333.01}
\|S^{\delta,+}_{t}u   -u\|_\infty \le \eps(t).
	\end{equation}
By \eqref{corr8.1}
	\begin{equation}
	\label{e6.333333.02}
\|\psi(\cdot,t)   -u\|_\infty \le \eps(t)
	\end{equation}
so that \eqref{corr8.2.00} is proved.
\qed

%
%
%

\renewcommand{\theequation}{\thesection.\arabic{equation}}
\setcounter{equation}{0}

\vskip2cm

\section{Proof of Theorem \ref{CDGPthm2.4.1}}
\label{ottavabis}

We shall now prove Theorem
 \ref{CDGPthm2.4.1}. The items below correspond to the items in
  Theorem
 \ref{CDGPthm2.4.1}.

 \medskip

 \begin{itemize}
 \item  (1)  is proved in Theorem \ref{IFBthm.1}.

  \item  (2) is   proved in  \eqref{4.3a.00}.

  \item  (3) with the + is proved in Theorem \ref{thm5.2.0.29} with the - in Theorem
  \ref{IFBthm.2}.

  \item (4) is proved in Corollary \ref{coronuovo35}.  Convergence in $L^1$
  follows from the convergence on the compacts and the uniform bound \eqref{e6.2.00}.


  \item (5) with the + is proved in \eqref{corr8.1.11}.  Monotonicity with the - has
  already been proved, see (3), and by (2) which has also been proved, the limit with the - is the same as with the +.

   \item (6) follows from \eqref{IFB.1} and (3)  which has also been proved.

   \item (7) is proved in Proposition \ref{propz6.3}.

    \item (8) follows by (5) and (1) which have   been already proved.

    \item (9) follows from Proposition \ref{propZ4.1.0.18} and (4)
    which has been already proved.

 \end{itemize}

\vskip.5cm

\renewcommand{\theequation}{\thesection.\arabic{equation}}
\setcounter{equation}{0}

\renewcommand{\theequation}{\thesection.\arabic{equation}}
\setcounter{equation}{0}

\chapter{Brownian motion and the heat equation}
\label{ch:5}
%
The proof of Theorem \ref{CDGPthm2.3.1}
uses extensively a representation of the solution of the heat equation
in terms of Brownian motions.  We will recall in this chapter the main properties and in particular we re-derive
a formula, \eqref{3.144444} below, for the
solution $\rho(r,t)$ of  \eqref{intro.2}--\eqref{intro.3} with initial datum   $\rho(r',s)$ at time $s$ in terms of Brownian motions.  We will write the Green function for \eqref{intro.2}--\eqref{intro.3} in terms of the first
exit time  distribution of a Brownian motion, \eqref{5.3.1bb}, and then relate  the exit time distribution density to the
derivative of the solution of the heat equation at the edge.  The latter gives
the rate of mass which is dissipated because of the Dirichlet boundary conditions thus the the mass loss is directly related to the exit probability of the Brownian motion.

By default in the sequel $X=\big(X_t, t\in [0,T]\big)$, is
a positive continuous  function piecewise $C^1$ and with right and left derivatives at all times.


\medskip

\renewcommand{\theequation}{\thesection.\arabic{equation}}
\setcounter{equation}{0}

\medskip

\section{Brownian motion on the line}
\label{ch:5.1}

We start from the heat equation on the whole $\mathbb R$.  We call
$Q_{r,s}$, $r\in \mathbb R$, $s\ge 0$, the law on $C(\mathbb R,[s,\infty))$ of the
Brownian motion $B_t$, $t\ge s$, which starts from $r$ at time $s$, i.e.\  $B_s=r$.
For each $t> s$ the law of $B_t$ is absolutely continuous with respect to the
Lebesgue measure and has a probability density  $G_{s,t}(r',r)$ which is the Gaussian $G_{t-s}(r',r)$ defined in \eqref{ch2.4.6}. Thus
	\begin{equation}
	\label{3.1.00.0}
E_{Q_{r',s}}[f(B_t)]=\int  G_{s,t}(r',r)f(r)dr,\quad f\in L^\infty(\mathbb R).
	\end{equation}
We can read \eqref{3.1.00.0} by saying that we
start a Brownian motion from $r'$ at time $s$ and run
it till time $t$.  We then compute $f$ at the final point and integrate over all samples: this is the same as integrating $f$ with the Green function $G_{s,t}(r',r)$.

Since $G_{s,t}(r',r)$ as a function of $(r,t)$ solves the heat equation for $t>s$, see \eqref{ch2.4.7.1-proof}, then by differentiating \eqref{3.1.00.0} with respect to $t$ we get
\begin{equation}
\label{3.1.00.1}
\frac{d}{dt}E_{Q_{r',s}}[f(B_t)]= \frac 12 E_{Q_{{r'},s}}[f''(B_t)],\quad f\in C^2(\mathbb R).
\end{equation}
By \eqref{ch2.4.7.1-proof}--\eqref{ch2.4.7.2} if $\rho(r',s)$ is a continuous function of $r'$, then
\begin{equation}
\label{3.1.00}
\rho(r,t) := \int \rho(r',s)G_{{s},t}(r',r) \,dr'
\end{equation}
solves the heat equation in $\mathbb R$ with datum $\rho(r',s)$ at time $s$
hence analogously to \eqref{3.1.00.0}
\begin{equation}
\label{3.1.00.9}
\int f(r)\rho(r,t)dr := \int \rho(r',s)E_{Q_{r',s}}[f(B_t)]  \,dr',\quad f\in L^\infty(\mathbb R).
\end{equation}

\medskip

\section{Reflected Brownian motion with mass injection}
\label{ch:5.2-cap}

We denote by
$P_{r,s}$, $r\ge 0$, $s\ge 0$,  the probability law on the space $C(\mathbb R_+,[s,\infty))$
of the
Brownian motion $B_t$, $t\ge s$, which starts from $r$ at time $s$, i.e.\  $B_s=r$, and which is
reflected at $0$,  $E_{r,s}$ denoting its expectation.
$P_{r,s}$ may be defined as the law of $|B_t|, t\ge s$, under $Q_{r,s}$.
Thus
for any $f\in L^\infty(\mathbb R_+)$,
\begin{equation}
\label{3.1.01}
E_{r',s}[f(B_t)] =  E_{Q_{r',s}}[ f(B_t) \mathbf 1_{B_t\ge 0}  +
f(-B_t) \mathbf 1_{B_t < 0}].
\end{equation}
Hence
\begin{equation}
\label{3.1.02}
G^{\rm neum}_{s,t}(r',r) := G_{s,t}(r',r)+G_{s,t}(r',-r) =
G_{s,t}(r',r)+G_{s,t}(-r',r)
\end{equation}
is the Lebesgue density of the law of the reflected Brownian motion $B_t$:
\begin{equation}
\label{3.1.02.00}
E_{{r'},s}[f(B_t)]=\int_{\mathbb R_+} G^{{\rm neum}}_{s,t}(r',r)f(r)dr,\quad f\in L^\infty(\mathbb R_+).
\end{equation}
Thus by Proposition \ref{CDGPprop2.4.1} the Lebesgue density of the law of the reflected Brownian motion $B_t$ solves the heat equation with Neumann conditions at 0 and
 if $ \rho(r',s)$ is a continuous function of $r'$,
	\begin{equation}
\label{3.1}
\rho(r,t) := \int_{\mathbb R_+}  \rho(r',s) G^{{\rm neum}}_{{s},t}(r',r) \,dr'
	\end{equation}
solves the heat  equation in $\mathbb R_+$ with
Neumann conditions at 0 and initial  datum $\rho(r',s)$ at time $s$.
Moreover if $\phi(r)$, $r\ge 0$, has a $C^2$ symmetric extension to $\mathbb R$,
i.e.\ $r \to \psi(r)= \phi(|r|)$
is   $C^2(\mathbb R)$, then by \eqref{3.1.01} and \eqref{3.1.00.1},
	\begin{equation}
\label{3.1.00.000}
\frac{d}{dt} E_{{r'},s}[\phi(B_t)] = \frac 12 E_{{r'},s}[\phi''(B_t)].
\end{equation}
\noindent
The operator in \eqref{ch2.4.7} can be written as
	\begin{equation}
\label{3.1a}
T_t\rho (r)dr= \int_{\mathbb R_+}  \rho(r') P_{r',0}(B_t=dr) \,dr' +j\int_0^t P_{0,s}(B_t=dr)ds.
	\end{equation}

\medskip

\section{Brownian motion with reflection at 0 and absorption at the edge}
\label{ch:5.3}
Let $B_t$, $t\ge s \ge 0$,  be the Brownian motion starting at $s$
from $r$ and with reflections at $0$, $P_{r,s}$ its law.
Recall that   $X=( X_t$, $t\ge 0,)$  is
a positive continuous  function piecewise $C^1$ and with bounded left and right derivatives at all $t$. Given
$s\ge 0$, we define
 \begin{equation}
\label{5.2.1}
 \tau^X_s = \inf\{t\ge s:  B_t \ge X_t\},\;\;\;\text{and $=\infty$ if the set is empty}
\end{equation}
and denote by $F^X_{r,s}(ds')$
the probability distribution  of  $\tau^X_s$ induced by $P_{r,s}$; it depends continuously on $r$ and $s$, other properties of  $F^X_{r,s}(ds')$    will be stated later.

  \medskip

  \begin{prop}
  \label{CDGP_7.3.0.1}
For any  $s\ge 0$ and $r'\in [0,X_s)$ the
function $(r,t) \to G^{X,\,{\rm neum}}_{s,t}(r',r)$, $\{(r,t):r \in [0,X_t),t>s\}$,
   \begin{equation}
\label{5.3.1bb}
G^{X,\,{\rm neum}}_{s,t}(r',r) = G^{{\rm neum}}_{s,t}(r',r)
- \int_s^t  F^X_{r',s}(ds') G^{{\rm neum}}_{s',t}(X_{s'},r)
\end{equation}
is   smooth and  for all $f\in L^{\infty}([0,X_t))$,
  \begin{equation}
\label{5.3.1bb.10}
E_{{r'},s}[ f(B_t);\tau^X_s>t]=  \int_{\mathbb R_+}  f(r) G^{X,\,{\rm neum}}_{s,t}(r',r)dr.
  \end{equation}
$G^{X,\,{\rm neum}}_{s,t}(r',r)$ solves the heat equation in  $\{(r,t):r \in [0,X_t),t>s\}$
with boundary conditions
  \begin{eqnarray}
\label{5.3.1bb.11}
&& \frac{\partial }{\partial r} G^{X,\,{\rm neum}}_{s,t}(r',r)\Big |_{r=0}=0 \qquad \text{for} \quad r'>0  \\
&&\hskip2cm \text{and} \qquad
G^{X,\,{\rm neum}}_{s,t}(r',X_t)=0 \qquad \text{for} \quad r'\ge 0\nn.
\end{eqnarray}
Finally if  $\rho(r',s)\in C([0,X_s),\mathbb R_+)$ then
   \begin{equation}
\label{5.3.1bbc}
\lim_{t\to s} \int \rho(r',s)G^{X,\,{\rm neum}}_{s,t}(r',r)dr' = \rho(r,s).
\end{equation}

 \end{prop}

 \medskip
 \noindent
 {\bf Proof.}  The smoothness of $G^{X,\,{\rm neum}}_{s,t}(r',r)$
is inherited from the smoothness of $G^{{\rm neum}}_{s,t}(r',r)$.
To prove \eqref{5.3.1bb.10} we first use the strong Markov property to write
  \begin{equation*}
E_{{r'},s}[ f(B_t)] =
E_{{r'},s}[ f(B_t);\tau^X_s>t] +
 \int_s^t  F^X_{r',s}(ds')  E_{{X_{s'}},s'}[ f(B_t)]
  \end{equation*}
and then \eqref{3.1.02.00}. Since  $G^{{\rm neum}}_{s,t}(r',r)$ solves the heat equation, then  $G^{X,\,{\rm neum}}_{s,t}(r',r)$
solves it as well. Similarly the Neumann boundary condition at 0
follows from the same property for  $G^{{\rm neum}}_{s,t}(r',r)$.  To prove
the Dirichlet condition at $X_t$ we will use the invariance of the law of
 Brownian motion under time reversal. Let $\delta>0$, $s^*:=s+\delta<t$, then, by the Markov property,
 \begin{equation}
\label{5.3.1bb.14}
 E_{{r'},s}[ f(B_t);\tau^X_s>t]=
 \int_{\mathbb R_+} h(r'') E_{{r''},s^*}[ f(B_t);\tau^X_{s^*}>t] dr''
  \end{equation}
with  $h(r'') = G^{X,\,{\rm neum}}_{s,s^*}(r',r'')$.
By the invariance of the law of
 Brownian motion under time reversal
  \begin{equation}
\label{5.3.1bb.13}
\int_{\mathbb R_+} h(r'')E_{{r''},s^*}[ f(B_t);\tau^X_{s^*}>t] dr''=   \int_{\mathbb R_+} f(r) E_{{r},s^*}[ h(B_t);\tau^{X'}_{s^*}>t] dr,
  \end{equation}
where $X'_{s^*+\si}= X_{t-\si}$, $\si\in [0,t-s^*]$; $h$ and $f$ are any two $L^{\infty}$ functions.
 By \eqref{5.3.1bb.13}
 \begin{equation}
\label{5.3.1bb.15}
 E_{{r'},s}[ f(B_t);\tau^X_s>t]=\int_{\mathbb R_+} f(r)  E_{{r},s^*}[ h(B_t);\tau^{X'}_{s^*}>t] dr
  \end{equation}
which  yields, for $f \ge 0$,
 \begin{equation}
\label{5.3.1bb.16}
E_{{r'},s}[ f(B_t);\tau^X_s>t]  \le \|h\|_\infty \int_{\mathbb R_+} f(r)  P_{{r},s^*}[ \tau^{X'}_{s^*}>t] dr
  \end{equation}
with $\|h\|_\infty \le (2\pi\delta)^{-1/2}$.
We fix $r \in [0,X_t)$ and take $f =f_\eps$, an approximate Dirac delta centered in $r$ with
support (for $\eps$ small enough) on $[r-\eps,r+\eps]$, so that (by the continuity of
$r\to G^{X,\,{\rm neum}}_{s,t}(r',r)$),
 \begin{equation}
\label{5.3.1bb.17}
G^{X,\,{\rm neum}}_{s,t}(r',r)= \lim_{\eps\to 0}E_{{r'},s}[ f_\eps(B_t);\tau^X_s>t]  \le (2\pi\delta)^{-1/2} P_{{r},s^*}[ \tau^{X'}_{s^*}>t]
  \end{equation}
which vanishes when $r \to X_t = X'_{s^*}$ recalling that $X_t$ is  Lipschitz (actually
piecewise $C^1$).

\vskip.2cm
\noindent
We will next prove \eqref{5.3.1bbc}. For the sake of brevity  we will write
$\rho(r)$ instead of $\rho(r,s)$.  It follows from \eqref{5.3.1bb.13} that
  \begin{equation}
\label{5.3.1bbc.1}
 \int_{\mathbb R_+} \rho(r')G^{X,\,{\rm neum}}_{s,t}(r',r)dr' =  E_{{r},s}[ \rho(B_t);\tau^{X'}_s>t].
\end{equation}
By Doob's inequality (see e.g. \cite{RY}),   and writing $ x_{t} = \min_{s'\in [s,t]} X'_{s'}$,
  \begin{equation}
\label{nn5.3.1.00}
P_{r,s} [ \tau^{X'}_{s}\le t] \le
P_{r,s}\Big[\max_{s'\in [s,t]} B_{s'}\ge  x_{t}\Big] \le
4 Q_{r,s}\Big[ B_{t}\ge x_{t}\Big]
\end{equation}
which vanishes in the limit $t\to s$.  Thus by \eqref{5.3.1bb},
   \begin{equation*}
\lim_{t\to s} \int_{\mathbb R_+} \rho(r',s)G^{X,\,{\rm neum}}_{s,t}(r',r)dr' =
\lim_{t\to s} \int_{\mathbb R_+} \rho(r',s)G^{{\rm neum}}_{s,t}(r',r)dr'
\end{equation*}
which is equal to $\rho(r,s)$ by \eqref{ch2.4.7} with $j=0$.
\qed

  \medskip

  \begin{coro}
Let $\rho(r',s)\in C([0,X_s),\mathbb R_+)$.  Then
		\begin{equation}
		\label{3.144444}
\rho(r,t) = \int_{\mathbb R_+}   \rho(r',s)  G^{X,\,{\rm neum}}_{{s},t}(r',r)\,dr'
+ j \int_s^t G^{X,\,{\rm neum}}_{s',t}(0,r) \,ds'  
		\end{equation}
 solves \eqref{intro.2}--\eqref{intro.3} with initial datum   $\rho(r',s)$ at time $s$ so that if $\phi$ is a bounded function
 		\begin{equation}
		\label{3.144444.0}
\int_{\mathbb R_+} \rho(r,t) \phi(r) \, dr = \int_{\mathbb R_+}  dr' \rho(r',s) E_{r',s}[\phi(B_t); \tau^X_s > t]
+ j\int_s^t  ds'  E_{0,s'}[\phi(B_t); \tau^X_{s'} > t].
		\end{equation}
 \end{coro}

 \medskip
 \noindent
 {\bf Proof.}  By Proposition \ref{CDGP_7.3.0.1} it follows that
$\rho(r,t)$ solves the heat equation \eqref{intro.2} with initial datum   $\rho(r',s)$ at time $s$
and that it vanishes at $X_t$. Using again  Proposition \ref{CDGP_7.3.0.1},
		\begin{equation*}
\frac {\partial}{\partial r}\rho(r,t)\Big|_{r=0}  
=
j \int_s^t  \frac {\partial}{\partial r} G^{{\rm neum}}_{s',t}(0,r)\Big|_{r=0}  ds'= -2j
		\end{equation*}
by \eqref{ch2.4.6} and \eqref{ch2.4.7.3}. \eqref{3.144444.0} follows from \eqref{3.144444} and
\eqref{5.3.1bb.10}.
 \qed

\vskip1cm

\section{Mass lost at the edge}
\label{ch:5.4}

Let $\rho(r,t), t\ge s$, be the solution of \eqref{intro.2}--\eqref{intro.3} which
at time $s$ is equal to $u\in L^\infty([0,X_s),\mathbb R_+)$.
We will give in Lemma \ref{7.4.0.3} a nice probabilistic representation for the mass
$\Delta^X_{I}(u)$, $I=[t_1,t_2]$, $t_1\ge s$, which has been lost in the time interval
$I$ by
$\rho(r,t)$, $\in I$.  The mass lost
$\Delta^X_{I}(u)$ is defined by
\begin{equation}
\label{5.3.1.1}
\Delta^X_{I} (u):=\int_{\mathbb R_+} \rho(r,t_1) dr-\int_{\mathbb R_+} \rho(r,t_2) dr  + j(t_2-t_1).
\end{equation}
Notice that if  $\Delta^X_{I}(u) = j|I|$ then the mass   is
conserved.
%
%

\medskip

     \begin{lem}[Mass loss]
     \label{7.4.0.3}
 With the above notation
 \begin{equation}
\label{5.3.1}
 \Delta^X_{I} (u) =   \int_{\mathbb R_+}   u(r') P_{r',s}\Big[\tau_s^X \in I \Big] \,dr'
 +j \int_s^{t_2}  P_{0,s'}\Big[ \tau^X_{s'}\in I\Big] \,ds'.
\end{equation}

\end{lem}

\medskip

\noindent
{\bf Proof.}
By integrating \eqref{3.144444} over $r$ we get 
  \begin{eqnarray*}
 \int_{\mathbb R_+} \rho(r,t) dr &=& \int_{\mathbb R_+} dr' {u(r')} P_{r',s}[\tau_s^X > t]
+j \int_{s}^{{t}}   P_{0,s'}[\tau_{s'}^X > t]\;ds' \\  &=&
 \int_{\mathbb R_+} dr' {u(r')} + j(t-s) -\int_{\mathbb R_+} dr' {u(r')} P_{r',s}[\tau_s^X \le t]
\\  &-& j \int_{s}^{{t}}   P_{0,s'}[\tau_{s'}^X \le t]\;ds'.
   \end{eqnarray*}
We use the above formula to compute $ \int \rho(r,t_2) dr - \int \rho(r,t_1) dr $.
We then use  the equality
   \begin{eqnarray*}
&& \int_s^{t_2}ds'\,P_{r',s'}[\tau_{s'}^X \le t_2]-\int_{s}^{t_1} ds' \,P_{r',s'}[\tau_{s'}^X \le t_1]
 = \int_s^{t_1} ds'\,  P_{r',s'}[\tau_{s'}^X \in I] \\&&\hskip2cm+ \int_{t_1}^{t_2} ds' \,P_{r',s'}[\tau_{s'}^X \le t_2].
    \end{eqnarray*}
We
 then
 get \eqref{5.3.1} after observing that $P_{r',s'}[\tau_{s'}^X \le t_2]=P_{r',s'}[\tau_{s'}^X \in I]$ for $s'\ge t_1$.
   \qed

\vskip.5cm
\noindent
Writing \eqref{5.3.1} in differential form we get

  \begin{equation}
  \label{5.3.1.2.2.001.a}
 \Delta^X_{I} (u) = \int_I\mu(dt)
   \end{equation}
where
  \begin{equation}
  \label{5.3.1.2.2.001.b}
\mu(dt)= \int_{\mathbb R_+} dr'\,  u(r')  F^X_{r',s}(dt) + j\int_s^t ds' F^X_{0,s'}(dt).
   \end{equation}
We also have
 \begin{equation}
  \label{5.3.1.2.2.001.c}
\Delta^X_{I} (u) = -\frac 12 \int_I \frac{\partial} {\partial r}\rho(r,t)\bigg|_{r=X_t} dt.
   \end{equation}
In fact by Theorem 2.6 in \cite{fasano}, $\frac{\partial} {\partial r}\rho(r,t)$ has a limit
when $r\to X_t$ under the assumption that the initial datum
is smooth  and that $X_t$ is Lipschitz.
We denote this limit
by $-2\la^X_{u,s}(t)$. Therefore from \eqref{5.3.1.2.2.001.a} and \eqref{5.3.1.2.2.001.c} we get
 \begin{equation}
  \label{5.3.1.2.2.001}
\int_{\mathbb R_+} dr' \, u(r')  F^X_{r',s}(dt) + j\int_s^t ds' F^X_{0,s'}(dt)=
  \la^X_{u,s}(t)
dt
   \end{equation}
and by \eqref{5.3.1bb} and \eqref{3.144444} the solution $\rho(r,t)$ can be written as
  \begin{equation}
\label{5.3.1.2.2.002}
\rho(r,t)   =  (T_{t-s}\rho(\cdot,s))(r) - \int_s^t \la^X_{u,s}(s') G^{\rm neum}_{s',t}(X_{s'},r)\, ds'
  \end{equation}
where   $T_{t-s}\rho(\cdot,s)$ is defined in \eqref{ch2.4.7}, namely
        \begin{equation*}
(T_{t-s}\rho(\cdot,s))(r)= \int_{\mathbb R_+} \rho(r',s)G_t^{\rm neum}(r',r)dr' + j
\int_s^t  G_{s',t}^{\rm neum}(0,r)\,ds'
   \end{equation*}
With $j=0$ this shows that the exit distribution of the Brownian has a density
with respect to Lebesgue
when the
starting point has a smooth distribution and $X_t$ is Lipschitz.  In \cite{PeSh}
it is proved that $F_{r,s}(dt)$ has a continuous density $g_{r,s}(t)$ if $X_t$ is $C^1$, the proof extends to our case when  $X_t$ is piecewise $C^1$ at all points with the possible exception of those where the derivative of $X_t$ is discontinuous.  Thus \eqref{5.3.1.2.2.001} becomes
  \begin{equation}
  \label{5.3.1.2.2.001.1}
\int_{\mathbb R_+} dr'\, u(r')  g_{r',s}(t) + j\int_s^t ds'  g_{0,s'}(t)=
  \la^X_{u,s}(t).
   \end{equation}

\renewcommand{\theequation}{\thesection.\arabic{equation}}
\setcounter{equation}{0}

\chapter{Existence of optimal sequences}
\label{ch:4b}

In this chapter we will prove
that there exist optimal sequences (see Definition \ref{defin2.1.3})
and in the following one we will conclude the proof of Theorem \ref{CDGPthm2.3.1}. The proofs in both chapters use extensively the representation of the solution of the heat equation
in terms of Brownian motions given in Chapter \ref{ch:4}.


\medskip

\renewcommand{\theequation}{\thesection.\arabic{equation}}
\setcounter{equation}{0}

\vskip2cm

\renewcommand{\theequation}{\thesection.\arabic{equation}}
\setcounter{equation}{0}

\section{The existence theorem}
\label{ch:5.2-section}

Recalling  Definition \ref{defin2.1.2} we  will prove in this chapter:

\medskip

\begin{thm}
\label{CDGPthm5.2.1}
For any $T>0$,   $\eps>0$ and $u \in \mathcal U$ there is an $\eps$-relaxed solution of the basic problem in $[0,T]$ with initial datum $u$, see  Definition \ref{defin2.1.2}.
\end{thm}

\medskip
\noindent
By the arbitrariness of $\eps$ Theorem
\ref{CDGPthm5.2.1}   proves the existence of optimal sequences.
 Since $\eps>0$ is fixed we will drop
it from the notation and simply write $X_t,\rho_0 $ for $X^{(\eps)}_t,\rho^{(\eps)}_0 $.
 We take $\rho_0 $ continuous, with compact support and such that $\int |u-\rho_0| \le \eps$.
Let then  $X_0$ be such  that $[0,X_0]$ contains the support of $\rho_0 $ and let
  \begin{equation}
  \label{CDGP5.2.1}
  X_t = X_0+\int_0^t ds V_s,\quad \text{$V_s$ piecewise constant in $[0,T]$}.
  \end{equation}
We will show that for a suitable choice of the  piecewise constant
velocity $V_s$ the solution of \eqref{intro.2}--\eqref{intro.3} is the  $\eps$-relaxed
solution we are looking for.  The proof is iterative, we introduce a time grid of length $t^*$, $t^* = j^{-1}\eps$, and prove that there is $V $ so that the solution
$\rho(r,t)$, $t\in [0,t^*]$, of \eqref{intro.2}--\eqref{intro.3} with $X_t=X_0+Vt$ is
such that
  \begin{equation}
  \label{CDGP5.2.2}
\bigg| \int_{\mathbb R_+} \rho(r,t) \, dr- \int_{\mathbb R_+} \rho_0(r)\, dr\bigg| \le \eps,\qquad \int_{\mathbb R_+} \rho(r,t^*)\, dr = \int_{\mathbb R_+} \rho_0(r) \, dr.
  \end{equation}
We will prove also uniformity on the initial datum to iterate.

\vskip.5cm

\section{The first step of the iteration}
\label{ch:4b.1}

Let
\[
V^* = - \frac {X_0}{t^*},\quad V > V^*,\quad X^{ V}_t=X_0+Vt
\]
and $u$ be a continuous, non-negative function with support in $[0,X_0]$ such that
$\int u = \int \rho_0$.
Let
 $u^{(V)}(r,t)$, $t\in[0,t^*]$, be defined as
\begin{equation}
\label{3.144444.00}
u^{(V)}(r,t) := \int_{\mathbb R_+}  G^{X,\,{\rm neum}}_{0,t}(r',r) u(r') \,dr'
+j \int_0^t G^{X,\,{\rm neum}}_{s,t}(0,r) \,ds.
\end{equation}
Then, see \eqref{3.144444}, $u^{(V)}(r,t)$
 is the solution of  \eqref{intro.2}--\eqref{intro.3}
with edge $X^{ V}_t$ and initial datum $u$.
We denote by $\Delta^{X^ V}_{[0,t]}(u)$ the mass lost in the time interval
$[0,t]$, see \eqref{5.3.1.1}.  The next lemma proves the intuitively evident fact that if $X_{t*}\to 0$ then all the
mass is taken out of the system, both that present initially and that injected through the origin.

\medskip

\begin{lem}
\label{CDGPlem5.2.1}

$\Delta^{X^ V}_{[0,t^*]}(u)$ converges to $jt^*+F(0;u)$ as $V\to V^*$.

\end{lem}

\medskip

\noindent
{\bf Proof.}
Let $V>V^*$ and shorthand $\delta= X^V_{t^*}=X_0+ Vt^* = (V-V^*)t^*$ so that $\delta\to 0$ as $V\to V^*$.
Then by \eqref{5.3.1.1} and \eqref{3.144444}
 \begin{eqnarray*}
  0 &\le & jt^*+F(0;u) - \Delta^{X^ V}_{[0,t^*]} (u) = \int_{\mathbb R_+} \rho(r,t^*)\, dr
\\&\le&    \int_{\mathbb R_+}   u(r') P_{r',0}\Big[B_{t^*} \le \delta\Big] \,dr'
 + \int_0^{t^*} jP_{0,s}\Big[ B_{t^*} \le \delta\Big] \,ds.
\end{eqnarray*}
By \eqref{3.1.02}, 
 \begin{eqnarray*}
&&
P_{r',0}\Big[B_{t^*} \le \delta\Big]  \le \frac {2\delta}{\sqrt{2\pi t^*}},\;\;\qquad P_{0,s}\Big[ B_{t^*} \le \delta\Big]
\le  \frac{ 2\delta}{\sqrt{2\pi (t^*-s)}}
\end{eqnarray*}
which yields
\begin{eqnarray*}
0&\le & jt^*+F(0;u)- \Delta^{X^{ V}}_{[0,t^*]} (u) 
\le    F(0;u) \cdot \frac {  2\delta}{\sqrt{2\pi t^*}} + \frac{  4j \delta \sqrt{t^*}}{\sqrt {2\pi}}.
\end{eqnarray*}
Thus  $\Delta^{X^V}_{[0,t^*]} (u) \to  jt^*+F(0;u)$ as $V\to V^*$ because $\delta = (V-V^*)t^*$.
\qed

\medskip
\noindent
Also the next lemma is quite evident as it claims that
there is no mass loss in the limit $V\to \infty$. These two lemmas  together with
 Lemma \ref{CDGPlem5.2.3}, which states that $\Delta^{X^ V}_{[0,t^*]}(u)$ depends continuously on $V$, will then show that there is a value of $V$ for which
 $\Delta^{X^ V}_{[0,t^*]}(u)=jt^*$.
 The second equality in \eqref{CDGP5.2.2}
 will then be proved.

\medskip

\begin{lem}
\label{CDGPlem5.2.2}

$\Delta^{X^ V}_{[0,t^*]}(u)$ converges to  0 as $V\to \infty$.

\end{lem}

\medskip

\noindent
{\bf Proof.}
Let $\zeta>0$,
$V = \zeta^{-\frac 34}$ and $r'< X_0 -\zeta^{\frac 14}$.  Call $t_k= k\zeta$ and $r_k = X_0+
Vt_k$, then
 \begin{eqnarray*}
&& P_{r',0}\Big[\tau_0^{X^{ V}} \le  t^*\Big] \le \sum_{k=1}^\infty P_{r',0}\Big[\max_{t\le t_k} B_t\ge r_{k-1}\Big].
\end{eqnarray*}
 Denoting by $Q_{r';0}$ the law of the Brownian motion on the whole $\mathbb R$ (i.e.\ without reflections at $0$), we have
 \[
 P_{r',0}\Big[\max_{t\le t_k} B_t\ge r_{k-1}\Big] \le 2  Q_{r',0}\Big[\max_{t\le t_k} B_t\ge r_{k-1}\Big].
 \]
 By Doob's inequality (see \cite{RY})  
  \begin{eqnarray}
\label{nn5.3.1}
P_{r',0}\Big[\max_{t\le t_k} B_t\ge r_{k-1}\Big] &\le& 4 Q_{r',0}\Big[ B_{t_k}\ge r_{k-1}\Big]
 \le 4  \int_{k\zeta^{\frac 14}}^\infty \frac {e^{-\frac {x^2}{2k \zeta}}}{\sqrt{2\pi k \zeta}}
\nn\\&\le&
{4 e^{-\frac{k}{4 \sqrt \zeta}} \int_{k\zeta^{\frac 14}}^\infty \frac {e^{-\frac {x^2}{4k \zeta}}}{\sqrt{2\pi k \zeta}} \le 4\sqrt 2 e^{-\frac{k}{4 \sqrt \zeta}}}
\end{eqnarray}
so that the first term on the right-hand side of \eqref{5.3.1} is bounded by
  \begin{eqnarray*}
&&  {\|u\|_{\infty}}  \zeta^{\frac 14} +4\sqrt 2  \, {\|u\|_{1}} \, \sum_{k=1}^\infty  e^{- \frac{{k}}{4\sqrt \zeta}}
\end{eqnarray*}
which vanishes as $\zeta\to 0$. An analogous argument (which is omitted) applies to  the second term on the right hand side of \eqref{5.3.1}.
\qed

\medskip

\begin{lem}
\label{CDGPlem5.2.3}

$\Delta^{X^ V}_{[0,t^*]}(u)$ depends continuously on $V$ in $(V^*,\infty)$.

\end{lem}

\medskip

\noindent
{\bf Proof.}
We consider the difference $\Delta^{X^ V}_{[0,t^*]}(u)-\Delta^{X^ {V'}}_{[0,t^*]}(u)$
with  $V^*<V<V'$ and  call $\delta=(V'-V)t^*$.  We need to prove that
the difference vanishes as $\delta\to 0$.
To make notation lighter
we shorthand  $X= \{X_t=X_0+Vt\}$ and $X'=\{X'_t=X_0+V't\}$. 
Then by  \eqref{5.3.1.2.2.001.a} and \eqref{5.3.1.2.2.001.b},
   \begin{eqnarray}
    \label{5.100}
&& \Big|\Delta^X_{[0,t^*]}(u) - \Delta^{X'}_{[0,t^*]}(u)\Big| \le  \int_0^{t^*-\delta} F(ds) \, P_{X_s,s} \Big[ \tau_{s}^{X'} >t^*\Big]  + R_\delta \\&&
F(ds) = \int_{\mathbb R_+} dr' u(r')F^X_{r',0}(ds) +j\int_0^s ds' \, F^X_{0,s'}(ds).\nn\\&&
R_\delta:= \int_{\mathbb R_+} dr' u(r')P_{r',0}\Big[\tau_0^{X} \in [t^*-\delta,t^*]\Big] \nn\\&& \hskip2cm
 +  j\int_0^{t^*}ds\,  P_{0,s}\Big[\tau_s^{X} \in \big[\max\{ s, t^*-\delta\},t^*\big]\Big].
  \nn
	\end{eqnarray}
We are going to prove that there is a function $o(\delta)$
which vanishes as $\delta\to 0$ so that
   \begin{equation}
    \label{5.101}
\sup_{0 \le s \le t^*-\delta} P_{X_s,s} \Big[ \tau_{s}^{X'} >t^*\Big] \le o(\delta).
	\end{equation}
Fix $s \le t^*-\delta$ and
define $\si_s:= \inf \{t \ge s \, : \: B_t \notin (X_s-\delta^{\frac 34}, X_s+\alpha \delta)\}$, with $\alpha > V'+1$,  then
   \begin{eqnarray}
    \label{5.102}
P_{X_s,s} \Big[ \tau_{s}^{X'} >t^*\Big] &\le&  P_{X_s,s} \Big[ \si_s > s +\delta\Big]+
P_{X_s,s} \Big[ B_{\si_s} < X'_{\si_s}; \si_s \le s +\delta\Big] \nn \\
&\le & P_{X_s,s} \Big[ \si_s > s +\delta\Big]+
P_{X_s,s} \Big[ B_{\si_s} = X_s- \delta^{\frac 34}\Big]
	\end{eqnarray}
because, by the choice of $\alpha$, if $B_{\si_s} = X_s +\alpha \delta$ then
$B_{\si_s} > X'_{\si_s}$, as  one can check that
$ X_s +\alpha \eps > X'_{s+\delta}$.
Since $P_{r;s} \Big[ B_{\si_s} = X_s- \delta^{\frac 34}\Big]$ is a linear function of $r$
which has value 1 at $r = X_s- \delta^{\frac 34}$ and is equal
to $0$ at $r=X_s+\alpha\delta$, it  follows that
   \begin{equation}
    \label{5.103}
P_{X_s,s} \Big[ B_{\si_s} = X_s- \delta^{\frac 34}\Big] \le \alpha \, \delta^{\frac 14}
	\end{equation}
Since the probability density of $B_{s+\delta}-X_s$ is  $e^{-x^2/(2\delta)}(2\pi \delta)^{-1/2}$ we have
   \begin{equation}
    \label{5.104}
P_{X_s,s} \Big[ \si_s >s + \delta\Big] \le P_{X_s,s} \Big[|B_{s+\delta} - X_s| \le \delta^{\frac 34}\Big]
\le  \frac{2}{\sqrt{2\pi }} \cdot \delta^{\frac 14}
	\end{equation}
so that \eqref{5.101} is proved.  We then have that the first term on the right-hand side of \eqref{5.100}
is bounded by:
\[
o(\delta) \int_0^{t^*} h(s)\, ds \le o(\delta) \cdot \big(F(0;u)+jt^*\big).
\]

\noindent
We shall next  bound the probabilities in $R_\delta$.  Call $Y=X _{t^*-\delta}= X_0 + V(t^*-\delta)$, then
   \begin{eqnarray}
    \label{5.105}
&& P_{r',0}\Big[\tau_0^{X} \in [t^*-\delta,t^*]\Big] \le  P_{r',0}\Big[B_{t^*-\delta} \in [Y-\delta^{\frac 14},Y]\Big]
\nn\\&& \hskip2cm + \sup_{r'' \le Y-\delta^{\frac 14}}
 P_{r'',t^*-\delta}\Big[ \max _{t\in [t^*-\delta,t^*] } B_t\ge Y\Big].
	\end{eqnarray}
As before we have
  \begin{equation*}
 P_{r',0}\Big[B_{t^*-\delta} \in [Y-\delta^{\frac 14},Y]\Big] \le \frac{\delta^{\frac 14}}{\sqrt{2\pi (t^*-\delta)}}.
 \end{equation*}
Now suppose $r'' \in [0, Y-\delta^{\frac 1 4}]$, then
 \begin{eqnarray*}
 && P_{r'',t^*-\delta}\Big[ \max _{t\in [t^*-\delta,t^*] }B_t\ge  Y\Big]
 \le P_{r'',t^*-\delta}\Big[ \max _{t\in [t^*-\delta,t^*] }(B_t - r'')
   \ge  \delta^{\frac 14} \Big].
  	\end{eqnarray*}
 	 By the same argument used in \eqref{nn5.3.1}, the latter is bounded by
     \begin{equation*}
  2
  P_{r',t^*-\delta}\Big[ B_{t^*} - r'\ge  \delta^{\frac 14}\Big] \le 2\int_{ \delta^{\frac 14}}^\infty
  \frac{e^{-\frac {x^2}{2\delta}}} {\sqrt{2\pi \delta}} \; dx
  \le 4\sqrt 2e^{-\frac{1}{4 \sqrt \delta}}
 	\end{equation*}
Analogous bounds are proved for $P_{0,s}\Big[ \tau^{X}_s\in [t^*-\delta,t^*]\Big]$, we omit the details.
We have thus proved that also $R_\delta$ is infinitesimal with $\delta$.

\qed

\vskip.5cm

\section{The iteration}
\label{ch:4b.2}

\medskip

\begin{coro}

There exists a  $V$ such that
\begin{equation}
\label{5.8888}
\Delta^{X^V}_{[0,t^*]}(u) = jt^*\qquad \text{and} \qquad \sup_{t\le t^*} |\Delta^{X^ V}_{[0,t]}(u) -jt| \le jt^*.
\end{equation}

\end{coro}

\medskip

\noindent
{\bf Proof.}
The equality in \eqref{5.8888} follows from
Lemmas \ref{CDGPlem5.2.1}--\ref{CDGPlem5.2.3}.  The last statement
holds because $\Delta^{X^ V}_{[0,t]}(u)$ is a non-decreasing function
of $t$ which is equal to $jt^*$ at $t=t^*$.  \qed

\medskip

By \eqref{5.3.1bb} and \eqref{3.144444.00}
the function $u^V(r,t^*)$ is continuous with support on $[0,X^V_{t^*}]$
and  by \eqref{5.8888}, $\int u^V(r,t^*)dr = \int \rho_0(r)dr$.
Thus $u^V(r,t^*)$ has the same properties as the initial $u$ and we can iterate the procedure constructing a function $\rho^X(r,t)$ with $X_t$ having constant velocity in each interval $[kt^*,(k+1)t^*)$ and such that
$|\int \rho^X(r,t)-\int \rho_0(r,t)| \le \eps$
 at all times $t\in [0,T]$.  Theorem \ref{CDGPthm5.2.1} is then proved.

\renewcommand{\theequation}{\thesection.\arabic{equation}}
\setcounter{equation}{0}

\chapter{Proof of the main theorem}
\label{ch:7}

In this chapter we will first prove Theorem  \ref{CDGPthm2.4.1.00} and then Theorem \ref{CDGPthm2.3.1}.  The main point will be to show that the elements of an optimal sequence are eventually squeezed between the upper and lower barriers which will be proved using  the representation of the solution of \eqref{intro.2} and \eqref{intro.4} in terms of Brownian motions, as discussed in Chapter \ref{ch:5}.  In Section \ref{ch:7.5} we will use this to prove Theorem  \ref{CDGPthm2.4.1.00} while  Theorem \ref{CDGPthm2.3.1} will  be proved in Section \ref{ch:7.6}.

\vskip1cm

\section{The key inequality}
\label{ch:7.1}
We fix $T>0$ and  $\rho_0 \in \mathcal U$.   Let $\delta_0>0$
be such that
$\rho_0 \in \mathcal U_{\delta_0}$, by default in the sequel $\delta < \delta_0$.
We also fix an optimal sequence in $[0,T]$ with  initial datum $\rho_0$, see
Definition \ref{defin2.1.3}.  We will prove:

\medskip

\begin{thm}
\label{CDGPthm9.1n}
Let $t\in (0,T]$, $\delta \in \{2^{-k}t, k \in \mathbb N\}$ with $k$ large enough.  Then
 \begin{equation}
\label{CDGP7.1.1}
S^{\delta,-}_{t} \rho^{(\eps_n)}(\cdot,0)  \preccurlyeq   \rho^{(\eps_n)}(\cdot, t) \preccurlyeq
S^{\delta,+}_{t}  \rho^{(\eps_n)}(\cdot,0) \;\;\;\qquad \text{\rm modulo}\quad  \frac t \delta \;2\eps_n.
  \end{equation}

\end{thm}

\medskip
The proof of Theorem \ref{CDGPthm9.1n} is reported in Section  \ref{ch:7.2}. We first use it to prove Theorem  \ref{CDGPthm2.4.1.00}.  

\vskip1cm

\section{Proof of Theorem  \ref{CDGPthm2.4.1.00}}
\label{ch:7.5}
From the key inequality \eqref{CDGP7.1.1} Theorem  \ref{CDGPthm2.4.1.00} easily follows. In fact
by  definition of optimal sequences, $\int |\rho^{(\eps_n)}(r,0) -   \rho_0(r)| dr\le \eps_n$, then
\[
\rho_0 \preccurlyeq \rho^{(\eps_n)}(\cdot,0) \quad \text{modulo}\; \eps_n,\qquad
\rho^{(\eps_n)}(\cdot,0)\preccurlyeq  \rho_0 \quad \text{modulo}\; \eps_n.
\]
Thus by \eqref{IFB.0} $S^{\delta,-}_{t} \rho_0  \preccurlyeq
S^{\delta,-}_{t} \rho^{(\eps_n)}(\cdot,0)$ modulo  $\eps_n$. By Lemma \ref{prop4.2} and
\eqref{CDGP7.1.1}, $S^{\delta,-}_{t} \rho_0  \preccurlyeq  \rho^{(\eps_n)}(\cdot, t)$ modulo
$\eps_n+ \frac t \delta \;2\eps_n$:  An analogous argument applies to
$S^{\delta,+}_{t} \rho_0$, hence
 \begin{equation}
\label{CDGP7.1.2}
S^{\delta,-}_{t} \rho_0  \preccurlyeq   \rho^{(\eps_n)}(\cdot, t) \preccurlyeq
S^{\delta,+}_{t}  \rho_0 \;\;\;\qquad \text{modulo}\quad  \frac t \delta\; 2 \eps_n + \eps_n.
  \end{equation}
We keep $\delta$ fixed in \eqref{CDGP7.1.2}
and let $\eps_n\to 0$:
 \begin{equation}
\label{CDGP7.1.3n}
F(r;S^{\delta,-}_{t} \rho_0) \le \liminf_{\eps_n\to 0} F(r; \rho^{(\eps_n)}(\cdot, t) )
\le \limsup_{\eps_n\to 0} F(r; \rho^{(\eps_n)}(\cdot, t)) \le
F(r;S^{\delta,+}_{t}  \rho_0).
  \end{equation}
By Theorem \ref{CDGPthm2.4.1} letting $\delta\to 0$,
 \begin{equation*}
F(r;S_{t} \rho_0) \le \liminf_{\eps_n\to 0} F(r; \rho^{(\eps_n)}(\cdot, t) )
\le \limsup_{\eps_n\to 0} F(r; \rho^{(\eps_n)}(\cdot, t)) \le
F(r;S_{t}  \rho_0)
  \end{equation*}
  which proves that 
   \begin{equation}
\label{CDGP7.1.4n}
 \lim_{\eps_n\to 0} F(r; \rho^{(\eps_n)}(\cdot, t) ) =F(r;S_{t} \rho_0).
  \end{equation}
This shows that $\rho^{(\eps_n)}(\cdot, t)$ converges in distribution
to $S_{t} \rho_0$ and hence it converges weakly as well.

%
%

\renewcommand{\theequation}{\thesection.\arabic{equation}}
\setcounter{equation}{0}
 \vskip1cm

\section{Proof of Theorem  \ref{CDGPthm9.1n}}
\label{ch:7.2}

To simplify notation we write
$\eps$ for $\eps_n$, $u$ for $\rho^{(\eps_n)}_0$, call $t = N\delta$,
 $u(r,k\delta)= \rho^{(\eps)}(r, k\delta)$.
Theorem  \ref{CDGPthm9.1n} then follows from showing that
for all $k \le N$:
  \begin{equation}
\label{CDGP7.1.3}
S^{\delta, - }_{k\delta}u    \preccurlyeq   u(\cdot, k\delta) \preccurlyeq
S^{\delta, + }_{k\delta}u    \;\;\;\qquad \text{modulo}\quad 2 k\eps
  \end{equation}
because \eqref{CDGP7.1.1} is \eqref{CDGP7.1.3} with $k =N$.  The proof is
by induction on $k$.  The case $k=1$
is notationally simpler and even if it
can be recovered by the induction procedure
when we start it from $k=0$ (for which \eqref{CDGP7.1.3} trivially holds), we will prove it explicitly to give an idea of the general case. The only difference
when treating the case $k=1$ is that \eqref{CDGP7.1.3} holds modulo $\eps$, while in the general case there is the extra factor 2: this is due to the fact that the approximate mass conservation gives:
	\begin{equation}
	\label{CDGP7.1.4}
|\Delta_{[0,t]}(u)-jt| \le \eps,\quad |\Delta_{[s,t]}(u)-j(t-s)| \le 2\eps.
\end{equation}

 \vskip1cm

\subsection{The first step of the induction}
\label{ch:7.3}

We will prove separately the two inequalities in
\eqref{CDGP7.1.3} with $k=1$.

\subsubsection{Lower bound}
\label{subsec:7.3.1}

We shorthand $S^{\delta, - }_{\delta}u =v^-(\cdot,\delta)$.  With this notation the lower bound in \eqref{CDGP7.1.3} for $k=1$ reads as
	\begin{equation}
	\label{CDGP7.3.1}
v^-(\cdot,\delta)= C_\delta T_\delta u \preccurlyeq   u(\cdot,\delta)  \quad \text{modulo}\;\;    \eps.
\end{equation}
By \eqref{3.1a}
 \begin{equation*}
 F(r;T_\delta u )=
  \int_{\mathbb R_+}   \;  u(r') \; P_{r',0}\Big[ B_\delta \ge r\Big]dr' +
  j \int_0^\delta  P_{0,s}\Big[
  B_\delta \ge r\Big]ds
 \end{equation*}
while, by \eqref{3.144444.0},
 \begin{eqnarray*}
F(r;  u(\cdot,\delta)) &=&
\int_{\mathbb R_+} \; u(r') \, P_{r'{,0}}\Big[\tau_0^X>\delta ;\; B_\delta \ge r\Big] dr'\nn\\ & +& j \int_0^\delta \, P_{0,s}\Big[\tau_s^X>\delta ; \;B_\delta \ge r\Big]ds.
\end{eqnarray*}
Since
\begin{equation*}
 \sup_{t\in [0,\delta]} \Big| \Delta^X_{[0,t]}(u)- jt\Big| \le \eps.
\end{equation*}
by \eqref{5.3.1} 
 \begin{eqnarray*}
F(r;  u(\cdot,\delta)) &\ge&
\int_{\mathbb R_+} \; u(r') \, P_{r'{,0}}\Big[ B_\delta \ge r\Big] dr' +
j \int_0^\delta \, P_{0,s} \Big[  B_\delta \ge r\Big]ds\\&-& (j\delta + \eps)
\ge F(r;T_\delta u ) -(j\delta + \eps).
\end{eqnarray*}
Thus
\begin{equation*}
T_\delta u \preccurlyeq  u(\cdot,\delta)  \quad\text{modulo}\;\; j\delta + \eps
 \end{equation*}
and therefore by \eqref{4.15a}
\begin{equation*}
C_\delta T_\delta u \preccurlyeq  u(\cdot,\delta)  \quad\text{modulo}\;\;  \eps
 \end{equation*}
which proves \eqref{CDGP7.3.1}.

\vskip.5cm

\subsubsection{Upper bound}
\label{subsec:7.3.2}
We shorthand $S^{\delta, +}_{\delta}u =v^+(\cdot,\delta)$;
 $u_1 = u-C_\delta u$; $u_0=u-u_1=C_\delta u$.  Then 
\begin{eqnarray*}
v^+(r,\delta) &=& \int_{\mathbb R_+}  u_0(r')G_{0,\delta}^ {\rm neum}(r',r)dr' +
j\int_0^\delta G_{s,\delta}^ {\rm neum}(0,r) ds= T_\delta u_0(r).
 \end{eqnarray*}
By   \eqref{5.3.1.2.2.002} and writing in the sequel $\la^X_{u_0}(s)$ for $\la^X_{u_0,0}(s)$,
\begin{equation*}
u(r,\delta)= T_\delta u_0(r) -\int_0^\delta \la^X_{u_0}(s) G_{s,\delta}^ {\rm neum}(X_{s},r)ds
+ \int_{\mathbb R_+} u_1(r') G^{X, \rm neum}_{0,\delta}(r',r)dr'.
 \end{equation*}
Calling
\begin{equation}
\label{CDGP7.3.3}
I(r):= F(r; v^+(\cdot, \delta))-F(r; u(\cdot, \delta))
\end{equation}
we then get
\begin{equation*}
I(r)=
\int_0^\delta ds  \la^X_{u_0}(s) P_{X_{s}, s}[B_\delta \ge r]
 -
\int_{\mathbb R_+} dr'  u_1(r')P_{r', 0}[B_\delta \ge r; \tau^X_0 > \delta]
\end{equation*}
which can be rewritten as
\begin{eqnarray}
\label{CDGP7.3.4}
I(r) &=&
\int_0^\delta ds \, \la^X_{u_0}(s) \, P_{X_{s},s}[B_\delta \ge r] \nn\\
 &-&
\int_{\mathbb R_+} dr'  \, u_1(r') \, P_{r',0}[\tau^X_0 > \delta] \,P_{r',0}[B_\delta \ge r\,|\,\tau^X_0 > \delta]
\end{eqnarray}
with
$$
P_{r',0}[B_\delta \ge r\,|\,\tau^X_0 > \delta] = \frac{P_{r',0}[B_\delta \ge r; \tau^X_0 > \delta]}{P_{r',0}[\tau^X_0 > \delta]}
$$
the conditional probability that $\{B_\delta \ge r\}$ given that
$\{\tau^X_0 > \delta\}$.
Let
\begin{equation}
\label{CDGP7.3.5}
D:= \int_0^\delta ds  \, \la^X_{u_0}(s) - \int_{\mathbb R_+} dr'  \, u_1(r') \, P_{r',0}[\tau^X_0 > \delta].
\end{equation}
$D$ can be rewritten and then bounded as follows:
\begin{equation}
\label{CDGP7.3.6}
D= \int_0^\delta ds \,  \la^X_{u}(s) - \int_{\mathbb R_+} dr'  \, u_1(r'), \quad |D| \le \eps.
\end{equation}
The inequality follows from the following facts:  $\int_0^\delta \la^X_u(s) ds= \Delta_{0,\delta}(u)$,  $\int dr'   u_1(r')= j\delta$, by the definition of $u_1$ and
$| \Delta_{0,\delta}(u) -j\delta| \le \eps$.

\vskip.2cm
\noindent
Let $q$ be such that
\begin{equation}
\label{CDGP7.3.7}
  \int_0^\delta ds  \, \la^X_{u_0}(s)  = \int_{\mathbb R_+} dr'  \, q \, u_1(r')\, P_{r',0}[\tau^X_0 > \delta].
\end{equation}
By \eqref{CDGP7.3.5}--\eqref{CDGP7.3.6}
\begin{equation}
\label{CDGP7.3.8}
     \Big|(1-q) \int_{\mathbb R_+} dr'  \,  u_1(r') \, P_{r',0}[\tau^X_0 > \delta]\Big|
     =|D| \le \eps
\end{equation}
so that
\begin{eqnarray}
\label{CDGP7.3.9}
&& \hskip-.5cm\bigg|I(r) - \bigg(
\int_0^\delta ds  \la^X_{u_0}(s) P_{X_{s},s}[B_t \ge r] \\
 &&\hskip.2cm
-\int_{\mathbb R_+} dr'  \, q \, u_1(r') \, P_{r',0}[\tau^X_0 > \delta]\,  P_{r',0}[B_t \ge r\,|\,\tau^X_0 > \delta]
 \bigg)\bigg| \le \eps. \nn
\end{eqnarray}
 We will prove   that for any $r$ the curly bracket is non-negative which by \eqref{CDGP7.3.3}  gives the desired upper bound
	\begin{equation*}
 u(\cdot,\delta) \preccurlyeq
 v^+(\cdot,\delta)  \quad \text{modulo}\;\;    \eps.
\end{equation*}
Since the measures $ \la^X_{u_0}(s) ds $ on $[0,\delta]$ and
$\{qu_1(r')P_{r',0}[\tau^X_0 > \delta]\} dr'$ on $[0,X_0]$ have
same mass, and since they are both non-atomic, by the theory of Lebesgue measures,
 see for instance Roklin, \cite{roklin},
there is a map $\Ga:[0,X_0]\to [0,\delta]$ such that
\begin{eqnarray}
\label{CDGP7.3.10}
&&\hskip-1cm\int_0^\delta ds \,  \la^X_{u_0}(s) \, P_{X_{s},s}[B_t \ge r] \nn\\&&
= \int_{\mathbb R_+} dr'\, 
q \, u_1(r')\, P_{r',0}[\tau^X_0 > \delta] \,
P_{ X_{\Ga(r')},\Ga(r')}\Big[B_\delta \ge r\Big].
	\end{eqnarray}
In the next subsection we will prove that
\begin{equation}
 P_{r',0}\Big[B_\delta \ge r\,\big|\, {\tau_0^X} >\delta\Big] \le
P_{X_{t},t}\Big[B_\delta \ge r\Big],\quad r' \in [0, X_0),\; t\in [0,\delta)
\label{CDGP7.3.11}
	\end{equation}
which completes the proof of the upper bound.

\medskip

\subsection{A stochastic inequality}
\label{subsec:7.3.3}
In this subsection we will prove  \eqref {CDGP7.3.11}, by using coupling between Brownian motions.
Let $r'$ and $t$ be as in \eqref {CDGP7.3.11}.  Recalling \eqref{5.3.1bb.14},
\begin{equation*}
 P_{r',0}\Big[B_\delta \ge r\,\big|\, {\tau_0^X} >\delta\Big]
= \int_0^{X_t} G^{X,\,{\rm neum}}_{0,t}(r',z) \, P_{z,t}\Big[B_\delta \ge r\,\big|\, {\tau_t^X} >\delta\Big]dz
	\end{equation*}
so that \eqref {CDGP7.3.11} will follow from
\begin{equation}
\label{CDGP7.3.12}
 P_{z,t}\Big[B_\delta \ge r\,\big|\, {\tau_t^X} >\delta\Big] \le
P_{X_{t},t}\Big[B_\delta \ge r\Big],\quad z \in [0, X_t),\; t\in [0,\delta)
	\end{equation}
which will be proved in the remaining part of this subsection.

\medskip

Let $\ga^{-1}$ be a positive integer (eventually $\ga\to 0$), $B^{(1)}_i$, $i=1,..,\ga^{-1}$
independent Brownian motions which start moving at time $t$  from $X_t$ and
denote by $P^{(1)}$ their law and by $E^{(1)}$ the corresponding
expectation.  We will use the identity:
\begin{eqnarray}
\label{CDGP7.3.13}
&& P_{X_{t},t}\Big[B_\delta \ge r\Big]=
E^{(1)}\Big[ \ga \sum_{i=1}^{\ga^{-1}} {\mathbf 1_{[r, +\infty)}(B^{(1)}_i(\delta))}\Big].
	\end{eqnarray}
We can proceed in an analogous way with $ P_{z,t}\Big[B_\delta \ge r\,\big|\,{\tau_t^X}>\delta\Big]$
which is now conveniently rewritten as
\begin{eqnarray}
\label{CDGP7.3.14}
&& P_{z,t}\Big[B_\delta \ge r\,\big|\, {\tau_t^X} >\delta\Big]=  P_{z,t}\Big[B_\delta \ge r\,;\,{\tau_t^X}  >\delta\Big]\nn\\&&\hskip1cm \times
\Big( 1 + \Big\{ \frac{1}{1-\alpha(z)} - 1\Big\} \Big),\quad \alpha(z) = P_{z,t}\Big[
\tau_t^X \le \delta\Big].
	\end{eqnarray}

Calling $N_\ga := $ the integer part of $\ga^{-1}\big\{ \frac{1}{1-\alpha(z)} - 1\big\}$,
we then consider  $B^{(2)}_i$, $i=1,..,\ga^{-1}+N_\ga$,
independent Brownian motions which start  at time $t$  from $z$ and are removed once they reach
the edge $X_t$. We denote by $P^{(2)}$ such a law and by $E^{(2)}$ the corresponding expectation. We have:
\begin{equation}
\label{CDGP7.3.15-first}
P_{z,t}\Big[B_\delta \ge r\,\big|\, {\tau_t^X} >\delta\Big]= \lim_{\ga\to 0}E^{(2)}\Big[ \ga \sum_{i=1}^{\ga^{-1}+N_\ga} {\mathbf 1_{[r,+\infty)}(B^{(2)}_i(\delta))}\Big].
	\end{equation}
The equality follows using \eqref{CDGP7.3.14}: it holds only in the limit because of the integer part in the definition of $N_\ga$.  We are going to couple the
Brownians $B^{(1)}_i(s)$ and $B ^{(2)}_i(s)$: this means that we will define a probability $P$ on all $B^{(1)}_i$ and $B^{(2)}_i$ such that the marginal law of the $B^{(1)}_i$ is $P^{(1)}$ and the marginal law of the $B^{(2)}_i$ is $P^{(2)}$.

At the initial time $t$ we have $\ga^{-1}+ N_\ga$ (2)-particles at $z$ and
$\ga^{-1}$ (1)-particles at $X_t$.  We say that the (2)-particle  with label $i\le \ga^{-1}$
is married with the (1)-particle  with the same label $i$.  The (2)-particles with
label $i> \ga^{-1}$ are called single.
We are going to couple the evolution of the married pairs in the following way. $B^{(1)}_i(s)$ and $B^{(2)}_i(s)$, $s\ge t$, $i\le \ga^{-1}$ move independently
of each other till when they meet, from then on they move in the same way (observe that
$B^{(2)}_i(s) \le B^{(1)}_i(s)$ because the inequality holds initially).
The coupling stops when $B^{(2)}_i(s)=X_s$  because at that time $B^{(2)}_i(s)$ must be erased. We let all married pairs move independently of each other and  of the single particles and this defines the coupled process till the first time $s$ when the (2)-particle, say with label $i$, in a married pair
reaches $X_s$. We define the process after time $s$ by redefining the broken pair: we take the single(2)-particle still alive at   time $s$ with smallest label, say $j$, and we say that at time $s^+$ the (2)-particle with label $j$ is married with the (1)-particle with label $i$.  The process is then continued with same rules till time $\delta$.  If it happens that
there are no longer (2) single particles, a broken pair cannot be reconstructed and there are (1)-particles which become single.
 We denote by $P$ the law of this coupled process and by $E$ the corresponding expectation.
The important features  of this construction are:

\begin{itemize}
\item In a married pair the position of the (1)-particle is always $\ge$ than the position of  the (2)-particle.

  \item  Single particles are all of type (2) till when the number of deaths
  of (2)-particles is $\le N_\ga$
  and are all of type (1) afterwards.

\end{itemize}

\noindent
Therefore
\begin{eqnarray}
\label{CDGP7.3.15-second}
&& P_{ X_{t},t}\Big[B_\delta \ge r\Big]- P_{z,t}\Big[B_\delta \ge r\,\big|\,{\tau_t^X}>\delta\Big] \ge
{-}\lim_{\ga\to 0} E\Big[\ga K_\ga \Big]
	\end{eqnarray}
where
\[
K_\ga = \max\Big\{0;N_\ga - \sum_{i=1}^{\ga^{-1}+N_\ga} \mathbf 1_{{B^{(2)}_i (s)=X_s, \;\text{for some $s\in [t,\delta]$}}}\Big\}.
\]
By the law of large numbers for independent variables,
for any $\zeta>0$,
\begin{eqnarray}
\label{CDGP7.3.18}
&&\lim_{\ga \to 0} P \Big[\; \Big| \ga\sum_{i=1}^{\ga^{-1}+N_\ga} \mathbf 1_{{B^{(2)}_i (s)=X_s, \;\text{for some $s\in [t,\delta]$}}}\nn\\&& \hskip2cm
- (1+\ga N_\ga) P_{z,t}[ {\tau_t^X} \le \delta] \Big| \le \zeta\Big] = 1.
	\end{eqnarray}
Recalling the definition of $\alpha(z)$ we have
\[
\lim_{\ga \to 0} (1+\ga N_\ga) P_{z,t}[ {\tau_t^X} \le \delta] = \frac{\alpha(z)}{1
-\alpha(z)} = \lim_{\ga \to 0} \ga N_\ga.
\]
Thus from \eqref{CDGP7.3.18}  we have that
\begin{equation*}
\lim_{\ga\to 0} P\Big[ \;\Big| \ga\sum_{i=1}^{\ga^{-1}+N_\ga} \mathbf 1_{{B^{(2)}_i (s)=X_s, \;\text{for some $s\in [t,\delta]$}}} -\ga N_\ga\Big|\le \zeta\Big] =1,
	\end{equation*}
hence $\dis{\lim_{\ga\to 0} P\big[\ga K_\ga \le \zeta \big]=1}$ which yields
$\dis{\lim_{\ga\to 0} E\big[\ga K_\ga  \big]=0}$, thus  the right-hand side of
\eqref{CDGP7.3.15-second} is equal to 0.

\qed

\renewcommand{\theequation}{\thesection.\arabic{equation}}
\setcounter{equation}{0}

 \vskip1cm

\subsection{The generic step of the induction}
\label{ch:7.4}

We suppose by induction that for all $n\le k$:
 \begin{equation}
\label{CDGP7.4.1}
S^{\delta, - }_{n\delta}u    \preccurlyeq   u(\cdot, n\delta) \preccurlyeq
S^{\delta, + }_{n\delta}u    \;\;\;\qquad \text{modulo}\quad 2 n\eps.
  \end{equation}
{\em The lower bound.}  Call $u^*(\cdot)=u(\cdot, k\delta)$.  Then
 \begin{equation}
\label{CDGP7.4.2}
S^{\delta, - }_{\delta}u^*    \preccurlyeq   u(\cdot, (k+1)\delta)  \;\;\;\qquad \text{modulo}\quad 2 \eps.
  \end{equation}
The proof of \eqref{CDGP7.4.2} is the same as that in Subsection \ref{subsec:7.3.1}, here we have a bound with $2\eps$ because unlike in Subsection \ref{subsec:7.3.1} we have
\begin{equation*}
 \sup_{t\in [0,\delta]} \Big| \Delta^X_{[0,t]}(u^*)- jt\Big| \le 2\eps.
\end{equation*}
By the induction hypothesis
 \begin{equation}
\label{CDGP7.4.3}
S^{\delta, - }_{k\delta}u    \preccurlyeq   u^*  \;\;\;\qquad \text{modulo}\quad 2 k\eps.
  \end{equation}
Then by Theorem \ref{IFBthm.1}
 \begin{equation}
\label{CDGP7.4.4}
S^{\delta, - }_{(k+1)\delta}u    \preccurlyeq  S^{\delta, - }_{ \delta} u^*  \;\;\;\qquad \text{modulo}\quad 2 k\eps
  \end{equation}
which by \eqref{CDGP7.4.2} yields
 \begin{equation}
\label{CDGP7.4.5}
S^{\delta, - }_{(k+1)\delta}u    \preccurlyeq u(\cdot, (k+1)\delta)  \;\;\;\qquad \text{modulo}\quad 2 (k+1)\eps.
  \end{equation}

\medskip
\noindent
{\em The upper bound.}  The same proof applies for the upper bound.  We just repeat it for the reader's convenience.  We have
 \begin{equation}
\label{CDGP7.4.6}
 u(\cdot, (k+1)\delta) \preccurlyeq
 S^{\delta, + }_{\delta}u^*   \;\;\;\qquad \text{modulo}\quad 2 \eps.
  \end{equation}
Using the same proof as that in Subsection \ref{subsec:7.3.2}, again the bound with $2\eps$ is
 due to  the bound   $|\Delta^X_{[0,t]}(u^*)- jt|\le 2\eps$.
By the induction hypothesis
 \begin{equation}
\label{CDGP7.4.7}
  u^*   \preccurlyeq
  S^{\delta, +}_{k\delta}u \;\;\;\qquad \text{modulo}\quad 2 k\eps.
  \end{equation}
Then by Theorem \ref{IFBthm.0}
 \begin{equation}
\label{CDGP7.4.8}
  S^{\delta, + }_{\delta} u^*  \preccurlyeq
S^{\delta, + }_{(k+1)\delta}u  \;\;\;\qquad \text{modulo}\quad 2 k\eps
  \end{equation}
which by \eqref{CDGP7.4.6} yields
 \begin{equation}
\label{CDGP7.4.9}
u(\cdot, (k+1)\delta) \preccurlyeq
S^{\delta, +}_{(k+1)\delta}u \;\;\;\qquad \text{modulo}\quad 2 (k+1)\eps.
  \end{equation}

\vskip2cm

\renewcommand{\theequation}{\thesection.\arabic{equation}}
\setcounter{equation}{0}

\section{Proof of Theorem \ref{CDGPthm2.3.1}}
\label{ch:7.6}

\begin{itemize}

\item (a) is proved in Theorem \ref{CDGPthm5.2.1}.

\item (b) is proved in Theorem \ref{CDGPthm2.4.1.00}.

\item  (c) is also proved in  Theorem \ref{CDGPthm2.4.1.00} where we identify a relaxed solution to the element $S_t\rho_0$ which separates the barriers.

\item  (d) follows from the identification theorem,  Theorem \ref{CDGPthm2.4.1.00},
and item (5) of Theorem \ref{CDGPthm2.4.1}.

\item  (e) follows from property (9) of Theorem \ref{CDGPthm2.4.1} (via Theorem \ref{CDGPthm2.4.1.00}).
    
    \item (f) follows from  Theorem \ref{CDGPthm2.4.1.00} and item (7) in Theorem \ref{CDGPthm2.4.1}.

\item (g) follows from (7) of  of Theorem \ref{CDGPthm2.4.1} (via Theorem \ref{CDGPthm2.4.1.00}).

\item (h) Let $\rho_0$ be a classical initial datum and let $u$ be the (local in time) solution whose existence has been proved in Theorem \ref{CDGPthm2.2.1}. Since $u$ can be regarded as an optimal sequence with $\eps_n=0$ for all $n$, then by (c) $u\equiv \rho$.

\end{itemize}

\renewcommand{\theequation}{\thesection.\arabic{equation}}
\setcounter{equation}{0}

\chapter{{\bf The basic particle model and its hydrodynamic limit}}
\label{ch:nuovo}
In this chapter we study the hydrodynamic limit of the particle version
of the basic model which has been introduced in
Chapter \ref{ch:1}.  We will prove in this chapter   convergence of the empirical density to the solution of the FBP of Part I, see
Theorem \ref{CDGPthm2.2.1}.  In Section \ref{secnuovo.1}
we recall the definition of the particle system and state the main result. In Section  \ref{secnuovo.2} we outline the strategy of the proof which is then given in the successive sections.

\renewcommand{\theequation}{\thesection.\arabic{equation}}
\setcounter{equation}{0}

 \vskip2cm

\section{The model and the main result}
\label{secnuovo.1}
We fix an initial ``macroscopic profile'' $\rho_0(r)$, $r\in \mathbb R_+$:
we suppose that $\rho_0(r)$ is smooth, has compact support and satisfies
the assumptions in Theorem \ref{CDGPthm2.2.1}, so that the FBP with initial datum
$\rho_0$ has a solution (at least for a positive time interval).

The $N$  particle ``approximation'' of $\rho_0$ consists of a system of
$N$ particles, with their positions, $x_1(0),..,x_N(0)$,
distributed independently with the same law
$\rho_0(r)dr$.  Their dynamics are defined by letting the particles move as independent Brownian motions (with reflections at the origin) till the first time $t_1$ of a Poisson point process on $\mathbb R_+$ of intensity $N$ (for notational simplicity we take here the parameter $j$ of Part I equal to $1$; we
are interpreting the events of the  Poisson point process as times).  At $t_1$
the rightmost particle is moved to the origin.  After  $t_1$ the particles move again as independent Brownian motions(with reflections at the origin)  till the second time  $t_2$ of the Poisson process when
the rightmost particle (at time $t_2^-$) is moved to the origin. The operation is repeated with the same rules
and the process is thus defined for all times (because with probability 1 the Poisson process in a compact has a finite number of events).  We denote by $\und x(t)=(x_1(t),..,x_N(t))$ the particle configuration at time $t$ and
by $P^{(N)}$
the law of  $\{\und x(t), t\ge 0\}$.

We finally
define the ``empirical mass density'' at time $t\ge 0$ as the probability  measure on $\mathbb R_+$ given by
 \begin{equation}
    \label{nuovo.1.1}
\pi_t^{(N)}(dr)= \frac 1N \sum_{i=1}^N \delta_{x_i(t) }(r)dr.
   \end{equation}
Our main result in this chapter is:

   \medskip

\begin{thm}
\label{thmnuovo.1.1}
For any $t\ge 0$ and any $\eps>0$,
 \begin{equation}
    \label{nuovo.1.2}
\lim_{N\to \infty}P^{(N)}\Big[ \sup_{r\ge 0}\Big|\int_r^\infty \pi_t^{(N)}(dr')
- \int_r^\infty S_t\rho_0(r')dr'\Big| > \eps\Big] = 0
   \end{equation}
where $S_t\rho_0$ is defined in Theorem \ref{CDGPthm2.4.1}.
\end{thm}

\medskip

$S_t\rho_0(r)$ coincides with the solution $\rho(r,t)$ of the FBP till when the latter exists, as it follows from Theorems \ref{CDGPthm2.4.1.00} and
item (f) of Theorem \ref{CDGPthm2.3.1}.

\renewcommand{\theequation}{\thesection.\arabic{equation}}
\setcounter{equation}{0}

 \vskip2cm

\section{Strategy of proof}
\label{secnuovo.2}

The proof of Theorem \ref{thmnuovo.1.1} follows the way we proved
Theorem \ref{CDGPthm2.3.1}.  The first step in fact is
to introduce stochastic upper and lower barriers $\und x^{\delta,\pm}(t)$
with the property that for all $t =k\delta, k\in \mathbb N$,
 \begin{equation}
    \label{nuovo.2.1}
    \und x^{\delta,-}(t) \preccurlyeq  \und x(t) \preccurlyeq \und x^{\delta,+}(t)
   \end{equation}
with $P^{(N)}$-probability 1.  The relation $ \preccurlyeq$ is defined as in \eqref{ch2.4.8},
namely two configurations $\und x$ and $\und y$ are ordered, $\und x\preccurlyeq \und y$,
if for any $r \ge 0$,
 \begin{equation}
    \label{nuovo.2.2}
   | \und x\cap  [r,\infty) | \le   | \und y\cap  [r,\infty) |,
        \end{equation}
having regarded $\und x$    and $\und y$ as subsets of $\mathbb R_+$.  \eqref{nuovo.2.2}
can also be stated in terms of the empirical mass densities: calling $\pi(dr)$ and $\pi'(dr)$ the probability measures associated to $\und x$ and $\und y$ via \eqref{nuovo.1.1}, then
    \eqref{nuovo.2.2} can be written as
     \begin{equation*}
 \int_r^\infty  \pi (dr') \le  \int_r^\infty  \pi' (dr').
   \end{equation*}
The definition of the stochastic barriers $\und x^{\delta,\pm}(t)$ is completely analogous
to the definition of the barriers $S^{\delta,\pm}_tu$ and it will be given in Section \ref{secnuovo.3} together with a proof of  \eqref{nuovo.2.1}.

The second step in the proof of Theorem \ref{thmnuovo.1.1} is to relate the stochastic and the deterministic barriers.  Fix $t>0$ and by default in the   sequel  $\delta \in \{2^{-n}t, n \in \mathbb N\}$.  We will prove that for any $\eps>0$,
 \begin{equation}
    \label{nuovo.2.3}
    \lim_{N\to \infty}P^{(N)}\Big[ \sup_{r\ge 0}\Big|
    | \und x^{\delta,\pm}(t)\cap  [r,\infty) |- N\int_r^\infty
       S ^{\delta,\pm}_t\rho_0(r')dr'\Big| > \eps \Big] =0.
        \end{equation}
The proof of \eqref{nuovo.2.3} is not too hard because
the processes
$\und x^{\delta,\pm}(t)$  are essentially independent Brownian motions (with reflections at the origin)
except at a  finite number of times, namely the times $k\delta \le t$.
\eqref{nuovo.2.3} is proved in Section \ref{secnuovo.4}.

The conclusion of the proof   of Theorem \ref{thmnuovo.1.1} is at this point a  three $\eps$ argument as we use \eqref{nuovo.2.1} to relate $\und x(t)$ to $\und x^{\delta,\pm}(t)$,
\eqref{nuovo.2.3} to relate  $\und x^{\delta,\pm}(t)$ to  $S ^{\delta,\pm}_t\rho_0$
and Theorem \ref{CDGPthm2.4.1} to relate  $S ^{\delta,\pm}_t\rho_0$ to  $S_t\rho_0$, the details are given in  Section \ref{secnuovo.5}.

\renewcommand{\theequation}{\thesection.\arabic{equation}}
\setcounter{equation}{0}

 \vskip2cm

\section{The stochastic barriers}
\label{secnuovo.3}

We fix $\delta >0$, $K \in \mathbb N$, and with probability 1 we may and will tacitly  suppose in the sequel that no Poisson event occurs at the times $k\delta$, $k\in \mathbb N$.
We define
the processes $\und x^{\delta,\pm}(t)$, $t\le K\delta$,
 iteratively.   We thus suppose to have defined $\und x^{\delta,\pm}(t)$ for $t\le k\delta$ and want to define it till time $t\le (k+1)\delta$.

 We start from
$\und x^{\delta,-}(t)$.
The particles of $\und x^{\delta,-}(t)$ move as independent Brownian motions
(with reflections at the origin) till time $t_1$ which is the first Poisson event after $k\delta$. At $t_1^+$ a new particle with label $N+1$ is added to $\und x^{\delta,-}(t_1^-)$
and
put at the origin.
The same rule is used at the successive times $t_n\in [k\delta,(k+1)\delta]$ of the Poisson process so that at time
$(k+1)\delta$   we will have a configuration $\und y$ with $N+m$ particles, $m$ the number of Poisson events in $[k\delta,(k+1)\delta]$.
$\und x^{\delta,-}((k+1)\delta)$ is then
obtained from $\und y$ by taking away the
rightmost $m$ particles and relabeling the remaining $N$ with labels $1,..,N$
in some arbitrary way.

The definition of the upper barrier $\und x^{\delta,+}(t)$ requires some more care as it will be defined for each $\delta$ only
in a subset whose probability however goes to 1 as $N\to \infty$.  Such a subset depends  only on the Poisson process: denote by $n_k$ the number of events of the Poisson process in the time
$[k\delta,(k+1)\delta]$;  we will then define $\und x^{\delta,+}(t)$, $t \le K\delta$, on the subset $\{ n_k < N, k=0,..,K-1\}$ observing that for any $K$ and any $\delta \in (0,1)$,
\begin{equation}
    \label{nuovo.3.1}
    \lim_{N\to \infty}P^{(N)}\Big[ \{ n_k < N, k=0,..,K-1\} \Big] = 1.
        \end{equation}
We next restrict to realizations of the Poisson process such that
$\{ n_k < N, k=0,..,K-1\}$, we  suppose iteratively to have defined $\und x^{\delta,+}(t)$ for $t\le k\delta$ and want to define it till time $t\le (k+1)\delta$. We start by taking away
from $\und x^{\delta,+}((k\delta)^-)$ the rightmost $n_k$ particles and let the remaining particles move as independent Brownian motions (with reflections at the origin)
till the first time $s_1$ of the Poisson event
in $[k\delta,(k+1)\delta]$.  At this time we add a new particle at the origin and keep repeating the above procedure till time $ (k+1)\delta$  where we have added $n_k$ particles, namely exactly the same number of particles we had taken away initially, so that
$\und x^{\delta,+}((k+1)\delta)^+)$ has again $N$ particles.

To prove the stochastic inequalities we will use the following notion: two
Brownian motions $x(t)$ and $y(t)$ with reflections at the origin
are {\em coupled increasingly} if:
\begin{itemize}

\item   $y(0) \ge x(0)$.

\item  $x(t)$ and $y(t)$   are independent B-motions (with reflections at the origin)
till the first time $\tau\ge 0$ when they meet.

\item   $x(\cdot)$  is a B-motion (with reflections at the origin) and $y(t)=x(t)$ for  $t\ge \tau $.

\end{itemize}

\noindent
The marginal laws of $x(t)$ and $y(t)$ are the laws of Brownian motions with reflections
at the origin.

\vskip.5cm

\subsection{Stochastic inequalities: lower bound}

We will prove here the first inequality in \eqref{nuovo.2.1} for all $t= k\delta$,
$k\le K$.  We suppose inductively to have proved that for $k\le n$ there is a relabeling of
$\und x^{\delta,-}(n\delta)$ such that
\begin{equation}
    \label{nuovo.3.2}
   x_i^{\delta,-}(n\delta) \le x_i(n\delta).
        \end{equation}
We couple increasingly each pair $x_i^{\delta,-}(t), x_i(t)$, $i=1,..,N$, (each pair being independent of the others)  till the first time $t$ of the
Poisson process in $[n\delta,(n+1)\delta]$.  If  $x_i(t^-)$ is the rightmost particle in $\und x(t^-)$ then  $x_i(t^+)=0$.  We then set
\[
x_i^{\delta,-}(t^+) = 0,\quad x_{N+1}^{\delta,-}(t^+) = x_i^{\delta,-}(t^-).
\]
We repeat this procedure for all Poisson times $t_1,..,t_m$ in $[n\delta,(n+1)\delta]$,
(with probability 1 we are supposing that no Poisson event occurs at the times $k\delta$).
Thus at time $t=((n+1)\delta)^-$
\[
x^{\delta,-}_i (t) \le x_i (t), \quad i =1,..,N.
\]
If there have been $m$ Poisson events in $[n\delta,(n+1)\delta]$ then
$\und x^{\delta,-}  (t^-)$ has   $m$ other particles, $x^{\delta,-}_i (t), i=N+1,..,N+m$.

$\und x^{\delta,-}(t^+)$ is obtained by removing from $\und x^{\delta,-}(t)$ its rightmost
$m$ particles. We do it iteratively. First we remove the
rightmost particle, if its label is $i>N$ we just take it away.  If instead  $i \le N$
we relabel particle $N+1$ as particle $i$, observing that
\[
x^{\delta,-}_i (t^+)= x^{\delta,-}_{N+1} (t^-) \le
x^{\delta,-}_i (t^-)\le x_i (t).
\]
Thus after the first removal the inequalities
$x^{\delta,-}_i (t^+) \le x_i (t)$, $i=1,..,N$ are preserved.  The same rule is used for the successive removals: if the rightmost particle at a step has label $>N$ we just remove it,
if instead it has label $i\le N$ we take the particle with the smallest label $>N$
and relabel it as particle $i$.  In this way we get
\begin{equation}
    \label{nuovo.3.3}
   x_i^{\delta,-}(((n+1)\delta)^+) \le x_i(n+1)\delta),\quad i=1,..,N.
        \end{equation}
Thus by induction \eqref{nuovo.3.2} is proved for all $n\le K$ hence
the first inequality in \eqref{nuovo.2.1}.

\vskip.5cm

\subsection{Stochastic inequalities: upper bound}

We will also use  induction to prove the second inequality in \eqref{nuovo.2.1}.  We thus suppose  to have proved that for $k\le n$ there is a relabeling of
$\und x^{\delta,+}(n\delta)$ so that
\begin{equation}
    \label{nuovo.3.4}
 x_i(n\delta) \le  x_i^{\delta,+}(n\delta), \quad i =1,..,N
        \end{equation}
and want to prove that the inequality remains valid at time $(n+1)\delta$.

$\und y(n\delta):=\und x^{\delta,+}((n\delta)^+)$  is obtained from
$\und x^{\delta,+}((n\delta)^-)$ by taking away
its $m$ rightmost particles, having
called $m$ the number of events in the Poisson process in the time interval
$[n\delta,(n+1)\delta]$.
We  paint in red the particles to be taken away and in blue the others so that
the system at time $(n\delta)^+$ is described by the triple $(\und x(n\delta),
\und y(n\delta), \und \si(n\delta))$, where $\si_i(n\delta) \in \{ R,B\}$, $i=1,..,N$, according to the color of $y_i$.  If $\si_i = R$ the particle $y_i(n\delta)$ is fictitious, it is just put for convenience, the only particles in $ x^{\delta,+}_i((n\delta)^+)$ are the blue ones, i.e.\ those with $\si_i=B$.

We will next define a joint process  $(\und x(t),\und y(t),\und \si(t))$, $t\in (n\delta,(n+1)\delta]$,  with the property that its marginal
$\und x(t)$ has the law of  the true process while  the marginal $y(t)$
once restricted to the blue particles has the law of $\und x^{\delta,+}(t)$.
We will also check that
     \begin{equation}
        \label{nuovo.3.5}
x_k(t) \le y_k(t), \quad 1\le k\le N,\quad t\in (n\delta,(n+1)\delta]
     \end{equation}
and prove that at the final time $(n+1)\delta$ no red particles are left, so that  the second inequality in \eqref{nuovo.2.1} will be proved.

We define the process iteratively and in such a way
that in between clock events each pair $x_i(t), y_i(t)$ is coupled increasingly and independently of the other pairs.  We thus need to check that at the clock events the inequalities are preserved.
Let $t$ be a clock event and suppose by induction that  $x_k(t^-) \le y_k(t^-), k=1,..N$.  Let \begin{itemize}

\item  $i:  x_i(t^-) = \max_k  x_k(t^-)$,

\item  $j:  x_j(t^-) = \max_{\si_k=R}  x_k(t^-)$.

\end{itemize}

\noindent
All colors $\si_k$ and positions $x_k$ and $y_k$ of particles with label  $k$ different from $i$ and $j$ do not  change  at $t$. When   $i=j$ (namely when $\si_i=R$) we set  $x_i(t^+)=y_i(t^+)=0$, $\si_i(t^+)=B$.
Instead when  $i\ne j$ we set
\begin{itemize}

\item  $x_i(t^+) =  y_i(t^+)=0$, $\si_i(t^+)=B$.

\item  $x_j(t^+) = x_j(t^-)$,  $y_j(t^+) = y_i(t^-)$,  $\si_j(t^+) = B$.

\end{itemize}

%
%
%
%

\noindent
We then have:

\begin{itemize}

  \item  The $x$-process is the true one.

 \item  The $y(t)$-process restricted to the blue particles has the same law as $\und x^{\delta,+}(t)$.

\item $x_k(t^+)\le y_k(t^+), k=1,..,N,$  so that the induction property is proved.

\item  The number of reds decreases by 1 at each clock event, so that at the end there are no red left and $\und y(n\delta)= x^{\delta,+}(n\delta)$. By \eqref{nuovo.3.5}
the second inequality in \eqref{nuovo.2.1} is proved.

\end{itemize}

\renewcommand{\theequation}{\thesection.\arabic{equation}}
\setcounter{equation}{0}

 \vskip2cm

\section{Hydrodynamic limit for the stochastic barriers}
\label{secnuovo.4}

In this section we will prove convergence in the limit $N\to \infty$ of the
stochastic barriers to the deterministic ones.  We will use extensively in the proof
the following semi-norms which are ``sort of weak $L^1$ norms''.

\vskip.5cm

\subsection{Semi-norms}
Let $\mathcal I$ be a partition of $\mathbb R_+$ into intervals of length $\ell>0$,
the generic interval $I\in \mathcal I$ being $[n\ell, (n+1)\ell)$. To be specific
 from now on we take $\ell = N^{-\beta}$, $\beta \in (0,1)$, and  write
$\mathcal I_N$ for $\mathcal I$.  Let $\mu$ and $\nu$ be positive, finite
measures on $\mathbb R_+$ with same total mass. We restrict in the sequel to the case where $\mu$ is the counting
measure associated to $\und x^{\delta,\pm}(t)$ and $\nu(dr)= N
S_t^{\delta,\pm} \rho_0(r)dr$, $t=k\delta$.
With this in mind we define for any subset $\mathcal A \subset \mathcal I_N$,
\begin{equation}
    \label{nuovo.4.1}
  \|\mu-\nu\|_{\mathcal A} = \sum_{I\in \mathcal A}\|\mu-\nu\|_I,\quad \|\mu-\nu\|_I=
  \{ \mu(I) - m_I + \nu(I) - m_I\}
        \end{equation}
where for each $I$:
\begin{equation}
    \label{nuovo.4.2}
  m_I= \sup\{ m\in \mathbb Z: m \le \min\big(\mu(I), \nu(I)\big) \}.
        \end{equation}
        Observe that $m_I \ge 0$ and that
        \begin{equation}
    \label{nuovo.4.2.00}
|\mu(I)-\nu(I)| \le \|\mu-\nu\|_I, \quad \mu(I) \le \|\mu-\nu\|_I + \nu(I).
        \end{equation}

We will derive upper bounds for   $\|\mu-\nu\|_{\mathcal I}$ by taking a real number
$m$ in \eqref{nuovo.4.2} which is $\le \mu(I)$ and $\le \nu(I)$ (namely not necessarily the best value $m_I$).  This is used in the proof of the next lemma:
\medskip

\begin{lem}
\label{lemmanuovo.4.0}
Suppose there are a real number $\alpha_I$, a subset $\mathcal A_0$ of $\mathcal I_N$ and $\zeta>0$
such that  $\alpha_I\le \mu(I)$, $\alpha_I \le \nu(I)$ and
\begin{equation}
    \label{nuovo.4.1.1}
 \sum_{I \in \mathcal A_0} \alpha_I \ge N -\zeta.
        \end{equation}
Then
\begin{equation}
    \label{nuovo.4.1.2}
  \|\mu-\nu\|_{\mathcal I_N}  \le 4\zeta +2 | \mathcal A_0|.
        \end{equation}

\end{lem}

\noindent
{\bf Proof.} We have
\[
N = \sum_{I\in \mathcal I_N}\mu(I) \ge N-\zeta + \sum_{I \notin \mathcal A_0}\mu(I),\quad
\sum_{I \notin \mathcal A_0}\mu(I) \le \zeta.
\]
Let $\beta_I$ be the largest integer $m \le \alpha_I$, then, since $m_I\ge 0$ and
$\beta_i \ge \alpha_i-1$,
\[
\sum_{I \in \mathcal I_N} [\mu(I)-m_I] \le \zeta +
\sum_{I \in \mathcal A_0} [\mu(I) - \beta_I] \le \zeta + N - \sum_{I \in \mathcal A_0} ( \alpha_I
-1) \le 2\zeta + | \mathcal A_0|
\]
having used \eqref{nuovo.4.1.1} in the last inequality.
An analogous bound holds for $\nu$, hence \eqref{nuovo.4.1.2}.
\qed

\medskip

We will use in the next subsection the above lemma with $\zeta = N^a$, $a\in (0,1)$.
We next state and prove some other elementary properties of the semi-norms   where
we are fixing $\delta>0$ and $N$, $\mu$ stands for  the counting measure relative to
a configuration $\und x$ with $N$ particles and
$\nu(dr) = N S_t^{\delta,\pm} \rho_0(r)dr$ for some $t=k\delta$.

\medskip

\begin{lem}

\label{lemmanuovo.4.1}
In the above setup there is $c>0$ so that for any $N$,
 \begin{equation}
    \label{nuovo.4.2.0.2}
 m_I \le  \nu(I) \le c N^{1-\beta}; \quad
 \mu(I) \le  c N^{1-\beta} + \|\mu-\nu\|_I.
        \end{equation}

\end{lem}

\noindent
{\bf Proof.} The first inequality holds by definition, the second one because $ \nu(I)= N \int_I S_t^{\delta,\pm} \rho_0(r)dr$  with
$\|S_t^{\delta,\pm} \rho_0\|_{L^\infty}$   bounded for all $t=k\delta$.  The last inequality follows from \eqref{nuovo.4.2.00}.  \qed

\medskip

\begin{lem}

\label{lemmanuovo.4.2}
Let $c$ be as in Lemma \ref{lemmanuovo.4.1} and let $\mu'(I) \le \mu(I)$,
$\nu'(I) \le \nu(I)$, then
\begin{equation}
    \label{nuovo.4.2.0.2.00}
  \|\mu'-\nu'\|_I \le  \|\mu-\nu\|_I + 2c N^{1-\beta}.
        \end{equation}

\end{lem}

\noindent
{\bf Proof.} Calling $m'_I \in [0, m_I]$ the quantity associated to  $\mu'(I)$ and
$\nu'(I)$,
\[
 \mu'(I) - m'_I + \nu'(I) - m'_I \le  \mu(I) - m'_I + \nu(I) - m'_I \le
 \|\mu-\nu\|_I + 2 m_I,
\]
because $m_I \ge m'_I$. 
\qed

\medskip
The next lemma bounds the distribution-distance in terms of the semi-norms and will be used in the proof of Theorem \ref{thmnuovo.1.1}.

\medskip

\begin{lem}

\label{lemmanuovo.4.2.1}
Let $c$ be as in Lemma \ref{lemmanuovo.4.1}. Then for any $r\ge 0$
\begin{equation}
    \label{nuovo.4.2.0.2.1}
  \bigg|\int_r^{\infty} \mu(dr')-\int_r^{\infty} \nu(dr')\bigg|  \le  \|\mu-\nu\|_{\mathcal I_N}
  + 2c N^{1-\beta}.
        \end{equation}

\end{lem}

\noindent
{\bf Proof.} Given $r\ge 0$ let $I_0$ be the interval which contains $r$ and $\mathcal A$ the set of all $I$ to the right of $I_0$. Call $\mu'(dr')= \mathbf 1_{r'\ge r}\mu(dr')$ and $\nu'(dr')= \mathbf 1_{r'\ge r}\nu(dr')$.
 Then
     \begin{equation*}
 |\int_r^{\infty} \mu(dr')-\int_r^{\infty} \nu(dr')|\le \|\mu-\nu\|_{\mathcal A}
 + \|\mu'-\nu'\|_{I_0}.
        \end{equation*}
By \eqref{nuovo.4.2.0.2.00} the right-hand side is bounded by
$\|\mu-\nu\|_{\mathcal I_N}+  2c N^{1-\beta}$.
 \qed

\medskip

In the next lemma $\mu'$  is the counting measure relative to $\und x'$ which is
obtained from $\und x$ by taking away the rightmost $N^* < N$ particles. Analogously
\[
\nu'(dr) = \nu(dr) \mathbf 1_{r\le R_\nu},\quad \int_{R_\nu}^{\infty} \nu(dr) = \delta N.
\]

\medskip

\begin{lem}

\label{lemmanuovo.4.3}
With the above notation
\begin{equation}
    \label{nuovo.4.2.0.3}
  \|\mu'-\nu'\|_{\mathcal I_N} \le  \|\mu-\nu\|_{\mathcal I_N} + 2c N^{1-\beta}
  + |\delta N-N^*|.
        \end{equation}

\end{lem}

\noindent
{\bf Proof.} Call $R_\mu$
the position of the leftmost particle erased from $\und x$ and suppose that $R_\mu<R_\nu$
(the opposite case is similar and its analysis omitted). Call $I_1$
and $I_2$ the intervals of $\mathcal I_N$ which contain $R_\mu$ and, respectively, $R_\nu$.
We call $\mathcal A_1$ the intervals
(of $\mathcal I_N$) to the left of $I_1$, $\mathcal A_3$ those   to the right of $I_2$
and  $\mathcal A_2$ those in between $I_1$ and $I_2$.  Then
\begin{eqnarray*}
 \|\mu'-\nu'\|_{\mathcal I_N} &=&  \|\mu-\nu\|_{\mathcal A_1}
 +  \|\mu'-\nu\|_{I_1} + \sum_{I\in \mathcal A_2}\nu(I) + \nu'(I_2)
 \\&\le&  \|\mu-\nu\|_{\mathcal A_1\cup I_1} + 2c N^{1-\beta}
 + \sum_{I\in \mathcal A_2}\nu(I) + \nu'(I_2).
\end{eqnarray*}
On the other hand
\begin{eqnarray*}
N^* &=& \sum_{I\in \mathcal A_2\cup \mathcal A_3\cup I_2} \{ m_I + [\mu(I) - m_I]\}
+ [\mu(I_1) -\mu'(I_1)], \\
\delta N  &=& \sum_{I\in \mathcal A_3}\{ m_i + [\nu(I) - m_I]\} + [\nu(I_2)-\nu(I'_2)].
\end{eqnarray*}
By taking their difference we get
\begin{equation*}
\sum_{I\in \mathcal A_2}  m_i \le |N^*-\delta N| + \sum_{I\in \mathcal A_3}[\nu(I) - m_I]
+[\nu(I_2)-  m_{I_2}] -\nu(I'_2).
\end{equation*}
Thus
\begin{eqnarray*}
 \|\mu'-\nu'\|_{\mathcal I_N} &\le&  \|\mu-\nu\|_{\mathcal A_1\cup I_1} + 2c N^{1-\beta}
+ \sum_{I\in \mathcal A_2\cup I_2\cup \mathcal A_3}[\nu(I) - m_I] +  |N^*-\delta N|.
\end{eqnarray*}
hence \eqref{nuovo.4.2.0.3}.  \qed

\medskip

We conclude this subsection by bounding $\|\mu_0-\nu_0\|_{\mathcal I_N}$, where
$\mu_0$ is the counting measure associated to the initial configuration $\und x(0)$
with $N$ particles and $\nu_0(dr)=Nu_0(r)dr$.  Let $I$ be an interval which has non-empty intersection with the support of $u_0$.  Since the particles $x_i(0)$ are distributed with law $u_0(r)dr$:
\begin{equation}
    \label{nuovo.4.2.1}
 \mu_0(I)\ge \nu_0(I) - \bigg| \sum_{i=1}^N \Big( \mathbf 1_{x_i(0)\in I} -
 P^{(N)}[x_i(0)\in I]\Big) \bigg|.
        \end{equation}
Since the  $x_i(0)$ are mutually independent,
\begin{equation}
    \label{nuovo.4.2.2}
 \lim_{N\to \infty} P^{(N)}\Big[ \sup_I | \sum_{i=1}^N \Big( \mathbf 1_{x_i(0)\in I} -
 P^{(N)}[x_i(0)\in I]\Big)|\ge N^{\alpha_0}\Big] =0
        \end{equation}
provided that
\begin{equation}
    \label{nuovo.4.2.3}
  \alpha_0 > \frac {1-\beta}{2}.
        \end{equation}
This yields (recalling that $u_0$ has compact support)
           \begin{equation}
    \label{nuovo.4.2.4}
  \|\mu_0-\nu_0\|_{\mathcal I_N} \le cN^{\beta +\alpha_0}
        \end{equation}
and since we want $N^{\beta +\alpha_0} < N$ we need
        \begin{equation}
   \label{nuovo.4.2.5}
  \alpha_0 +\beta  <1,\quad  \frac {1-\beta}{2} <  \alpha_0 < 1-\beta.
        \end{equation}

\vskip.5cm

\subsection{The key estimate}

We fix $\delta>0$ and a positive integer $K$.  We call
$\mu^{\pm}_k$, $k=0,..,K$, the counting measure associated to $\und x^{\delta,\pm}(k\delta)$
and $\nu^{\pm}_k(dr)= S_{k\delta}^{\delta,\pm} \rho_0(r)dr$.  We call
$\mathcal I_N$ the partition $\mathcal I$ when $\ell =N^{-\beta}$.


\medskip

\begin{thm}
\label{thmnuovo.4.1}
There are $\alpha $ and $\beta$ in $(0,1)$ and constants $c_k$ so that
    \begin{equation}
        \label{nuovo.4.2a}
\lim_{N\to \infty}P^{(N)}\Big[\bigcap_{k=0}^{K}\Big\{ \| \mu^{\pm}_k-\nu^{\pm}_k   \|_{\mathcal I_N} \le c_k N^{\alpha} \log^k N\Big] = 1.
     \end{equation}

\end{thm}

\medskip

\noindent
{\bf Proof.}  In the course of the proof we will introduce several parameters.

{ \bf Choice of parameters}.
The main parameters are $\alpha$ and $\beta$:  all $\beta$ small enough will work
 (in particular $\beta< 1/2$) while  $\alpha$
should then be $\alpha > 1 -\beta/3$.   We fix the parameter $\alpha_0 $ in \eqref{nuovo.4.2.3} as
$\alpha_0 = 1/2$.
 Other parameters: $\alpha_2=\alpha_3 > (1-\beta)/2$ and such that
$\alpha> \beta+\alpha_3$.  Finally $\ga = \beta/3$.

 \medskip

\noindent
As the proofs are similar we will only check \eqref{nuovo.4.2a}
for $\mu^{+}_k$ and $\nu^{+}_k$.  Call $P^{(N)}_{\und x^{\delta,+}(k\delta)}$
the law of the  process after time $k\delta$ conditioned on having
$\und x^{\delta,+}(k\delta)$ at time $k\delta$.  We can then write
    \begin{eqnarray}
        \label{nuovo.4.3}
&& P^{(N)}\Big[\bigcap_{k=0}^{K}\Big\{ \| \mu^{+}_k-\nu^{+}_k   \|_{\mathcal I_N} \le c_k N^{\alpha}\log^k N\Big]\nn\\&&\hskip1cm
= E^{(N)}\Big[\prod_{k=0}^{K-1}\mathbf 1_{ \| \mu^{+}_k-\nu^{+}_k   \|_{\mathcal I_N} \le c_k N^{\alpha} \log^k N}\nn\\
&& \times
P^{(N)}_{\und x^{\delta,+}((K-1)\delta)}\Big[ \{ \| \mu^{+}_K-\nu^{+}_K   \|_{\mathcal I_N} \le c_K  N^{\alpha} \log^K N\Big]\Big].
     \end{eqnarray}
We will prove that for any $k\le K-1$ if $\und x^{\delta,+}(k\delta)$
 is such that $\| \mu^{+}_k-\nu^{+}_k   \|_{\mathcal I_N} \le c_k N^{\alpha}\log^k N$ then
    \begin{equation}
        \label{nuovo.4.4}
P^{(N)}_{\und x^{\delta,+}(k\delta)}\Big[ \{ \| \mu^{+}_{k+1}-\nu^{+}_{k+1}   \|_{\mathcal I_N} \le c_{k+1}  N^{\alpha} \log^{k+1} N\}\Big] \ge 1 - \eps_{k,N}
     \end{equation}
where $\lim_{N\to \infty} \eps_{k,N} =0$ for all $k$.  Applied to \eqref{nuovo.4.3} it gives
    \begin{eqnarray}
        \label{nuovo.4.5}
&& \hskip-1cm P^{(N)}\Big[\bigcap_{k=0}^{K}\Big\{ \| \mu^{+}_k-\nu^{+}_k   \|_{\mathcal I_N} \le c_k N^{\alpha}\log^k N\Big\}\Big] \nn\\&&
\ge P^{(N)}\Big[\bigcap_{k=0}^{K-1}\Big\{ \| \mu^{+}_k-\nu^{+}_k   \|_{\mathcal I_N} \le c_k N^{\alpha}\log^k N\Big\}\Big] - \eps_{K,N}
     \end{eqnarray}
and by iteration
    \begin{eqnarray}
        \label{nuovo.4.6}
&& \hskip-1cm P^{(N)}\Big[\bigcap_{k=0}^{K}\Big\{ \| \mu^{+}_k-\nu^{+}_k   \|_{\mathcal I_N} \le c_k N^{\alpha}\log^k N\Big\}\Big] \nn\\&&
\ge P^{(N)}\Big[\Big\{ \| \mu^{+}_0-\nu^{+}_0   \|_{\mathcal I_N} \le c_0 N^{\alpha}\Big] -
\sum_{k=1}^{K-1}\eps_{k,N}.
     \end{eqnarray}
 \eqref{nuovo.4.2a} follows from \eqref{nuovo.4.6} and
\eqref{nuovo.4.2.4} choosing $\alpha > \alpha_0$.

We are thus left with the proof of  \eqref {nuovo.4.4}.
The first operation is the cutting.  Call $N^*$ the number of events of the Poisson process in the interval $[k\delta,(k+1)\delta]$.  Since $N\delta$ is the intensity of the Poisson process given any $\alpha_1 \in (1/2, \alpha)$ there are for any $n>0$ constants $b_n$ so that
    \begin{equation}
        \label{nuovo.4.7}
\lim_{N\to \infty}P^{(N)}\Big[|N^*-\delta N| > N^{\alpha_1}\Big] \le b_n N^{-n}.
     \end{equation}
Then by Lemma \ref{lemmanuovo.4.3} calling $\mu'$ and $\nu'$ the measures
$\mu_k$ and $\nu_k$ after the cutting, we may restrict to the case
\begin{equation}
    \label{nuovo.4.8}
  \|\mu'-\nu'\|_{\mathcal I_N} \le  c'_k N^\alpha \log^k N,
        \end{equation}
provided $\alpha > {1-\beta}$, $\alpha> \alpha_1$ and with $c'_k$ suitable constants.  We start
with $\nu_{k+1}$ and using a gaussian bound, 
	\begin{equation}
    \label{nuovo.4.9}
  \sum_{I \not \subset [0,\log N]} \nu_{k+1}(I) \le e^{-b \log^2 N},\quad b>0.
        \end{equation}
We partition the time interval $[k\delta,\infty)$ into intervals of length $N^{-\beta}$, the partition $\mathcal J$ obtained in this way is $\mathcal I_N$ shifted by $k\delta$.  We denote by $J$ the elements of $\mathcal J$. Let $\ga \in (0,2\beta/3)$ (for the sake of definiteness $\ga = \beta/3$, see the paragraph {\em Choice of parameters} at the beginning of the proof) and $t_\ga$ the  endpoint of the last $J$ in $[0,\delta-N^{-\ga}]$.
Then  for any $I\subset [0,\log N]$,
         \begin{equation}
    \label{nuovo.4.10}
 \nu_{k+1}(I) = \sum_{I'\in \mathcal I_N} A_{I,I'}+ \sum_{J\subset [0,\delta-N^{-\ga}]} B_{I,J} +R_I
        \end{equation}
where
         \begin{equation}
    \label{nuovo.4.11}
  A_{I,I'} =\int_{I}dr \int_{I'}\nu'(dr')G^{\rm neum}_{0,\delta}(r',r),\quad
    B_{I,J} =N\int_{I}dr \int_{J} dt\, G^{\rm neum}_{t,\delta}(0,r),
        \end{equation}
         \begin{equation}
    \label{nuovo.4.12}
 R_I = N\int_{I}dr \int_{t_\ga}^\delta dt G^{\rm neum}_{t,\delta}(0,r).
        \end{equation}
Call $x_I$
and $t_J$ the centers of the intervals $I$ and $J$, then
         \begin{eqnarray}
    \label{nuovo.4.13}
&& |G^{\rm neum}_{0,\delta}(r',r) -  G^{\rm neum}_{0,\delta}(x_{I'},x_I)| \le c N^{-\beta},\; r\in I,\,r'\in I',\nn\\&&
| G^{\rm neum}_{t,\delta}(0,r)- G^{\rm neum}_{t_J,\delta}(0,x_I)| \le c  N^{-\beta + 3\ga/2},\;\: t\in J,\nn
          \end{eqnarray}
        \begin{equation}
    \label{nuovo.4.14}
\sum_I R_I \le c N^{1-\ga},
        \end{equation}
where $c$ is a suitable constant. Denoting by $m'_{I'}$ the left-hand side of \eqref{nuovo.4.2} when $\mu(I) \to \mu'(I')$ and $\nu(I) \to \nu'(I')$, we get
         \begin{eqnarray}
    \label{nuovo.4.15}
        &&| \nu_{k+1}(I) - M_I | \le \Ga+ R_I,\\&&
 M_I= \sum_{I'}m'_{I'}G^{\rm neum}_{0,\delta}(x_{I'},x_I)N^{-\beta}+
\sum_{J\subset [0,\delta-N^{-\ga}]} G^{\rm neum}_{t_J,\delta}(0,x_I)N^{-2\beta}, \nn
\\&& \Ga= c N^{-\beta} \| \mu'-\nu'\|_{\mathcal I_N} + c  N^{1-2\beta}
+c  N^{1-2\beta + 3\ga/2}. \nn
        \end{eqnarray}
For $\mu_{k+1}$ we will only need lower bounds which will be obtained with
similar arguments. The analysis however
will require
probability estimates involving the realization of the Poisson process and the motion of the Brownian particles.  We start from the former.  Call $\und t$ the realizations of the process in $[k\delta,(k+1)\delta]$ then
         \begin{equation}
    \label{nuovo.4.16}
\lim_{N\to \infty}P^{(N)}\Big[\sup_{J\subset [0,\delta-N^{-\ga}]} | |\und t \cap J | - N^{1-\beta}|\le N^{\alpha_2}\Big] = 1,
        \end{equation}
provided $\alpha_2 > (1-\beta)/2$.  We can thus restrict to $\und t$ as in \eqref{nuovo.4.16}.
We thus have $N$ Brownian particles: those in $\und x'$ which start moving at time
$k\delta$ and $N^*$ Brownians which start from the origin  at times $\und t$.  Call $y_i$ the position at time $(k+1)\delta$ of the particle $i$ and given $I$ call $\langle y_i\rangle$ the probability that $y_i$ is in $I$. By the independence of the motion of the particles
we get:
         \begin{equation}
    \label{nuovo.4.17}
\lim_{N\to \infty}P^{(N)}\Big[\sup_{I\subset [0,\log N]} | \sum_{i=1}^N \mathbf
[ 1_{y_i\in I}-\langle y_i\rangle]|\le N^{\alpha_3}\Big] = 1,
        \end{equation}
provided $\alpha_3 > (1-\beta)/2$. We will thus work in the set where \eqref{nuovo.4.16}--\eqref{nuovo.4.17} both hold.
If the label $i$ refers to a particle  t $x'_i$ of $\und x'$ then
\[
\langle y_i\rangle = \int_I G^{\rm neum}_{0,\delta}(x'_i,r)dr.
\]
Analogously, if the label $i$ refers to a particle created at time $t_i$, then
\[
\langle y_i\rangle = \int_I G^{\rm neum}_{t_i,\delta}(0,r)dr.
\]
We use \eqref{nuovo.4.13} and get  a lower bound  
        \begin{eqnarray}
    \label{nuovo.4.18}
        && \mu_{k+1}(I) \ge  M_I - \Delta,
\\&& \Delta= N^{\alpha_3}+c N^{-\beta} \| \mu'-\nu'\|_{\mathcal I_N} + c  N^{1-2\beta}
+c  N^{1-2\beta + 3\ga/2}.  \nn
        \end{eqnarray}
To conclude the proof we use Lemma \ref{lemmanuovo.4.0} choosing
\[
\alpha_I = M_I - \Ga- R_I - \Delta,\quad  \mathcal A_0 =\{I \subset [0,\log N]\}.
\]
We have $| \mathcal A_0| \le N^\beta \log N$.
By \eqref{nuovo.4.15}
\begin{equation*}
 \sum_{I \in \mathcal A_0} \alpha_I \ge   \sum_{I \in \mathcal A_0} \nu_{k+1}(I)
 - \sum_{I } R_I -2(\Ga+ \Delta) N^\beta \log N.
        \end{equation*}
        By \eqref{nuovo.4.9},
\[
 \sum_{I \in \mathcal A_0} \nu_{k+1}(I) \ge N - e^{-b \log N^2}
 \]
 so that using \eqref{nuovo.4.14}
 \begin{equation*}
 \sum_{I \in \mathcal A_0} \alpha_I \ge N - e^{-b \log N^2} -c N^{1-\ga}
  -2(\Ga+ \Delta) N^\beta \log N.
         \end{equation*}
Thus by \eqref{nuovo.4.1.2}
\begin{equation*}
  \|\mu_{k+1}-\nu_{k+1}\|_{\mathcal I_N}  \le  2N^\beta \log N
  + 4\{e^{-b \log N^2} +c N^{1-\ga}
  +2(\Ga+ \Delta) N^\beta \log N \}.
        \end{equation*}
\qed

\renewcommand{\theequation}{\thesection.\arabic{equation}}
\setcounter{equation}{0}

 \vskip2cm

\section{Proof of Theorem \ref{thmnuovo.1.1}}
\label{secnuovo.5}

We fix $t > 0$ and   $\eps>0$ and choose $\delta$ in $\{2^{-n}t, n \in \mathbb N\}$ such
that $\delta \le \eps^2$ and $K:K\delta=t$. As in the previous section we shorthand by
$\mu^{\delta,+}_K$ the counting measure associated to the upper barrier
$\und x^{\delta,+}_{K\delta}$.  By \eqref{nuovo.2.1}
 \begin{equation}
    \label{nuovo.5.1}
\int_r^\infty \pi_t^{(N)}(dr') \le \int_r^\infty N^{-1}\mu^{\delta,+}_K(dr').
   \end{equation}
In the set $\| \mu^{\delta,+}_K-\nu^{\delta,+}_K  \|_{\mathcal I_N} \le c_K N^{\alpha}\log^K N$, where $\nu^{\delta,+}_K(dr)= N S_{K\delta}^{\delta,+} \rho_0(r)dr$, we have by Theorem
 \ref{thmnuovo.4.1},
 \begin{equation}
    \label{nuovo.5.2}
\int_r^\infty N^{-1}\mu^{\delta,+}_K(dr') \le \int_r^\infty S_{K\delta}^{\delta,+} \rho_0(r')dr' +
c_K N^{\alpha-1} \log^K N.
   \end{equation}
 \begin{eqnarray}
    \label{nuovo.5.3}
\int_r^\infty S_{K\delta}^{\delta,+} \rho_0(r')dr'  &\le &
\int_r^\infty S_{K\delta} \rho_0(r')dr' \nn\\&+&
 \{\int_r^\infty S^{\delta,+}_{K\delta} \rho_0(r')dr'-\int_r^\infty S_{K\delta} \rho_0(r')dr'\}.
   \end{eqnarray}
By \eqref{cambi.0} the latter is bounded by
 \begin{eqnarray}
    \label{nuovo.5.4}
 && \hskip-2cm   \int_r^\infty S^{\delta,+}_{K\delta} \rho_0(r')dr'-\int_r^\infty S_{K\delta} \rho_0(r')dr' \nonumber\\
 &\le &
\int_r^\infty S^{\delta,+}_{K\delta} \rho_0(r')dr'-\int_r^\infty S^{\delta,-}_{K\delta} \rho_0(r')dr' \le 2 \eps^2.
   \end{eqnarray}
   By taking $N$ large enough, $c_K N^{\alpha-1} \log^K N \le \eps^2$, so that
   for all $r\ge 0$,
  \begin{equation*}
\int_r^\infty \pi_t^{(N)}(dr') \le
  \int_r^\infty S_t\rho_0(r')dr'+ 2 \eps^2 \le   \int_r^\infty S_t\rho_0(r')dr' + \eps
   \end{equation*}
(for $\eps$ small enough) in the set $\| \mu^{\delta,+}_K-\nu^{\delta,+}_K  \|_{\mathcal I_N} \le c_K N^{\alpha}\log^K N$.  By Theorem  \ref{thmnuovo.4.1} this set has full measure in the limit $N\to \infty$
hence the upper bound in  Theorem \ref{thmnuovo.1.1}. The lower bound is proved in an analogous way.


\part{Variants of the basic model}

\renewcommand{\theequation}{\thesection.\arabic{equation}}
\setcounter{equation}{0}

\chapter{{\bf Introduction to Part II}}
\label{ch:1bbb}

In part I we developed a general approach to study problems
with injection and removal of mass. We showed that such an approach
can be applied to the model in the continuum (using
deterministic mass transport inequalities), as well as to interacting
particle systems (where the inequalities hold point-wise for
almost all random trajectories).  Indeed, it is precisely this
common structure that allowed us in Chapter \ref{ch:nuovo}
to prove that -- in the hydrodynamic limit -- the empirical mass density
of the basic particle model converges to the classical solution
of the free boundary problem defined by \eqref{intro.2} and \eqref{intro.4}.

\vspace{.5cm}
In part II we discuss several problems that can, or possibly could, be
studied using the general approach of part I. We  address the following issues.
\begin{itemize}
\item [i)]
We start by considering a model of particles that move as continuous-time independent random walkers  in the interval $[0,N]\cap\mathbb{Z}$
(with reflecting boundary conditions). In addition, there is injection
of particles at the origin and removal of particles at the rightmost 
occupied site at the event time of two independent Poisson point 
processes, both of intensity $j/N$. It is well know that in the absence
of the injection/removal mechanism the empirical density field converges in the
diffusive scaling limit to the solution of the heat equation on $[0,1]$
with Neumann boundary condition.
We  argue that the scaling limit holds true also with injection/removal
of particles. Namely, in the diffusive scaling the density field of
{\em independent random walkers with current reservoirs}  converges
to the solution of the free boundary problem \eqref{intro.2} and \eqref{intro.4}
now defined in the interval $[0,1]$.
The hydrodynamic limit of this process process was considered
in \cite{CDGP}. We discuss in Chapter \ref{hydro-random-walkers}
the main differences with respect to the spatial setting considered in part I
(where particles could move instead on the half-line).

\vspace{.5cm}
\item [ii)]
Next we address the consequences of 
having two independent Poisson processes ruling 
the injection and removal of mass. Obviously in this case
mass is no longer conserved  at microscopic level. However,
since the intensity of creation and removal of mass is $j/N$,
one needs to go beyond the diffusive scaling to see relevant 
mass fluctuation. We will see that indeed one needs to consider
a {\em super-hydrodynamic limit} (where time is speed-up
by a factor $N^3$ and space is rescaled by a factor $N$) 
to find a meaningful scaling for this second time scale. 
For the case of independent random
walkers with current reservoirs this was considered in \cite{CDGP2}
and it will be discussed in Chapter \ref{chap-super}.

\vspace{.5cm}
\item [iii)]
The last Chapter is devoted to the discussion of several models
with different mechanisms for creation and annihilation of particles.
This includes models with a diffuse injection of mass
(which extends the model with creation of particles at the origin), the
Brunet-Derrida model (which is a model for a population with 
Darwinian selection), as well as the Durrett-Remenik model.
Next we will consider models with two species of particles
whose macroscopic behavior is described by systems
of free boundary problems. Last we will briefly discuss models 
with only mass removal, in which the total particle number decreases to zero. In this context the edge follows a monotonous trajectory and thus there is a better control of
the solution of the corresponding FBP, in particular 
the classical solutions are global in time.

\end{itemize}

\renewcommand{\theequation}{\thesection.\arabic{equation}}
\setcounter{equation}{0}

\chapter{Independent walkers with current reservoirs}
\label{hydro-random-walkers}

In this chapter we consider the model introduced in \cite{CDGP}, 
consisting of independent particle moving as continuous time 
random walkers on a finite lattice, including injection of particles
at the origin and removal from the rightmost occupied site.
We discuss similarities and differences with the setting
deve\-lo\-ped in Part I.

\vskip1cm

\section{Introduction}

The basic problem that was discussed in part I
is rooted in non-equilibrium 
statistical physics. Indeed the derivation of macroscopic 
laws of transport from microscopic models of interacting particles
is a central theme in the mathematical
physics literature. For instance the heat equation
arises as the hydrodynamic limit of a large class 
of models with diffusive behavior.
When the microscopic system is open  there are different possibilities 
to model the interaction with the exterior.
Traditionally the system is coupled to so-called
{\em density reservoirs} that impose
a given density-field at the boundary of a fixed domain.  
As explained in the Introduction of Part I it is of interest
to consider the situation in which one would rather
like to fix a current-field at the boundary.

The idea of {\em current reservoirs} has been 
introduced in a series of recent papers  (see e.g. 
\cite{CDGP,CDGP2,dptv1,dptv2,dpt,dfp, DF}). 
The main difference -- compared to the traditional setting of density reservoirs --
lies in the topological nature of the interaction among particles. 
In systems with density reservoirs the addition/removal mechanism is of a metric and local nature (only particles at boundary sites interact with the reservoirs).  In the setting of current reservoirs the interaction is topological and highly non-local, indeed the determination of the particle to be removed requires  knowledge of the entire configuration.

\medskip
In this chapter we shall investigate another interacting particle model (somewhat similar to the model in Chapter \ref{ch:nuovo}) whose hydrodynamic limit is again related
to the basic problem of part I. The main differences will be the following.
\begin{itemize}
\item The microscopic dynamics of each single particle will be given by a continuous time random walk. A system of independent random walkers is a more detailed description
of the microscopic particle dynamic and thus it better serves the aim of being a physical model for heat conduction. On the other hand this modification will require
an additional diffusive scaling limit, that was not needed for particles moving
as Brownian motion. 
\smallskip
\item Furthermore, to model a finite system, we will restrict the dynamics to 
a finite interval $[0,N]\cap \mathbb{Z}$, with $N$ an integer. 
The creation of particles will always occur at the origin, whereas the removal of particles will be at $N$ if a particle is present there, or in the rightmost occupied site if
the site $N$ is empty. 
\smallskip
\item We will relax the assumption of particle number conservation at microscopic level,
by using two independent exponential clocks for the creation and annihilation
of particles. As a consequence the macroscopic mass will be conserved 
in the diffusive scaling limit, whilst it will fluctuate on a longer time scale,
which will be called the {\em super-hydrodynamic} limit (see Chapter \ref{chap-super}).
\end{itemize}

\noindent
This  model has been named in \cite{CDGP} as {\em independent random walkers
with current reservoir}. 
Calling  $j>0$ the parameter that controls the amount of the imposed 
current, the system evolves according to the 
following simple rules (for a precise definition see 
the following section):
\begin{enumerate} 
\item[i)] 
particles move as independent, symmetric random walks on 
a finite interval of size $N$ with reflections at the boundaries;
\item[ii)]
new particles are created at rate $j/N$ at the left boundary while 
the rightmost particle is killed also at rate $j/N$.
\end{enumerate}

\begin{center}
\begin{figure}[ht]
\label{figura1}
\centering
\includegraphics[angle=0, width=10cm]
{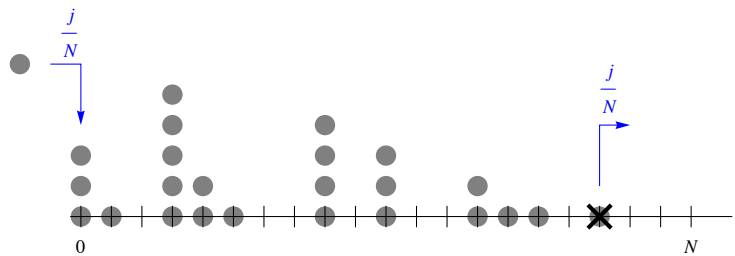}
\caption{ Current reservoirs:   particles are injected at the origin at rate $j/N$ and also the rightmost particle
is removed with the same rate (two independent clocks are used).}
\end{figure}
\end{center}

\noindent
See Figure \ref{figura1} for a pictorial description.

It is worth  observing that the replacement of Brownian particles of Chapter \ref{ch:nuovo} with continuous time random walkers allows us to interpret the system as a queuing model with a spatial structure \cite{ANS,SRI}. Namely, customers enter the queue at the origin following a Poisson process, stay in the queue by changing randomly their position and leave the queue when they reach the rightmost site (being served
at the event time of another Poisson process). The load of the queue
at any given time is thus given by the total particle number, whereas the waiting time 
before being served is related to the location of the rightmost occupied site.

\section{Definition of the model}
\label{definiscimi}

We consider a Markov process $\{\xi_t,t\ge 0\}$ on the space $\Omega$ of particles configurations
$\xi = (\xi(x))_{x\in [0, N]}$, the component $\xi(x)\in\mathbb{N}$ is interpreted as
the number of particles at site $x$. The generator $L$ of the process, working on functions $f: \Omega \to \mathbb{R}$, 
is the sum of three contributions
\footnote{The three terms above have a volume dependence, however 
the  dependence on $N$ is not made explicit.}
 \begin{equation}
\label{generatore}
L   = L^0   +  L^{in}  + L^{out}\;.
\end{equation}
The first term $L^0$ is the generator
of the independent random walks process   
\begin{equation}
\label{0.1-free}
L^0 f(\xi) =  \frac 12 \sum_{x=0}^{N-1} L^0_{x,x+1} f(\xi).
    \end{equation}
    \begin{equation}
\label{generatore-bulk}
L^0_{x,x+1} f(\xi) = \xi(x) \left(f(\xi^{x,x+1}) -f(\xi)\right) + \xi(x+1) \left(f(\xi^{x+1,x}) -f(\xi)\right)
    \end{equation}
where
$\xi^{x,y}$ denotes the
configuration obtained from $\xi$ by removing
one particle from site $x$ and putting it at site $y$, i.e.,
   \begin{equation*}
\xi^{x,y}(z) = \left\{
\begin{array}{ll}
\xi(z)  & \text{if } z \ne x,y,\\
\xi(z) -1  & \text{if } z = x,\\
\xi(z)+1  & \text{if } z=y.
\end{array} \right.
   \end{equation*}
$L^0$ describes
independent symmetric random walks which jump
with equal probability
after an
exponential time of mean 1 to the nearest
neighbour sites, the jumps leading outside $[0,N]$
being suppressed (reflecting boundary conditions).

The term $L^{in}$ in \eqref{generatore} is given by
   \begin{equation}
\label{generatore-birth}
L^{in} f(\xi) = \frac{j}{N} \left(f(\xi^{+}) -f(\xi)\right),\quad \xi^{+}(x) = \xi(x) + \delta_{x,0}\;,
   \end{equation}
where $\delta_{x,y}$ denotes the Kronecker delta.   
For $j>0$,  it describes the action of throwing into the system new particles at rate $\frac{j}{N}$, which then
land at  site  0. Instead $L^{out}$ removes particles and is defined as
    \begin{equation}
\label{generatore-death}
L^{out}f(\xi) =  \frac{j}{N} \left(f(\xi^{-}) -f(\xi)\right),\quad \xi^{-}(x) = \xi(x) - \delta_{x,\bar R_\xi}
    \end{equation}
    where
       \begin{equation}
  \label{3.3.1.00}
 \text{ $\bar{R}_\xi$ is such that:}\;\;\;  \begin{cases}
  \xi (y) >0 & \text{for $y= \bar{R}_\xi$} \\
  \xi (y) =0 & \text{for $y>\bar{R}_\xi$} \;. \end{cases}
   \end{equation}
We also impose $L^{in}f(\xi)=0$ if $\bar{R}_\xi$ does not exist, i.e.\ if $\xi\equiv 0$
is the empty configuration.

\section{Hydrodynamic limit}

The paper \cite{CDGP} proves the existence of the
hydrodynamic limit for independent random walkers with
current reservoirs on a finite macroscopic volume, 
i.e. the existence under diffusive space-time scaling 
of a well-defined function $\rho(r,t)$ 
describing the evolution of an initial profile $\rho_{\rm init}(r)$.
In addition, in \cite{CDGP} it is also proved
that $\rho(r,t)$ is the unique
separating element of suitably defined barriers.

\medskip
\noindent
In this section we recall the main results of  \cite{CDGP}.
While we refer to the original paper for the proofs,
we provide here the main ideas that are used in the proofs.
We shall denote by $P^{(N)}_\xi$ the law of the  process $\{\xi_t,t\ge 0\}$
in the interval $[0,N]$ with generator $L$ given in \eqref{generatore}
and started at time 0 from a configuration  $\xi$.
We consider initial macroscopic profiles $\rho_{\rm init}(r)$ that, similarly to part I, 
belong to the set 
$$
\hat{\mathcal{U}}:= \left\{u\in L^{\infty}([0,1],\mathbb R_+)\cap L^{1}([0,1],\mathbb R_+): \int_0^1 u(r)dr >0\right\}.
$$
The configuration $\xi$ from which the process is started must approximate 
the initial macroscopic profile in the sense of local averages.
That is, the following assumptions are made on the initial particle configuration.
Fix $0< a,b <1$  and denote denote by $\ell$ the integer part of $N^{b}$.
Then we assume that for any $N$ the initial configuration $\xi$ verifies
 \begin{equation}
  \label{3.3.1}
\max_{x\in [0,N-\ell+1]}\Big| \frac 1{\ell} \sum_{y=x}^{ x+\ell-1} \xi(y)-
\frac {N}{\ell} \Big( \int_{x/N}^{x/N+\ell/N}  \rho_{\rm init}(r)dr \Big)\Big| \le \frac{1}{N^a}
\end{equation}
where $\rho_{\rm init}\in\mathcal U$.
Moreover, defining the edge of $\rho_{\rm init}$ as 
	\begin{equation}
	\label{ee8.1bis}
R(\rho_{\rm init})=\inf\{r \in [0,1]: \int_r^1 \rho_{\rm init}(r')dr' = 0\}
	\end{equation}
we also suppose that 
 \begin{equation}
  \label{3.3.1.00.bis}
  \Big| \frac{\bar{R}_\xi}{N}-R(\rho_{\rm init})\Big|\le \frac{1}{N^a}
   \end{equation}
with $\bar{R}_\xi$ the position of the rightmost particle, see \eqref{3.3.1.00}.

The first result in \cite{CDGP} is the following:

\medskip
\begin{theorem}[Existence of hydrodynamic limit, \cite{CDGP}]
 \label{Teo:Hydro}
Let $\rho_{\rm init}\in \hat{\mathcal U}$ and  $\xi$
an approximation in the sense described above.
Then there exists a non-negative, continuous function $\rho(r,t)$
defined on $[0,1]\times \mathbb{R}_+$ 
such that for any $t>0$ and $\zeta>0$,
   \begin{equation}
   \label{eq:Hydro}
\lim_{N \to \infty}P^{(N)}_\xi\Big[\max_{x\in[0,N]} \Big|\frac{1}{N}  \sum_{y=x}^N \xi_{N^{2}t}(y) - \int_{x/N}^1 \rho(r',t) dr'\Big| \le \zeta \Big] = 1   
 \end{equation}
and such that
for any $r\in [0,1]$,
   \begin{equation}
   \label{eq:Hydrobis}
\lim_{t \to 0} \int_r^1 \rho(r',t)dr'  = \int_r^1 \rho_{\rm init}(r')dr'.
    \end{equation}

\end{theorem}

It is easy to see that the above convergence also implies weak convergence of the density field against smooth test functions  $\phi$, i.e. for all $\zeta>0$,
$$
\lim_{N\to \infty} P^{(N)}_\xi\left[\Big|\frac{1}{N} \sum_{x=0}^N \phi\left(\frac{x}{N}\right) \xi_{N^{2}t}(x) - \int_{0}^1 \phi(r)\rho(r,t) dr\Big| \le \zeta\right] = 1\;.
$$
The proof of Theorem \ref{Teo:Hydro} 
(see Figure \ref{Figura222a} for a pictorial representation)
follows closely the path used to prove
the hydrodynamic limit of particle basic problem,
Theorem \ref{thmnuovo.1.1}.
However, in the present setting there are additional
difficulties due to the fact that one needs also
to consider a diffusive scaling for the microscopic
dynamics.
We recall the main steps below, commenting on the main
differences.
\begin{enumerate}
\item The key idea is to define stochastic barriers. These processes, called
$\{ \xi^{(\delta,-)}_{t}, t \ge 0\}$ and $\{ \xi^{(\delta,+)}_{t}, t \ge 0\}$, 
satisfy inequalities with respect to the partial order induced by mass transport and they provide
lower and upper bounds for the original process. 
\smallskip
\item
The proof proceeds by considering discrete time 
intervals  of width $\delta N^{2}$.
 The stochastic barriers $\{ \xi^{(\delta,\pm)}_{k\delta N^{2}}, k \in \mathbb{N}\}$ 
converge weakly as $N\to \infty$ 
to macroscopic objects given, respectively, by the lower barrier $\{\hat{S}_{k\delta}^{(\delta,-)}, k \in \mathbb{N}\}$ 
and the upper barrier $\{\hat{S}_{k\delta}^{(\delta,+)}, k \in \mathbb{N}\}$.
These barriers are defined similarly to those of part I, however they are slightly
different (see below). 
\smallskip
\item
The proof is concluded by observing that in the limit $\delta \to 0$ the upper and lower barriers converge 
to the same limit given by the barriers separating element. Thus also the process $\{\xi_t,t\ge 0\}$,
that is squeezed between the two stochastic barriers, converges
to such limiting object. 
\end{enumerate}
\begin{center}
\begin{figure}[h!]
\begin{equation*}
\begin{array}{lllllc}
\xi^{(\delta,-)}_{k N^{2}\delta} & \preccurlyeq & \xi_{k N^{2}\delta} & \preccurlyeq & \xi^{(\delta,+)}_{k N^{2}\delta} &\\
\\
\Bigg\downarrow & & \Bigg\downarrow & & \Bigg\downarrow & \quad \text{as }\quad N \to \infty \\
\\
\hat{S}^{(\delta,-)}_{k \delta} & \preccurlyeq & \hat{S}_{k\delta} & \preccurlyeq &  \hat{S}^{(\delta,+)}_{k \delta} &  \quad   $$
\end{array}
\end{equation*}
\caption{Scheme of the proof of existence and characterization of the hydrodynamic limit.}
\label{Figura222a}
\end{figure}
\end{center}
The main difference between the barriers $\hat{S}^{\delta,\pm}_{k\delta}$ considered in \cite{CDGP} and the barriers $S^{\delta,\pm}_{k\delta}$ defined in part I are the following:
\begin{itemize}
\item In the definition of the free evolution operator the finite volume 
setting of \cite{CDGP} required to consider the Green function 
for the heat equation in $[0,1]$ with Neumann boundary conditions:
\begin{equation}
\label{e4.3333}
\hat{G}_t^{\rm neum}(r,r') = \sum_{k\in \mathbb Z} G_t(r,r'_k)
\end{equation}
$r'_k$ being the images of $r'$ under
repeated reflections of the interval $[0,1]$ to its right and left, $G_t(r',r)$ as in \eqref{ch2.4.6}.
\smallskip
\item Furthermore, in the definition of the barriers $\hat{S}^{\delta,\pm}_{k\delta}$   
the injection of mass occured at the discrete times. Namely, the free evolution
operator was defined as 
\begin{equation}
\hat{T}_{\delta} u (r) = \hat{G}_t^{\rm neum}* u(r)
\end{equation}
and the cut operator was defined as
\begin{equation}
\hat{C}_{\delta} u (r) = {C}_{\delta} u (r)  + j\delta D_0
\end{equation}
where ${C}_{\delta}$ is the same as Definition \ref{CDGPdefin2.4.1}
and $D_0$ denotes the Dirac delta at zero.
\end{itemize}
Despite these differences, the same analysis of part I could be
carried out. In particular the existence
of a unique separating element of the barriers could
be proved. As a consequence the second result in \cite{CDGP} was 
the following.     
\begin{theorem}[Characterisation of hydrodynamic limit, \cite{CDGP}]
\label{Teo:Char}
Let $\rho_{\rm init}\in \mathcal U$, 
then the hydrodynamic limit $\rho(r,t)$ of Theorem \ref{Teo:Hydro} is
the unique separating element of barriers   $\hat{S}_{n\delta}^{\delta,\pm}(\rho_{\rm init})$,
i.e.,
$$\rho(r,t) = (\hat{S}_t \rho_{\text{init}}) (r)\;.
$$
\end{theorem}
   
\begin{remark}
By the same arguments of part I, we argue that  
the hydrodynamic limit $\rho(r,t)$ is given by 
the solution of the FBP associated to the basic model on the domain $[0,1]$. 
This would essentially be the ``restriction'' of the basic FBP of part I 
(which was defined on the whole half-line $\mathbb{R}_+$)
with the additional constraint that the edge is bounded, i.e. $X_t\le 1$.
\end{remark}

\renewcommand{\theequation}{\thesection.\arabic{equation}}
\setcounter{equation}{0}

\chapter{Beyond diffusive scaling}
\label{chap-super}

In this chapter we analyze the consequences of having
two independent Poisson process for the injection
and removal of mass. We use again the 
model introduced in \cite{CDGP},  for
which we describe the super-hydrodynamic limit.

\vskip1cm

\section{Introduction}
As already remarked in Section \ref{ch:2.2a}, there exist 
stationary classical solutions of the basic FBP on $\mathbb{R}_+$.
Furthermore, we argued at the end of the previous chapter 
that the hydrodynamic limit $\rho(r,t)= \hat{S}_t\rho_{\rm init}(r)$ 
of interacting random walkers with current reservoirs is provided 
by the solution of the basic FBP restricted to the interval $[0,1]$.
It is natural then to conjecture that $\rho(r,t)$ converges as
$t\to\infty$ to the stationary solutions of the basic FBP
on the interval $[0,1]$. However the stationary classical 
solutions of the basic FBP is not unique, 
there exists an entire manifold of stationary linear profiles
labeled by the mass $M$. 
As a result, the following questions naturally arise.

\begin{itemize}
\item What is the basin of attraction (i.e. the set of initial  conditions
that will be attracted to a given stationary solution in the course of time)?
\item Given the existence of infinitely many stationary linear profiles 
of the FBP,  which one will be selected by the microscopic dynamics?
\end{itemize}

\noindent
In this chapter we shall discuss these questions, following the results
obtained in \cite{CDGP2}. We start by describing in section 
\ref{stat-01} the stationary profiles of the basic FBP
on the interval $[0,1]$ and then we describe their basin of attraction.
The second question leads to the identification
of a multi-scale phenomenon whose origin is explained
in section \ref{mass-law} and whose formulation is given in 
section \ref{super-sec}.

\section{Stationary density profiles}
\label{stat-01}

The stationary solutions of the basic FBP restricted to the interval $[0,1]$
are similar to those of the  basic FBP on $\mathbb{R}_+$ already
described in Section \ref{ch:2.2a}. The difference is that, due to the finite
volume, now we need to consider the case of linear profiles truncated at $r=1$, 
i.e. trapezium-shaped profiles.
It is immediate to verify that they are given by
\begin{equation}
\rho_{\rm stat}^{(M)}(r)  =
\begin{cases}
(-2jr + 2\sqrt{Mj}) \mathbf 1_{0\le r \le \sqrt{M/j}} & \qquad\text{if } M \le j,\\
(-2jr + M + j)  \mathbf 1_{0\le r \le 1} & \qquad \text{if } M > j\;.
\end{cases}
\end{equation}
As in Section \ref{ch:2.2a} they are labeled by the value of the
total mass M via the relation
\begin{equation}
\label{e2.6}
\int_0^1 \rho_{\rm stat}^{(M)}(r) dr= M.
\end{equation}
For later convenience we also define $\rho_{\rm stat}^{(0)}\equiv 0$.

The following result is proved in \cite{CDGP2}. It identifies the basin of attraction of
the linear profiles through the analysis of their stability.

\medskip

  \begin{theorem} [Convergence to the stationary profiles, \cite{CDGP2}]
   \label{thme2.3}
 For $r\in[0,1]$, $t > 0$ 
 let $\rho(r,t)= \hat{S}_t\rho_{\rm init}(r)$,  be  the hydrodynamic
limit of the process defined in section \ref{definiscimi}  
with initial configuration  $\xi$  approximating the initial profile $\rho_{\rm init}$ such that
  $\int _0^{1}\rho_{\rm init}(r)dr = M$. Then
   \begin{equation}
   \label{e2.666}
\lim_{t\to \infty}\sup_{r\in [0,1]} 
\Big|\int_{r}^1 \rho(r',t)dr'  -  \int_{r}^1\rho_{\rm stat}^{(M)}(r')dr'\Big|=0\;.
    \end{equation}

 \end{theorem}
%


\section{The law of the total mass}
\label{mass-law}

For the system of independent random walkers with current reservoirs defined
in section \ref{definiscimi}, we consider  the process  $\{|\xi_t| , t \ge 0\}$ yielding 
the particles' number at time $t$, i.e.,
\begin{equation}
|\xi_t| = \sum_{x=0}^{N} \xi_t(x)\;.
\end{equation}
The next theorem shows that this process is very simple, 
despite the complexity of the full process $\{\xi_t , t \ge 0\}$.

\vskip.5cm

\begin{theorem}[Number of particles]
\label{thme3.6a}
The process $\{|\xi_t|, t\ge 0\}$ has the law
of a simple symmetric random walk 
on $\mathbb N$ which jumps with equal probability 
by $\pm 1$ after an {exponential} time of  {parameter} $\frac{2j}{N}$,
the jumps leading to $-1$ being suppressed.
\end{theorem}

\noindent
{\bf Proof.}  From the generator \eqref{generatore} we deduce the generator
of the particle's number process $\{|\xi_t|, t\ge 0\}$  :
\begin{equation}\label{Gen}
{\cal L}f(n) = \frac{j}{N} \Big \{\big( f(n+1)-f(n)\big) + (1-\delta_{|\xi|,0})
\big( f(n-1)-f(n)\big)\Big \}
\end{equation}
where $f$ denotes a bounded function $f: \mathbb N \to \mathbb{R}$.
This coincides with the generator of the simple symmetric random walk
on $\mathbb{N}$ that jumps at rate  $\frac{j}{N}$ and  is reflected at the origin.  
This uniquely characterize the law of $\{|\xi_t|, t\ge 0\}$.
\qed

An immediate consequence of the previous Theorem is the following
\begin{corollary}[Scaling limit]
Let $\rho_{\rm init}\in\hat{\mathcal{U}}$ and  $\xi\in\Omega$
such that
\begin{equation}
M:=\lim_{N\to \infty} \frac{|\xi|}{N} = \int_0^{1} \rho_{\rm {init}}(r) dr.
\end{equation}
Let  $\{\xi_t, t\ge 0\}$  be the processis initialized from $\xi$, 
then the following scaling limits hold:
\begin{equation}
\frac{|\xi_{N^{2}t}|}{N} \to M \qquad\qquad \text{as} \quad N\to \infty ,
\end{equation}
\begin{equation}
\frac{|\xi_{N^{3}t}|}{N} \to B_{jt} \qquad\qquad \text{as} \quad N\to \infty 
\end{equation}
where the converge is in law and
$\{B_t , {t\ge 0}\}$ denotes the Brownian motion on $\mathbb R_+$ 
with reflections at the origin which starts from $B_0 = M$.
\end{corollary}

\section{Super-hydrodynamic limit}  
\label{super-sec}

Hydrodynamics  describes the behavior of the system on times $N^2t$ 
in the limit when $N\to\infty$.  
Hydrodynamics predicts convergence to equilibrium as in Theorem    \ref{thme2.3}.
As a consequence of \eqref{e2.666} we have that 
for any $\zeta>0$,
\begin{equation}
 \label{e2.6666x}
\lim_{t\to \infty} \lim_{N\to \infty}P_\xi^{(N)}
\Big[ \max_{x\in [0,N]}\; 
\Big|\frac{1}{N} \sum_{y=x}^N \xi_{N^2t}(y) - \int_{x/N}^{1} \rho_{\rm{stat}}^{(M)}(r')dr'\Big| {\ge \zeta\Big] =0}
\end{equation}
where
$ M = \int_{0}^1\rho_{\rm init}(r)dr$. \eqref{e2.6666x} shows
convergence of the macroscopic density field to the invariant
profiles.  
Thus, on the hydrodynamic time-scale, the profile that is selected
by the system is dictated by the total mass, which is a conserved 
quantity on the time scale $N^{2}t$.

However, if one inverts the order of the two limits in \eqref{e2.6666x}
then a different result would be obtained.
Indeed, due to the result in the previous section, on a longer time scale 
over which fluctuations of the total mass are allowed,
there is not anymore a privileged profile. 
The investigation of the long time behavior requires  the
study of the process at times  ${N^{2}t_N}$ where $t_N \to \infty$ as $N\to \infty$.
If in this limit we obtain something different from
\eqref{e2.6666x} then we can conclude that there are other significant time-scales
beyond the hydrodynamical one. This has been proved in \cite{CDGP2},
from which we quote the following

\medskip

\begin{theorem}
[Super-hydrodynamic limit, \cite{CDGP2}]
 \label{Teo:Local}
 Let $\xi$ be a sequence such that $\frac{|\xi|}{N}\to M>0$ as $N\to \infty$.  Let $t_N$
 be an increasing, divergent sequence, then the process $\xi_{N^{2}t_N}$ has two regimes:

 \begin{itemize}

 \item Subcritical. If $N t_N\to 0$, then
 \begin{equation}
 \label{e2.6666xx}
\lim_{N\to \infty}P_\xi^{(N)}
\Big[ \max_{x\in [0,N]}\; 
\Big|\frac{1}{N} \sum_{y=x}^N \xi_{N^2t_N}(y) - \int_{x/N}^{1} \rho_{\rm{stat}}^{(M)}(r')dr'\Big| {\le \zeta\Big] =1}.
\end{equation}

\item Critical. Let $t_N=Nt$ then
\begin{equation}
\label{SuperHydro}
\lim_{N\to \infty}P_\xi^{(N)}
\Big[ \max_{x\in [0,N]}\; 
\Big|\frac{1}{N} \sum_{y=x}^N \xi_{N^3t}(y) - \int_{x/N}^{1} \rho_{\rm{stat}}^{(M_t^{(N)})}(r')dr'\Big| {\le \zeta\Big] =1}
\end{equation}
where $M^{(N)}_t:= \frac{|\xi_{N^{3}t}|}{N} $ converges in law as $N\to \infty$ to
  $B_{jt}$, where  {$(B_t)_{t\ge 0}$} is the Brownian motion on $\mathbb R_+$ with reflections at the origin started from $B_0=M$.

 \end{itemize}

\end{theorem}

\medskip

\noindent
We refer to \cite{CDGP2} for the proof of the theorem.
We conclude this section with the following comment.
On a first time scale, i.e.\ the subcritical regime,
the process behaves deterministically and it is
attracted to the invariant linear profile with
mass M (the mass at time zero).  
However on longer times of the order $N^{3}t$
it starts moving stochastically on the manifold of the
linear profiles where it performs a
Brownian motion (with reflection at 0 since the mass can not
become negative). Thus a random behavior arises again
on the super-hydrodynamic time scale.

\renewcommand{\theequation}{\thesection.\arabic{equation}}
\setcounter{equation}{0}

\chapter{Other models}
\label{ch:10}

In this chapter we discuss very briefly other models which have several features in common with our basic model.  The interaction at the particle level has in fact in these models
a topological nature as the rightmost and/or the leftmost particles act differently from the others.  At the macroscopic level this is reflected by a PDE with a free boundary where
the evolution of the edges has to be determined via the outgoing or incoming flux.
The models we present have these features and they can be studied (or have been studied) using barrier inequalities as in Part I.  The strategy is thus the same but the mathematical
problems in its implementation can be quite different.

The models we are going to present have a natural biological motivation.
Particles represent cells, particles positions the states of the cells.  The
natural order in $\mathbb R$ is used to express the fitness of a cell state, for instance the rightmost  cell could be the best fitted in the whole population (of course same analysis would apply when we exchange right and left).  Cells are not clever, they do not know what is best for them and mutate  by exploring all possible nearby states, this is modeled by the cells performing independent Brownian motions.  Cells also duplicate independently of each other this is modeled by adding a new particle say at rate 1 in the same state of the duplicating cell (or in one nearby).  The body which contains the cells cannot support any number of cells, we suppose that a saturation point has been reached for which the number of cells, say $N$, does not change in time.  This means that when a cell duplicates then another cell must be removed from the system.  Here nature imposes its Darwinian law for which the cell removed is the less fitted, the weakest one.  A model with these features has been introduced by Brunet-Derrida and studied by several authors as we will discuss in the sequel.

The question we address and partially answer are: (1) hydrodynamic limit, i.e.\ the FBP associated to the particle model; (2) validity of barrier inequalities; (3) existence and features of stationary states (or traveling waves).

\renewcommand{\theequation}{\thesection.\arabic{equation}}
\setcounter{equation}{0}

 \vskip2cm

\section{Cells evolution in an active environment}
\label{ch:10.1}

Here we consider a variant of the Brunet-Derrida model described above where (1) cells are Brownian particles  in $\mathbb R_+$ (with reflections at the origin),
the cell states (as in the Brunet-Derrida case) are the positions of the particles, but here their fitness decreases when moving to the right, so that $0$ is the best fitted state; (2) it is the environment which creates new cells so that we have an a-priori given probability density $f(r)$, $r\in \mathbb R_+$, with compact support and the state of a new born cell is randomly
distributed with law $f(r)dr$; (3) to preserve, as in the Brunet-Derrida model, the total number  $N$ of cells when a new cell is created the rightmost cell (i.e.\ the weakest, less fitted) is removed; (4) the rate at which a new cell is added is set equal to
$N$ (which corresponds to the rate in  Brunet-Derrida because in that case each particle duplicates at rate 1 so that the intensity for a new particle to appear is $N$).

If $f(r)dr$ is replaced by a delta function at 0 then this becomes the basic model we have studied in Part I (with the parameter $j$ set equal to 1). If instead
\[
f(r)dr = \frac 1N\sum_{i=1}^N \delta_{x_i}(dr)
\]
where $\und x= (x_1,..,x_N)$ is the actual configuration of the cells, then this would be the Brunet-Derrida model (in $\mathbb R_+$).
In the diffuse case (where $f(r)$ is a fixed true function) it may happen that the state of a new born cell is to the right of all the others. In such a case the new cell is also the rightmost one and it is thus removed right away.  This leads to the conjecture that the hydrodynamic limit for this system is ruled by the following FBP:
   \begin{equation}
    \label{10.1.2}
\frac{\partial \rho}{\partial t}= \frac 12 \frac{\partial^2 \rho}{\partial r^2}
+  f,\quad r\in (0,X_t)
   \end{equation}
with an initial datum $\rho_0$, Neumann
boundary condition  at 0
   \begin{equation}
    \label{10.1.3.a}
\frac{\partial \rho(r,t)}{\partial r}\Big|_{r=0}=0
   \end{equation}
while,
at the edge $X_t$, $\rho(X_t,t)=0$ and
   \begin{equation}
    \label{10.1.3}
 -\frac 12\frac{\partial \rho(r,t)}{\partial r}\Big|_{r=X_t}= \phi(X_t),\; \phi(x) := \int_{0}^{x} f(r).
   \end{equation}

%
Namely denoting as in Chapter \ref{ch:nuovo} by $ \pi_t^{(N)}(dr)$ the empirical particles density we conjecture that for any $\eps>0$:
\begin{equation}
    \label{10.1.1}
\lim_{N\to \infty}P^{(N)}\Big[ \sup_{r\ge 0}\Big|\int_r^\infty \pi_t^{(N)}(dr')
- \int_r^\infty  \rho(r',t)dr'\Big| > \eps\Big] = 0
   \end{equation}
where $\rho(r,t)$ is the classical or relaxed solution of \eqref{10.1.2}--\eqref{10.1.3.a}--\eqref{10.1.3}

The proofs of Part I should extend to this case at least when the edge (of the approximating barriers) is to the right of the support of $f$, the analysis of the general case may be more delicate.

Stationary profiles are  analogous  to those of \eqref{ch2.2a.1} (to which they reduce when $f$ is a delta)
\begin{equation}
\label{10.1.5}
\rho^{(st)}(r |M ) = a(M) - 2j \int_0^r \phi(r')dr' \mathbf 1_{a(M)  - 2j \int_0^r \phi(r')dr' \ge 0},\quad \int\rho^{(st)}(r |M ) = M
\end{equation}
and are parameterized by the mass $M$ which is conserved.

\renewcommand{\theequation}{\thesection.\arabic{equation}}
\setcounter{equation}{0}

 \vskip2cm

\section{The Brunet-Derrida evolution-selection mechanism}
\label{ch:10.2}

The Brunet-Derrida model is the one described in the beginning of this chapter.
Namely the cell evolution is described by independent Brownian motions on $\mathbb R$, each cell duplicates independently of the others at rate 1 creating a new Brownian particle
in its same state, simultaneously the rightmost particle is deleted.
The conjectured hydrodynamic limit in this system is
   \begin{equation}
    \label{10.2.1}
\frac{\partial \rho}{\partial t}= \frac 12 \frac{\partial^2 \rho}{\partial r^2}
+ \rho,\quad  {\rm in}\;\; [L_t,+\infty)
   \end{equation}
with initial state $\rho_0(r)$ and boundary conditions at the free boundary $L_t$ given by $\rho(L_t,t)=0$ and
   \begin{equation}
    \label{10.2.2}
\frac{\partial \rho(r,t)}{\partial r}\Big|_{r=0}=0;\quad \frac 12\frac{\partial \rho(r,t)}{\partial r}\Big|_{r=L_t}= M,\; M := \int_{0}^{X_t} \rho(r,t).
   \end{equation}
There is a paper in preparation by A. De Masi,
P. Ferrari, E. Presutti and N. Soprano-Loto which goes in this direction, namely  that
the evolution of the cells in this model is described in the
hydrodynamic limit  by the above FBP, the analysis follows the strategy described in Part I.

%
%

On the whole line there are no longer stationary solutions but there are traveling waves.  One can in fact check that $\rho(r,t) =
\rho^{(tw)}(r-Vt )$ solves  \eqref{10.2.1}--\eqref{10.2.2} where
\begin{equation}
\label{10.2.3}
\rho^{(tw)}(r ) =  M V^2 r e^{-Vr},\quad V^2=2.
\end{equation}
This is not the only traveling wave (with mass $M$) but it is the one with the minimal velocity.  We refer the reader for more details and for related models to the review article  \cite{GJ}.

An interesting (perhaps basic) question is whether there is a stationary, or in this case traveling, measure for the particle system when $N$ is fixed and if its velocity (in the case of traveling waves) is close to the one found by solving the analogous problem for the hydrodynamic equation. This is not at all obvious as the latter
describe the behavior of the system when time is fixed and the number $N$ of particles goes to infinity while we want first to take $t\to \infty$ and then $N\to \infty$.  The analysis of the motion of the system at finite $N$ has been investigated thoroughly by by P. Maillard, \cite{M}, who
has studied the system at times $\log ^3 N$ determining the law of the fluctuations of the edge.

\renewcommand{\theequation}{\thesection.\arabic{equation}}
\setcounter{equation}{0}

 \vskip2cm

\section{The Durrett and Remenik model}
\label{ch:10.4}

Durrett and Remenik, \cite{durrett}, have studied a variant of the Brunet-Derrida model
where cells do not change their states, but they duplicate in a non-local way.
Namely each cell (independently of the others)
creates at rate 1 a new cell in a state $r\in \mathbb R$ which is randomly chosen with probability $\kappa(r-r') dr$, if $r'$ is the state of the generating cell;
$\kappa$ is a smooth probability kernel.  As in Brunet-Derrida simultaneously to the creation of the new cell the leftmost one is erased.

In \cite{durrett} it is shown that for suitable initial data there is a hydrodynamic limit
described by the equation
\begin{equation}
\label{10.4.1}
\frac{\partial}{\partial t} \rho(r,t) = \int \kappa(r'-r) \rho(r',t)
 dr'
\end{equation}
with $\rho(r,0)= \rho_0(r)$, $\int \rho_0=:M$, and
boundary conditions at the left edge $L_t$:
  \begin{equation}
    \label{10.4.2}
\rho(L_t^-,t)=0,\quad  \int_{L_t}^\infty \rho(r,t)  = M.
   \end{equation}
The proof uses barriers in a way similar to that in Part I and traveling fronts are determined.

\renewcommand{\theequation}{\thesection.\arabic{equation}}
\setcounter{equation}{0}

 \vskip2cm

\section{Models with two species}
\label{ch:10.7}

A natural extension of the previous models is when there are two species of cells, say $R$ and $B$
(red and blue). The cells live  both in $\mathbb R$ whose points give their degree of fitness (like in the Brunet-Derrida model).  However  for the red particles fitness increases to the right while for the blue to the left. We suppose that
there
are $N$ red and $N$ blue particles, their number being conserved.
As in the previous models the cells move like
independent Brownian motions however at rate $N$ the weakest red (i.e.\ the leftmost red particle) becomes blue and the weakest blue (i.e.\ the rightmost blue particle) becomes red.

In \cite{DF} for a similar model it has been proved that in the hydrodynamic limit
the system is described by the following FBP:
	\begin{eqnarray}
	\nn
&&
\hskip-.5cm
u_t=\frac 12 u_{rr}+j\,\delta_{V_t},\; r < U_t; \; u(r,0)=u (r), \;u(U_t,t)=0, \; -\frac 12 u_r(U^-_t,t)=j,
\\   \label{eq11}
\\&&
\hskip-.5cm v_t=\frac 12 v_{rr}+j\,\delta_{U_t},\; r > V_t; \; v(r,0)=v (r),\; v(V_t,t)=0, \;   -\frac 12 v_r(V^+_t,t)=-j, \nonumber
	\end{eqnarray}
(under the assumption that this has a classical solution).
\eqref{eq11} is a system of two free boundary equations as the domains
$(-\infty, U_t)$ where $u(r,t)$ is defined
and $(V_t,\infty)$ where $v(r,t)$ is defined are also unknowns to be determined.

\vskip.5cm
So far we have considered models where the particles move independently, Brownian motions or independent random walks.  In the next model there is an interaction between particles.  This is still a two species model but particles are on $\mathbb Z$ with a constraint: at each site there is one particle either red or blue.  Their motion is defined by the ``stirring process'', namely at rate $1/2$
each pair $x,x+1$ of successive points of $\mathbb Z$ exchange their content independently of all the other pairs, so that if at $x,x+1$ we had $R,B$ after the stirring
we have $B,R$.  If instead we had $B,B$ or $R,R$ the stirring does not produce any effect.  Thus the particles are no longer independent, when a particle jumps from $x$ to $x+1$ it forces the opposite jump of another particle.

Allowed configurations are those where there is a rightmost blue and a leftmost red particle.  We may describe the system by giving the positions of only the blue particles (because if at a site there is not a blue particle then there is a red particle), we thus introduce a variable $\eta(x,t)$ equal to 1 when at $x,t$ there is a blue particle and equal to 0 otherwise.
 Allowed
configurations are therefore those where $\eta(x)=1$ definitively as $x \to -\infty$ and $\eta(x)=0$ as $x \to +\infty$. The evolution is determined by the stirring process described earlier
and by a ``selection mechanism'' which here is defined by saying that
at rate $\eps$ the leftmost $0$ becomes 1 and the rightmost $1$ becomes a 0.
In \cite{dfp}
 it is proved that under suitable assumptions on the initial distribution of particles,
when space is scaled as $\eps^{-1}$ and time as   $\eps^{-2}$ the empirical density of 1's
 converges to a limit which is conjectured to satisfy the FBP:
\begin{eqnarray}
    \label{0.1}
&& \frac{\partial \rho}{\partial t}= \frac 12  \frac{\partial^2 \rho}{\partial r^2},  \quad
r\in (L_t,R_t),\\&& L_0, R_0, \rho(r,0) \;\text{given} \nonumber\\&&
\rho(L_t ,t)=1,\;\;\rho(R_t ,t)=0;\quad \frac {\partial\rho}{\partial r}(L_t ,t)= \frac {\partial\rho}{\partial r} (R_t ,t)=-2j. \nonumber
  \end{eqnarray}
This is proved using the same strategy as in Part I, actually \cite{dfp}  is the paper where such a strategy has been introduced.
Thus there exist lower and upper barriers which squeeze in between the actual evolving configuration as described in Part I.  To prove convergence to \eqref{0.1} we would need to reproduce the analysis of Chapter \ref{ch:7}.

%
%
%
%
%
%

\renewcommand{\theequation}{\thesection.\arabic{equation}}
\setcounter{equation}{0}

 \vskip2cm

\section{Interface models}
\label{ch:10.6}

The particle system described at the end of the previous section has also an interpretation in terms of moving interfaces.  The interface is a graph in $\mathbb Z^2$ determined by the particles configuration: if $(x,y)$ belongs to the interface and $\eta(x)=1$ then
$(x+1,y-1)$ also belongs to the interface, while if $\eta(x)=0$ then the next point of the interface is $(x+1,y+1)$.  The correspondence is one to one if we fix for instance the height of the interface at $0$.  The evolution of the particles determines the motion of the interface.

\cite{lacoin}   studies the stochastic evolution   of  interfaces over a ``sticky substrate'', we refer to \cite{lacoin} for the exact definition of the model.  The paper contains a full proof of the hydrodynamic limit for this system, the limit hydrodynamic equation written in terms of the underlying particle system is the FBP
\begin{equation}
\label{10.5.1}
\frac{\partial}{\partial t} \rho(r,t) = \frac 12 \frac{\partial^2}{\partial r^2} \rho(r,t),\quad r\in (L_t,R_t)
\end{equation}
with given initial condition $\rho(r,0)= \rho_0(r)$ and boundary conditions at the free boundaries:
\begin{equation}
\label{10.5.2}
\frac{\partial}{\partial r} \rho(r,t)\Big|_{r = R_t} = -\frac 12
\frac{\partial}{\partial r} \rho(r,t)\Big|_{r = L_t} = \frac 12.
\end{equation}
Global existence of the classical solution of \eqref{10.5.1} (till extinction)
is also proved in \cite{lacoin}.

\eqref{10.5.2} appears also in the analysis of propagation of fire, see for instance Caffarelli-Vazquez, \cite{Caffarelli}. In the $d\ge 1$ setup the FBP in its classical formulation is:
\begin{equation}
\label{10.5.1-new}
\frac{\partial}{\partial t} \rho  = \frac 12 \Delta \rho ,\quad r\in \Om_t
\end{equation}
with  boundary conditions $\rho=0$  and  $\nabla \rho \cdot n = -\frac 12$
on   $\partial\Om_t$, ($n$ the outward normal unit vector to $\Om_t$
at the boundary $\partial\Om_t$).  We refer to the literature for an analysis of this FBP and of other related models.

%
%
%
%
%
%
%
%


%
%
%

\printindex


\end{document}